\documentclass[11pt,a4paper]{article}
\usepackage{authblk}
\usepackage{amsmath}
\usepackage{amsfonts}
\usepackage{amssymb}
\usepackage{amsthm}
\usepackage{xcolor}
\usepackage{placeins}
\usepackage{multirow}
\usepackage{float}
\usepackage{graphicx}
\usepackage{caption}
\usepackage{subcaption}
\usepackage{hyperref}
\usepackage{fullpage}
\usepackage{multirow}
\newtheorem{thm}{Theorem}[section]

\numberwithin{equation}{section}
\theoremstyle{definition}

\newtheorem{rem}[thm]{Remark}

\newcommand{\rev}[1]{\medskip{\color{magenta}#1}}

\title{Reduced basis approximation of parametric eigenvalue problems in presence of clusters and intersections}

\author[K,P]{Daniele Boffi}
\author[K,H]{Abdul Halim}
\author[K,PP]{Gopal Priyadarshi}
\affil[K]{King Abdullah University of Science and Technology, Saudi Arabia}
\affil[P]{Dipartimento di Matematica ``F. Casorati'', University of Pavia, Italy}
\affil[H]{Department of Mathematics, Memari College, West Bengal, India}
\affil[PP]{Department of Mathematics, S.M.D. College, Patliputra University, India}

\begin{document}
\maketitle

\begin{abstract}
In this paper we discuss reduced order models for the approximation of parametric eigenvalue problems. In particular, we are interested in the presence of intersections or clusters of eigenvalues. The singularities originating by these phenomena make it hard a straightforward generalization of well known strategies normally used for standards PDEs. We investigate how the known results extend (or not) to higher order frequencies.
\end{abstract}

\section{Introduction}

Reduced order methods are nowadays a classical tool for the efficient and effective approximation of partial differential equations. The interested reader is referred, for instance, to the monographs~\cite{Quarteronietal16,Hesthavenetal16} and to the references therein.

We aim to investigate the use of reduced order models for the approximation of parametric eigenvalue problems. In the pioneer works~\cite{Madayetal99,Machielsetal00} and~\cite{Prudetal02a,Prudetal02b} the approximation of the first fundamental mode was explored. Investigations for higher (isolated) modes started in~\cite{Pau07a,Pau07b,Pau08}. In~\cite{exploratory} we exploited an idea originating from~\cite{Buchanetal13,GermanRagusa19} of adding a fictitious time variable for the solution of a non parametric eigenvalue problem.

The state of the art of the theory concerning reduced basis approximation of parametric eigenvalue problems stems from the results of~\cite{Fumagallietal16} and~\cite{Horgeretal17}. In particular,~\cite{Fumagallietal16} deals with the first isolated eigenvalue, while~\cite{Horgeretal17} approximates at once a fixed number of eigenmodes, with eigenvalues separated from the rest of the spectrum, starting from the first one.
From those papers it is clear that trouble can come from the intersection of eigenvalues or when eigenvalues are not well separated from each other. Some preliminary discussion about these issues and about the tracking of different modes is contained in~\cite{eccomas,Moataz} and an example of application to the Maxwell eigenvalue problem has been investigated in~\cite{UQ}.

In this paper we are highlighting problematic situations, with the hope that these investigations can open the way towards robust solutions to the presented issues. For this reason, we will discuss rather simple examples where the parametric space is one or two dimensional.

After examining examples where the results of~\cite{Fumagallietal16} and~\cite{Horgeretal17} are confirmed, we investigate the situation when higher frequencies are looked for. For instance, a naive generalization of~\cite{Fumagallietal16} would deal with an isolated eigenvalue which is not the first fundamental one. Analogously, a direct generalization of~\cite{Horgeretal17} involves a cluster of eigenvalues well separated from the rest of the spectrum which does not contain the first fundamental mode. It turns out that the extension of those results to such situations is not true; in particular, the solution computed with a standard POD approach does not satisfy a min-max property with respect to the high fidelity solution. At a first glance we were surprised by the obtained results, even if there are clear explanations for the fact that the conclusions of~\cite{Fumagallietal16} and~\cite{Horgeretal17} do not extend trivially to higher modes.
As a side note, we also try to suggest cheaper solutions for the simultaneous approximation of the first eigenmodes presented in~\cite{Horgeretal17}.
We also show how the obtained solutions can be sensitive to the dimension of the POD reduced basis.

We are confident that the presented results will be useful for future investigations with the aim of optimizing the strategy for the approximation of parametric eigenvalue problems.

In Section~\ref{se:setting} we present the general setting of our elliptic parametric eigenvalue problem, together with its finite element discretization and the standard approach for the construction of a reduced order model based on POD and reduced basis. In Section~\ref{se:overview} we describe the main outcomes of our investigations, which are reported in Section~\ref{se:1D} and~\ref{se:2D} for the one and two dimensional cases, respectively.

\section{Problem setting}
\label{se:setting}
\subsection{Parametric  elliptic eigenvalue problem}
Let $\Omega \subset \mathbb{R}^2 $ be a bounded and polygonal domain and $\mu \in \mathcal M \subset \mathbb{R}^P,$ where $\mathcal M$ is the parameter space. Our goal is to find $(\lambda(\mu), u(\mu))$ such that 
\begin{equation}
\left\{
\aligned
\label{1}
&-\text{div} (A(\mu) \nabla u(\mu)) 
= \lambda(\mu) u(\mu) &&\text{in }\Omega\\
&u(\mu) = 0 &&\text{on }\partial \Omega,
\endaligned
\right.
\end{equation}
where the diffusion matrix $A(\mu) \in \mathbb{R}^{2\times2}$ is symmetric and positive definite for all values of the parameter $\mu$.

It is well known from the theory of compact elliptic eigenvalue problems that the above problem is well-posed and all the eigenvalues are strictly positive.

It is also well known that the regularity of the eigenmodes in terms of the parameter $\mu$ play a crucial role for the study of the solution. In particular, in this paper we are interested in those parametric eigenvalue problems where intersecting eigenvalue phenomena occur.
We have initiated some discussion about this issue in~\cite{eccomas} and we go deeper into it in this work.

\subsection{Finite element approximation}
The finite element approximation of~\eqref{1} is based on the following variational formulation. Let $V$ be the Sobolev space $H_{0}^1(\Omega)$ and $a, m : V \times V \times \mathcal M \to \mathbb{R}$ be the parameter dependent bilinear forms defined by 
\begin{equation}
    \label{eq:a_b_example}
    \begin{aligned}
        a(u,v;\mu)& = \int_\Omega (A(\mu)\nabla u(\mu))\cdot\nabla v\, dx\\
        m(u,v;\mu)& = \int_\Omega u(\mu) v\, dx.
    \end{aligned}
\end{equation}
The weak formulation of~\eqref{1} is defined as: given $\mu \in \mathcal M$, find $(\lambda(\mu),u(\mu))\in \mathbb R_+\times V\setminus{\{0\}}$ such that
\begin{equation}\label{3}
    a(u,v; \mu) = \lambda(\mu) m(u, v; \mu)\quad \forall v\in V.
\end{equation}
We note that $a(\bullet,\bullet;\mu)$ is symmetric and coercive, with coercivity constant depending on the Poincar\'e constant and on $\mu$, and that $m(\bullet,\bullet;\mu)$ coincides with the $L^2(\Omega)$-inner product for all $\mu$. Actually, our setting would allow for a varying bilinear form on the right hand side of~\eqref{3} even if in our examples we take it constant.

Let us consider the finite element subspace $V_{h} \subset V$ of dimension $N_h$
\[
V_{h} = \{v_{h} \in C^{0}(\bar{\Omega}): v_{h}|_K \in \mathbb{P}_{1}(K)\ \forall K \in \mathcal T_{h}\ v_{h}|_{\partial \Omega} = 0 \},
\]
where $\mathcal T_{h}$ is a conforming triangulation of $\Omega$ \rev{with $N_h$ internal vertices} and $\mathbb{P}_{1}(K)$ is the set of polynomials with degree less than or equal to one on $K \in \mathcal T_{h}$.
The finite element or \emph{high fidelity} problem is defined as: given $\mu \in \mathcal M$, find $(\lambda_{h}(\mu),u_{h}(\mu))\in \mathbb R_+\times V_{h}\setminus{\{0\}}$ such that
\begin{equation}\label{4}
    a(u_{h},v_{h}; \mu) = \lambda_{h}(\mu) m(u_{h}, v_{h}; \mu)\quad \forall v_{h}\in V_h.
\end{equation}
Let $\{\phi_i\}_{i=1}^{N_h}$ be a basis of $V_h,$ then the finite element solution $u_h(\mu)$ can be expressed as
\begin{equation*}
    u_{h}(\mu) = \sum_{i = 1}^{N_h} u_{h}^{i}(\mu)\phi_{i}
\end{equation*}
where $\{u_{h}^{i}(\mu)\}_{i = 1}^ {N_h}$ are the finite element coefficients.
Substituting the expression of $u_h(\mu)$ in~\eqref{4}, we obtain the following linear algebraic equation
\begin{equation*}
    A_{h}(\mu)U_{h}(\mu) = \lambda_h(\mu)M_h(\mu)U_h(\mu)
\end{equation*}
where the stiffness matrix $A_{h}$ and mass matrix $M_{h}$ are given by
\[(A_h(\mu))_{i, j} = a(\phi_{j}, \phi_{i}; \mu)\]
\[(M_h(\mu))_{i, j} = m (\phi_{j}, \phi_{i}; \mu)\]
with $i, j = 1, 2, \dots, N_h$ and where $U_{h}(\mu)$ is the column vector formed by the finite element coefficients.

\section{Reduced basis approach}
In the last few years, reduced basis method has emerged as a very powerful method to solve parameterized PDEs. In this method, the original PDE is projected onto a reduced subspace whose basis is obtained from the high fidelity solution evaluated at few suitable sample points.

We start with a set of finite element solutions for various values of the parameter $\mu$, the so called snapshots, and generate a set of $N$ basis functions, the so called reduced basis, $\{\zeta_1, \zeta_2, \dots, \zeta_N\}$. We define the $N$-dimensional reduced basis space as follows
\begin{equation*}
  V_N = \text{span} \{\zeta_i : i = 1, 2, \dots, N\}.
\end{equation*}
There are mainly two approaches based on which reduced basis can be constructed. The first one is the so called greedy approach in which snapshot vectors are chosen based on a posteriori error estimator, whereas the other approach is the proper orthogonal decomposition (POD) in which the reduced basis is chosen based on a singular value decomposition of the snapshot matrix. In our discussion we are going to follow the latter approach.

The reduced basis eigenvalue problem is given by: for $\mu\in\mathcal{M}$, find $(\lambda_{N}(\mu),u_{N}(\mu))\in \mathbb R_+\times V_{N}\setminus{\{0\}}$, such that
\begin{equation}\label{5}
    a(u_{N},v_{N}; \mu) = \lambda_{N}(\mu) m(u_{N}, v_{N}; \mu)\quad\forall v_{N}\in V_{N}.
\end{equation}
Since $\{\zeta_{i}\}_{i=1}^{N}$ form a basis for the reduced space $V_N$, any $u_{N} \in V_{N}$ can be expressed as
\begin{equation*}
    u_{N}(\mu) = \sum_{i = 1}^{N} u_{N}^{i}(\mu)\zeta_{i},
\end{equation*}
where $\{u_{N}^{i}(\mu)\}_{i = 1}^N$ are the reduced basis coefficients.
Substituting the expression of $u_N(\mu)$ in~\eqref{5}, we obtain the following linear algebraic equation 
\begin{equation*}
    A_{N}(\mu)U_{N}(\mu) = \lambda_N(\mu)M_N(\mu)U_N(\mu),
\end{equation*}
where 
\begin{equation}\label{RB_system}
(A_N(\mu))_{i, j} = a(\zeta_{j}, \zeta_{i}; \mu) \mbox{ and } (M_N(\mu))_{i, j} = m (\zeta_{j}, \zeta_{i}; \mu)\quad(i, j = 1, 2, \dots, N)
\end{equation}
and where $U_{N}(\mu)$ is the column vector of the the reduced basis coefficients.

\begin{rem}
\label{re:mult}
In the above discussion it is hidden a crucial difference between the analysis of eigenvalue problems and source problems.
Namely, in a typical source problem there is only one high fidelity solution associated with the source term and any reduced order model is trying to approximate that unique solution.
However, an eigenvalue problem has typically infinitely many solutions: in our case, a sequence of increasing eigenvalues which correspond to (finite dimensional) eigenspaces. The high fidelity problem~\eqref{4} has $N_h$ solutions (counted with their multiplicities) and the classical theory guarantees that, for $h$ small enough, the $k$-th discrete eigenvalue $\lambda_{h,k}(\mu)$ converges to the corresponding continuous one $\lambda_k(\mu)$. Analogously, the discrete eigenfunctions converge to the continous ones, according to a definition that should take into account the multiplicity of the eigenspaces and the fact that convergence involves the entire eigenspace and not just a basis of it (typically the definition of convergence is stated in terms of the gap between subspaces of a Hilbert space).

It turns out, first of all, that the choice of snapshots should take into account which eigenmode(s) we want to approximate. After the choice of the reduced basis, the reduced problem~\eqref{5} has $N$ solutions (counted with their multiplicities) and the main question is whether the $k$-th eigenmode of~\eqref{5} has some similarity with the $\ell$-th eigenmode of~\eqref{4} and, if so, for what values of $k$ and $\ell$.

\end{rem}

\subsection{Relation between the reduced and the high fidelity systems}
Since $V_{N} \subset V_{h}$, the reduced basis $\zeta_i$ can be written  in terms of the finite element basis as
\begin{equation}\label{zeta_exp}
    \zeta_i=\sum\limits_{j=1}^{N_h} \zeta_i^j \phi_j, \quad i=1,\dots, N.
\end{equation}
Let us denote by $\pmb{\zeta}_i=(\zeta_i^1, \dots,\zeta_i^{N_h})^\top\in \mathbb{R}^{N_h}$ the nodal vectors corresponding to the basis $\zeta_i$ and
let us consider the matrix
$$\mathbb{V}=\big[\pmb{\zeta}_1 | \cdots | \pmb{\zeta}_N\big] \in \mathbb{R}^{N_h \times N}.$$
Using the expression \eqref{zeta_exp} we have
\begin{align*}
    a(\zeta_p,\zeta_q;\mu)=\sum\limits_{i=1}^{N_h}\sum\limits_{j=1}^{N_h}\zeta_p^j a(\phi_j,\phi_i;\mu)\zeta_q^i, \quad m(\zeta_p,\zeta_q;\mu)=\sum\limits_{i=1}^{N_h}\sum\limits_{j=1}^{N_h}\zeta_p^j m(\phi_j,\phi_i;\mu)\zeta_q^i,
\end{align*}
 for any $1\leq p,q \leq N$. \\
 
 Thus, we have the following relation between the matrices of the high fidelity and of the reduced systems
 \begin{align*}
     A_N(\mu)=\mathbb{V}^\top A_h(\mu)\mathbb{V},\quad  M_N(\mu)=\mathbb{V}^\top M_h(\mu)\mathbb{V}.
 \end{align*}

 \subsection{Online/offline paradigm}
In order to develop computationally efficient reduced order models, we need some further assumptions on our parametric bilinear forms.
We assume as usual that the bilinear forms $a(\bullet,\bullet ; \mu)$ and $m(\bullet,\bullet; \mu)$ are affine dependent from the parameter, i.e.
\begin{equation*}
    a(\bullet,\bullet; \mu) = \sum_{k= 1}^{S}\theta_{k}(\mu)a_{k}(\bullet,\bullet)
\end{equation*}
and
\begin{equation*}
    m(\bullet,\bullet; \mu) = \sum_{k= 1}^{Q}\Theta_{k}(\mu)m_{k}(\bullet,\bullet),
\end{equation*}
where the bilinear forms $a_{k}(\bullet,\bullet)$ and $m_k(\bullet,\bullet)$ are parameter independent. Hence, these bilinear forms can be assembled only once for all the computations, which leads to a huge computational reduction in the reduced order method.

Let $A^k_h$ and $M^k_h$ be the matrices corresponding to the bilinear forms $a_k(\bullet,\bullet)$ and $m_k(\bullet,\bullet)$, respectively, then the matrix form for the high fidelity problem will be
$$A_h(\mu)= \sum_{k= 1}^{S}\theta_{k}(\mu)A_h^k, \quad M_h(\mu)= \sum_{k= 1}^{Q}\Theta_{k}(\mu)M_h^k$$
and the matrices corresponding to the reduced system will be
$$A_N(\mu)= \sum_{k= 1}^{S}\theta_{k}(\mu)A_N^k=\sum_{k= 1}^{S}\theta_{k}(\mu)\mathbb{V}^\top A_h^k\mathbb{V}, \quad M_h(\mu)= \sum_{k= 1}^{Q}\Theta_{k}(\mu)M_N^k=\sum_{k= 1}^{Q}\Theta_{k}(\mu)\mathbb{V}^\top M_h^k\mathbb{V}.$$
Since the matrices $\mathbb{V}^\top A_h^k\mathbb{V}$ and $\mathbb{V}^\top M_h^k\mathbb{V}$ are parameter independent, these are calculated in the offline stage and in the online stage the reduced matrices are formed by just evaluating the parameter dependent function $\theta_k$ ($k=1,\dots,S$) and $\Theta_k$ ($k=1,\dots,Q$) at the given parameter $\mu$. Hence, in the online stage we only have to solve the reduced system.
 
\subsection{Construction of the POD basis functions}
In this subsection we recall the standard technique for the construction of the reduced basis using POD.

Given a snapshot matrix $S=\big[\pmb{s}_1 | \dots | \pmb{s}_{N_s}\big]$ of $N_s$ snapshots, we want to find a reduced basis $\zeta_1,\dots, \zeta_N$ $(N<< N_h)$ that spans the reduced space $V_N$. The reduced basis can be obtained using the SVD of the snapshot matrix. Applying the SVD to the snapshot matrix $S$, we get
\begin{align}\label{svd_S}
    S=W \Sigma Z^\top
\end{align}
where
$$W=\big[\pmb{\zeta}_1,\dots,\pmb{\zeta}_{N_h}\big] \in \mathbb{R}^{N_h\times N_h} \quad \text{and} \quad {Z}=\big[\pmb{\psi}_1,\cdots,\pmb{\psi}_{N_s}\big] \in \mathbb{R}^{N_s\times N_s}$$
are orthogonal matrices and $\Sigma=\mathrm{diag}(\sigma_1,\dots,\sigma_r,0,\dots,0)\in \mathbb{R}^{N_h \times N_s}$ contains the singular values of $S$ with $\sigma_1\geq \sigma_2\geq \cdots \geq \sigma_r>0$, where $r$ is the rank of the matrix $S$. The columns of $W$ are the left singular vectors of $S$, the columns of the matrix $Z$ are the right singular vectors of $S$, and $\sigma_1,\dots,\sigma_r$ are the singular values of $S$. The first $N$ columns of $W$ are the best choice for the $N$ dimensional basis. The meaning of this claim is that this choice minimizes the sum of squares of the errors between each snapshot vector $\pmb{s}_i$ and its projection onto any $N$-dimensional subspace; this is a consequence of the Schmidt--Eckart--Young theorem which we recall for the readers convenience.

\begin{thm}[Schmidt--Eckart--Young. See~\cite{Quarteronietal16}] \label{SEY}
Given a matrix $S \in \mathbb{R}^{N_h \times N_s}$ of rank $r$, the matrix
$$
S_k=\sum_{i=1}^k \sigma_i \pmb{\zeta}_{i} \pmb{\psi}_i^\top\quad (1\leq k\leq r)
$$
satisfies the optimality property
\begin{equation}
\|{S}-{S}_k\|_F=\min_{\substack{B_{N_h \times N_s}\\\mathrm{rank}({B})\leq k}} \|S-B \|_F =\sqrt{\sum_{i=k+1}^r \sigma_{i}^2},
\end{equation}
where $\|\cdot\|_F$ is the Frobenius matrix norm.
\end{thm}
In general, $N_h \gg N_s$ and the following procedure can be adopted in order to compute the reduced basis, that is the first $N$ left singular vectors.
From \eqref{svd_S} we deduce
\begin{align}
    S Z=W \Sigma \quad \text{and} \quad S^\top W= Z \Sigma^\top
\end{align}
and we can write 
\begin{align}\label{singular_eqn}
    S \pmb{\psi}_i=\sigma_i \pmb{\zeta}_i \quad \text{and} \quad S^\top\pmb{\zeta}_i= \sigma_i \pmb{\psi}_i \quad i=1,\dots,r.
\end{align}
The first relation of \eqref{singular_eqn} gives
\begin{equation}\label{left_sing}
    \pmb{\zeta}_i=\frac{1}{\sigma_i} S \pmb{\psi}_i \quad i=1,\dots,r.
\end{equation}
Using \eqref{left_sing} in the second relation of \eqref{singular_eqn} we get
\begin{align}
     S^\top S\pmb{\psi}_i= \sigma_i^2 \pmb{\psi}_i \quad i=1,\dots,r.
\end{align}
Thus we need to find the first $N$ eigenvectors corresponding to the largest eigenvalues of the symmetric matrix $S^\top S\in \mathbb{R}^{N_s\times N_s}$. Once we find the right eigenvectors $\pmb{\psi}_i$ ($i=1,\dots, N$) we get the first $N$ left eigenvectors using \eqref{left_sing}. The advantage of this procedure is that we need to find eigenvalues of a matrix whose size is $N_s\ll N_h$. The matrix $S^\top S$ is called correlation or Gram matrix.
\begin{rem}
In this paper, we consider snapshot matrices using different eigenvectors and their combinations at the sample parameters (see also Remark~\ref{re:mult}). For example, we consider a snapshot matrix $S$ using the first three eigenvectors at the sample points $\mu_1,\mu_2,\dots,\mu_{n_s}$, that is
$$S=\big[u_{1,h}(\mu_1)|u_{2,h}(\mu_1)|u_{3,h}(\mu_1)|u_{1,h}(\mu_2)|\dots |u_{3,h}(\mu_{n_s})\big].$$
In this case $N_s=3n_s$, and the columns are renamed as $\pmb{s}_i$, $i=1,2,\dots,N_s$.
\end{rem}

From the above discussion it is clear that a crucial aspect will be the selection of the reduced basis dimension $N$.
A viable strategy consists in selecting a tolerance $\epsilon_{tol}$ and in using the following formula.

Choose $N$ as the smallest integer such that
\begin{equation}\label{criterion}
\frac{\sum\limits_{i = 1}^N \sigma_{i}^2}{\sum\limits_{i = 1}^r \sigma_{i}^2} \geq 1 - \epsilon_{tol}
\end{equation}
where $r$ is the rank of the snapshot matrix.

A reasonable choice for our computations seems to be $\epsilon_{tol}$ = $10^{-8}$. However, within the same example, we might vary $N$ in order to see how the results are sensitive to this choice.

\subsection{Parameter sampling technique}
The sampling of the parametric space is very important for the success of reduced order modeling, especially when the dimension of the parameter domain is high and where curse of dimensionality should be tackled. There are several techniques for sampling the parameters such as uniform tensorial grid, Latin hypercube sampling, random sampling, sparse Smolyak sampling with Clenshaw-Curtis points, Monte-Carlo sampling, etc. LHS sampling is a special type of random sampling. In this paper we used uniform sampling as the focus of the paper is to investigate the behavior of reduced order model in relation to the choice of eingenfunctions for the construction of the snapshot matrix.

\begin{figure}
    \centering
     \begin{subfigure}{0.4\textwidth}
        \includegraphics[height=5cm,width=6cm]{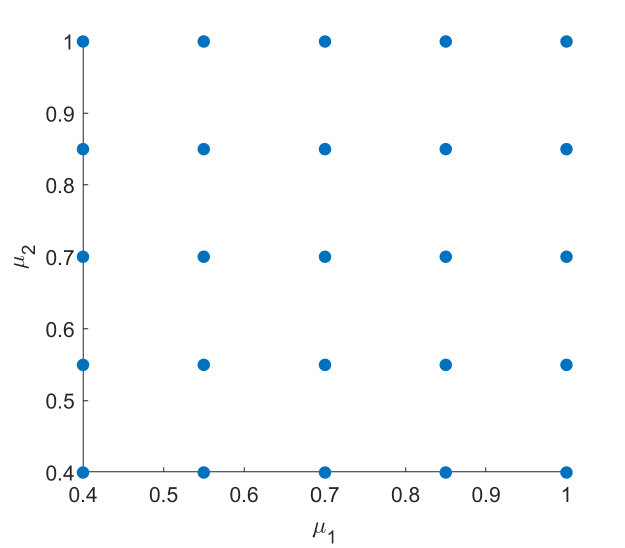}
       \caption{Uniform sampling with 25 points}
     \end{subfigure}
     \begin{subfigure}{0.4\textwidth}
          \includegraphics[height=5cm,width=6cm]{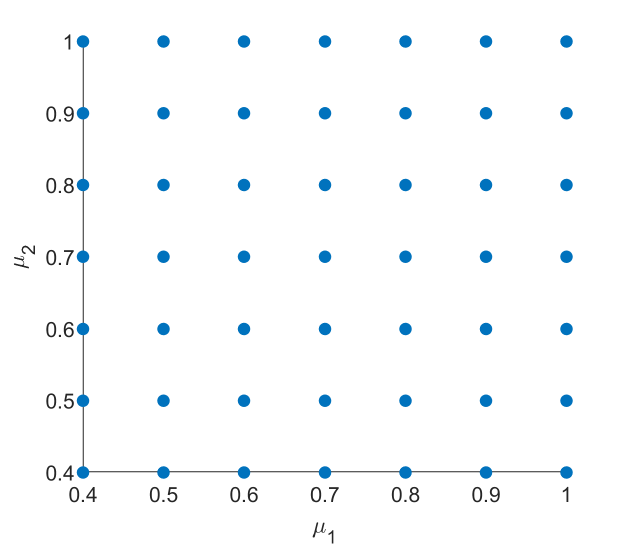}
         \caption{Uniform sampling with 49 points}
     \end{subfigure}
  \caption{Uniform sampling of the parameter domain $[0.4,1]^2$ with 25 and 49 points.}
        \label{sampling1}
  \end{figure}

In Figure~\ref{sampling1} we have displayed the sample points with blue dots using a uniform sampling technique for the parameter domain $\mathcal M=[0.4,1]^2$. Four test points are used for testing the results of our reduced order model (ROM) in Subsection~\ref{se:2D}.

\section{Overview of the numerical investigations}
\label{se:overview}

We want to investigate numerically the behavior of POD based ROM methods for the approximation of parametric eigenvalue problems in which intersecting eigenvalue phenomena occur. Our main questions are the following ones:

\begin{enumerate}

\item How to form the snapshot matrix? Which eigenfunctions should be included in presence of crossings?
\item How to choose the dimension $N$ of the reduced basis?
\item Is it always the case that if $N$ is large (eventually up to the size of the snapshot matrix) then the solutions of the reduced model converge towards the high fidelity solutions?

\end{enumerate}

From the analysis of the existing literature and from our intuition we were expecting the following behavior related to question number 1.

\begin{enumerate}

\item[a)] If we want to approximate an eigenmode with an eigenvalue $\lambda(\mu)$ which has no intersections with other modes in the range of the parameter $\mu\in\mathcal{M}$, then the snapshot matrix could be computed by using only high fidelity  eigenfunctions associated with $\lambda(\mu)$.
\item[b)] If we want to approximate an eigenmode with an eigenvalue that has intersections with other eigenvalues and if all these eigenvalues are well separated from the rest of the spectrum in the range of the parameter $\mu\in\mathcal{M}$, then the snapshot matrix should be computed by considering high fidelity eigenfunctions associated with all the modes that cross each other.

\end{enumerate}

Concerning question number 2 we thought that the heuristic strategy described with formula~\eqref{criterion} could be a good starting point for our numerical tests.
Consequently, we were expecting the following behavior in connection with question number 3 depending on whether we are in presence of cases a) or b) above.

\begin{enumerate}

\item[a)] The first eigenvalue of the reduced system converges toward the isolated high fidelity eigenvalue $\lambda(\mu)$ as $N$ increases. The same applies to the corresponding eigenspace.

\item[b)] If there are $k$ eigenmodes belonging to the cluster that we are considereing, then the first $k$ eigenvalues of the reduced system converge towards the eigenvalues of interest of the high fidelity problem. The same applies to the corresponding eigenspaces.

\end{enumerate}

The results of our computation are in agreement with our expectations \emph{provided we are dealing with the first eigenvalues in the spectrum}. This means that things work well in case a) if $\lambda(\mu)$ is the first fundamental mode for all values of $\mu\in\mathcal{M}$. We will see that things may go wrong if we form the snapshot matrix of an isolated eigenmode which is not the first one.
Analogously, in presence of clusters of $k$ eigenvalues, it is in general a good strategy to form the snapshot matrix by considering all the corresponding eigenfunctions if we are dealing with the first $k$ eigenmodes (with possible crossings). On the other hand, things may go wrong if we perform the same strategy with clusters in the higher part of the spectrum.

We point out that this behavior is in perfect agreement with the theoretical results of~\cite{Fumagallietal16} (first isolated eigenvalue) and~\cite{Horgeretal17} (first $k$ eigenvalues isolated from the next ones). However, our computations show that these results cannot be generalized to higher order modes.

To deal with these  problems, we consider all $n$ eigenfunctions  (corresponding to the first smallest $n$ eigenvalues) simultaneously in the snapshot matrix in order to obtain the first smallest $n$ eigenvalues. It is shown through various examples that using this strategy the reduced order method provides all $n$ eigenvalues (whether intersecting or non-intersecting) and the corresponding eigenfunctions correctly. It is also shown that as we increase the ROM dimension, eigenvalue based on ROM converge to the eigenvalue obtained through finite element approximation. 

If we are interested in finding the first $n$ eigenvalues, then putting all the $n$ eigenvectors in the snapshot matrix at the sample points, we get good results which also maintain the order of the eigenvalues but the number of snapshots increased by multiple of $n$. In order to reduce the size of the snapshot matrix, we consider the following \textbf{alternate} strategy for getting the first $n$ eigenvalues simultaneously. Instead of taking all the $n$ eigenvectors, we choose some linear combination
\begin{equation}\label{com_snap}
\tilde{u}(\mu)=c_1u_1(\mu)+\cdots+c_nu_n(\mu)\quad (c_j\neq 0\ \forall j)
\end{equation}
of them in the snapshot matrix. Our tests show that this can be a cheaper alternative in order to compute $n$ simultaneous approximations of the eigenvalues and of the corresponding eigenfunctions in the correct order. In our numerical tests we use the sum of the first $n$ eigenfunctions in the snapshot matrix.

We present two sets of numerical examples: the first one is for a one dimensional parameter space and the second one is for a two dimensional parameter space.


\section{Numerical results for eigenvalue problems depending on one parameter}
\label{se:1D}

We consider the following eigenvalue problem: for $\mu\in(-\sqrt{2},\sqrt{2})$, find $(\lambda(\mu),u(\mu))\in\mathbb{R}^+\times V\setminus{\{0\}}$ such that
\begin{equation}
\left\{
\aligned
\label{2}
&-\text{div} (A(\mu) \nabla u(\mu)) 
= \lambda(\mu) u(\mu)&&\text{in }\Omega = (-1, 1)^2\\
&u(\mu) = 0&&\text{on }\partial \Omega,
\endaligned
\right.
\end{equation}
where the diffusion matrix $A(\mu) \in \mathbb{R}^{2\times2}$ is given by
$$A(\mu) = 
\begin{bmatrix}
1 & \mu\\
\mu & 2
\end{bmatrix}.
$$
It is easy to check that $A$ is symmetric and its eigenvalues are strictly greater than zero when
$-\sqrt{2}< \mu < \sqrt{2}$.
\begin{figure}[H]
\centering
\includegraphics[scale=0.4]{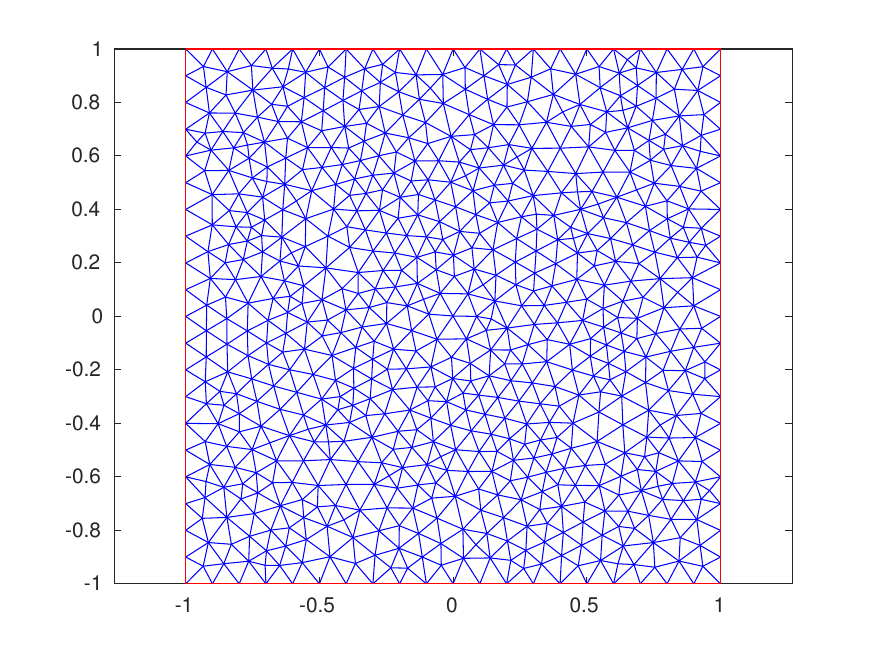}
\includegraphics[scale=0.4]{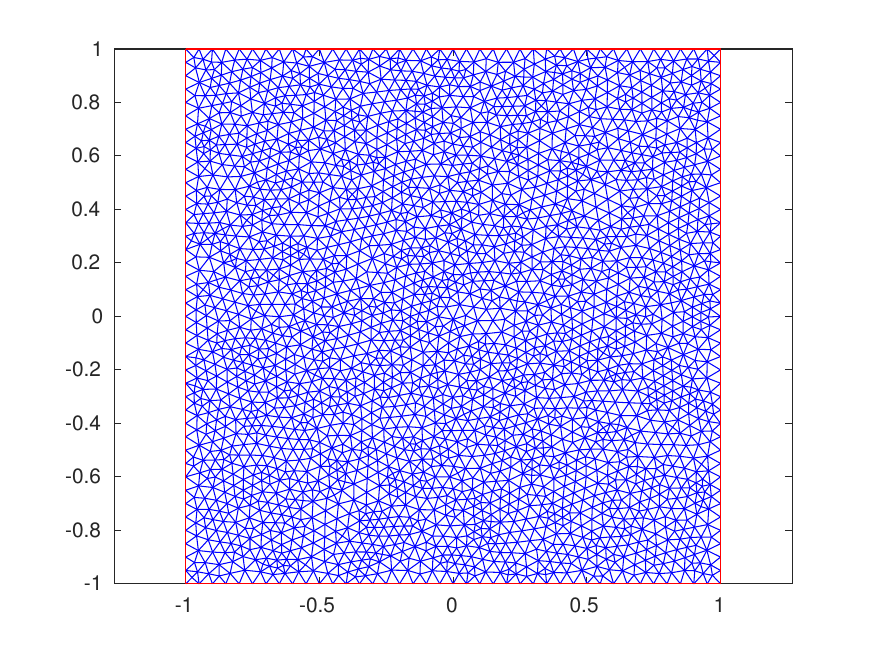}
\caption{Sequence of unstructured meshes ($h = 0.1$ and $0.05$)}
\label{fg:meshes}
\end{figure}
We use the unstructured meshes shown in Figure~\ref{fg:meshes} for finite element approximation, giving the high fidelity solutions reported in Figure~\ref{fig3}.
We have plotted the first six eigenvalues with a legenda and a color code that refers to a local sorting for each value of the parameter. It is clear that the eigenvalues cross each others and that the eigenspaces may have jump discontinuities when a crossing occurs. This phenomenon can be seen, for instance, by looking at the eigenfunctions reported in Figures~\ref{fig:1par_nve} and~\ref{fig:3rdev_1par_pve}.

\begin{figure}
\centering
\includegraphics[scale=0.4]{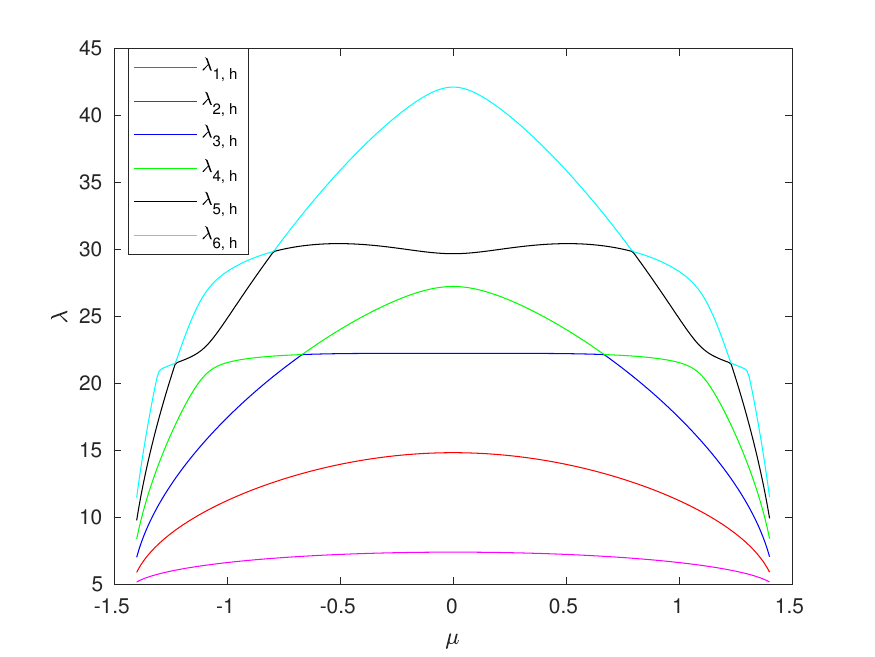}
\caption{First six sorted  eigenvalues when $h = 0.05$ and $\mu = -1.4:0.01:1.4$.} 
\label{fig3}
\end{figure}
\begin{figure}
     \begin{subfigure}{0.32\textwidth}
         \includegraphics[height=4.5cm,width=5.5cm]{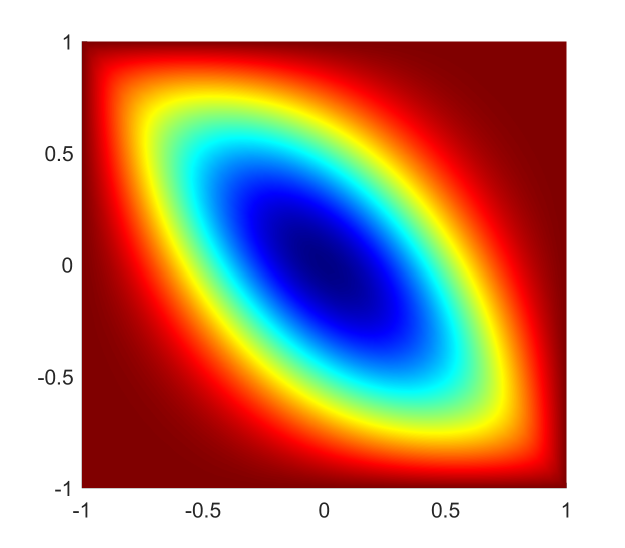}
         \caption{First EV ($\mu=-1.25$)}
     \end{subfigure}
     \begin{subfigure}{0.32\textwidth}
         \includegraphics[height=4.5cm,width=5.5cm]{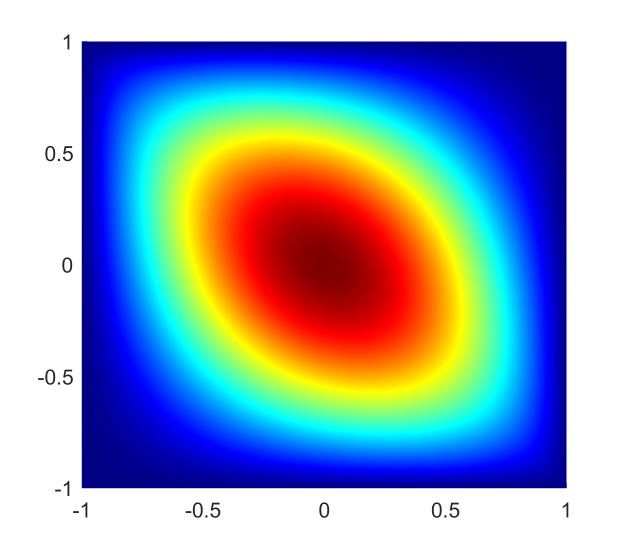}
         \caption{First EV ($\mu=-0.75$)}
     \end{subfigure}
     \begin{subfigure}{0.32\textwidth}
         \includegraphics[height=4.5cm,width=5.5cm]{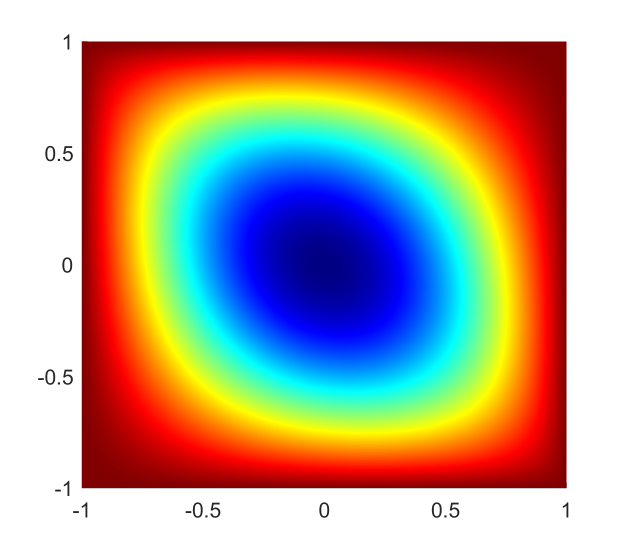}
          \caption{First EV ($\mu=-0.5$)}
     \end{subfigure}
     \begin{subfigure}{0.32\textwidth}
         \includegraphics[height=4.5cm,width=5.5cm]{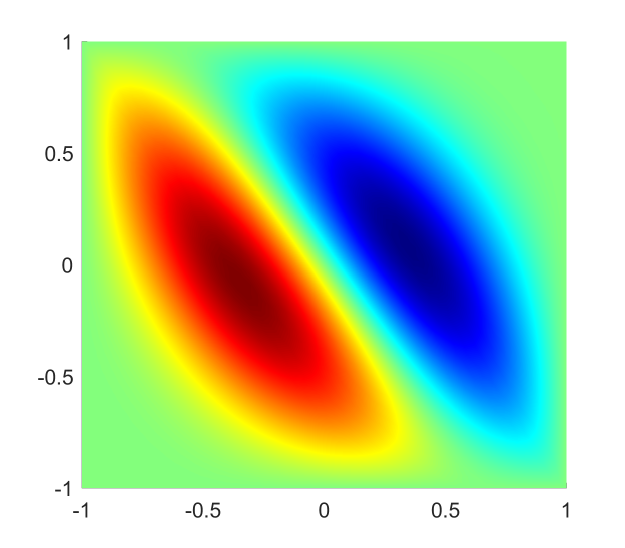}
         \caption{Second EV ($\mu=-1.25$)}
     \end{subfigure}
     \begin{subfigure}{0.32\textwidth}
         \includegraphics[height=4.5cm,width=5.5cm]{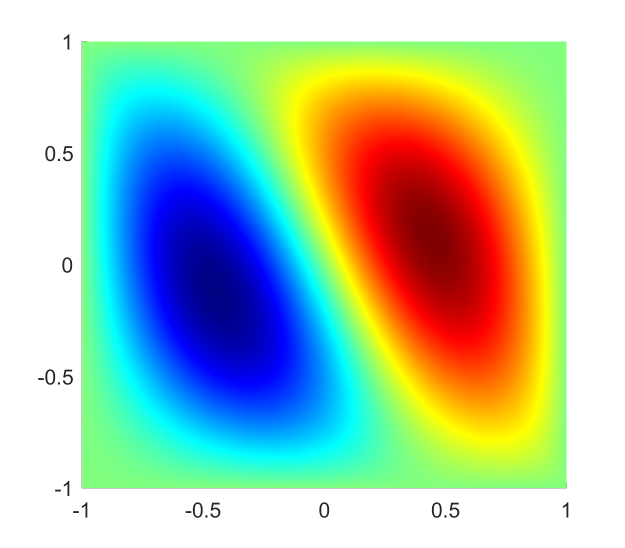}
         \caption{Second EV ($\mu=-0.75$)}
     \end{subfigure}
     \begin{subfigure}{0.32\textwidth}
         \includegraphics[height=4.5cm,width=5.5cm]{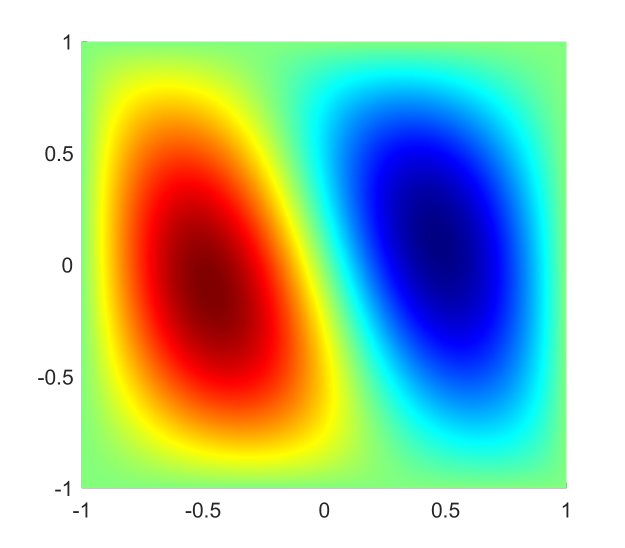}
          \caption{Second EV ($\mu=-0.5$)}
     \end{subfigure}
     \begin{subfigure}{0.32\textwidth}
         \includegraphics[height=4.5cm,width=5.5cm]{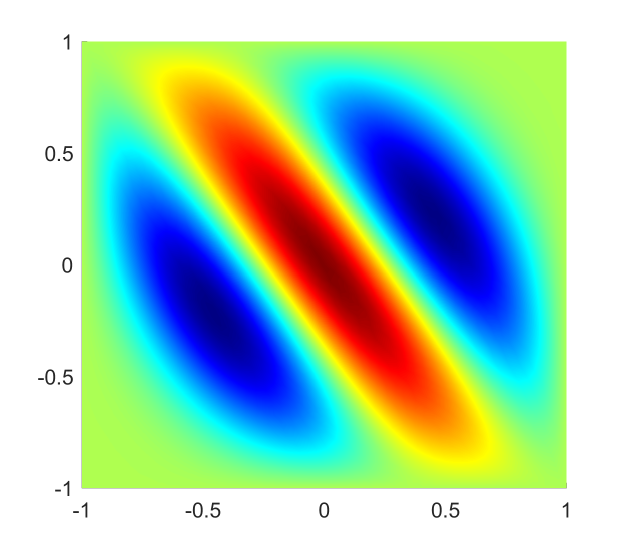}
         \caption{Third EV ($\mu=-1.25$)}
     \end{subfigure}
     \begin{subfigure}{0.32\textwidth}
         \includegraphics[height=4.5cm,width=5.5cm]{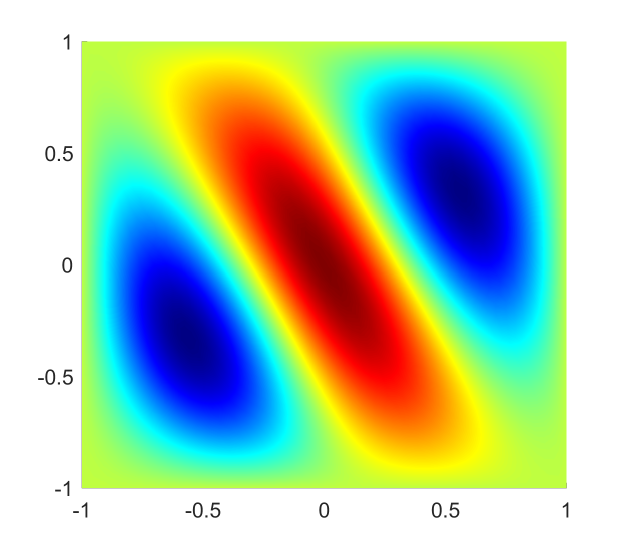}
         \caption{Third EV ($\mu=-0.75$)}
     \end{subfigure}
     \begin{subfigure}{0.32\textwidth}
         \includegraphics[height=4.5cm,width=5.5cm]{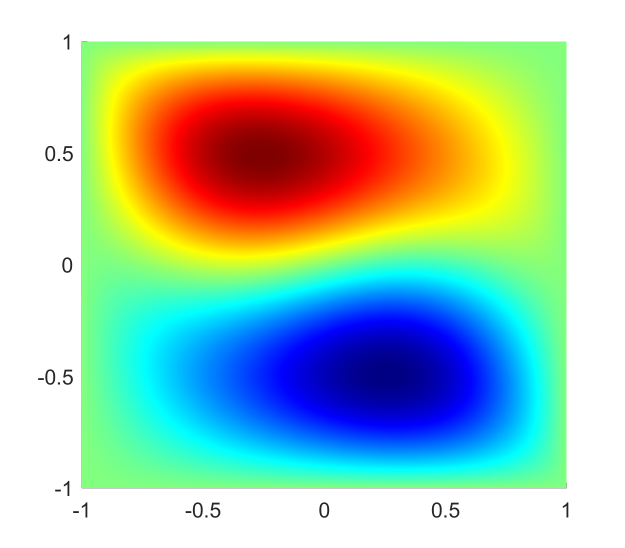}
          \caption{Third EV ($\mu=-0.5$)}
     \end{subfigure}
     \begin{subfigure}{0.32\textwidth}
         \includegraphics[height=4.5cm,width=5.5cm]{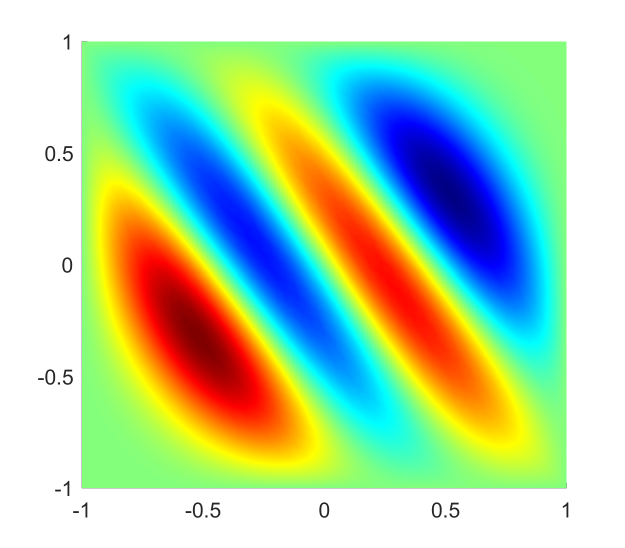}
         \caption{Fourth EV ($\mu=-1.25$)}
     \end{subfigure}
     \begin{subfigure}{0.32\textwidth}
         \includegraphics[height=4.5cm,width=5.5cm]{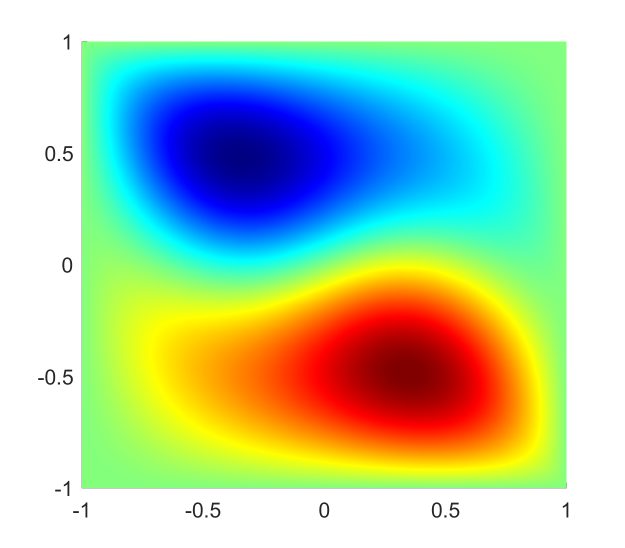}
         \caption{Fourth EV ($\mu=-0.75$)}
     \end{subfigure}
     \begin{subfigure}{0.32\textwidth}
         \includegraphics[height=4.5cm,width=5.5cm]{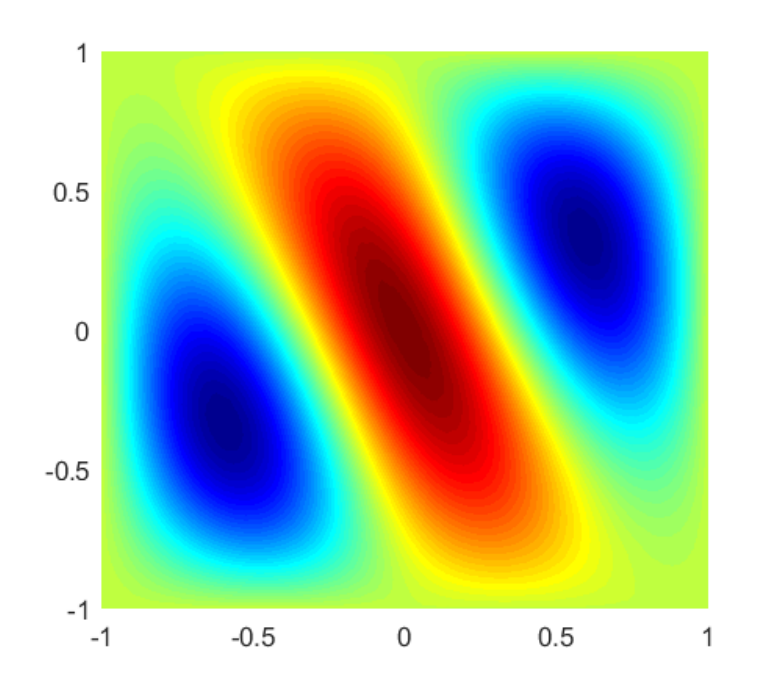}
          \caption{Fourth EV ($\mu=-0.5$)}
     \end{subfigure}
        \caption{First four eigenvectors at $\mu=-1.25,-0.75,-0.5$ using FEM.}
        \label{fig:1par_nve}
\end{figure}

\begin{figure}
     \begin{subfigure}{0.32\textwidth}
         \includegraphics[height=4.5cm,width=5.5cm]{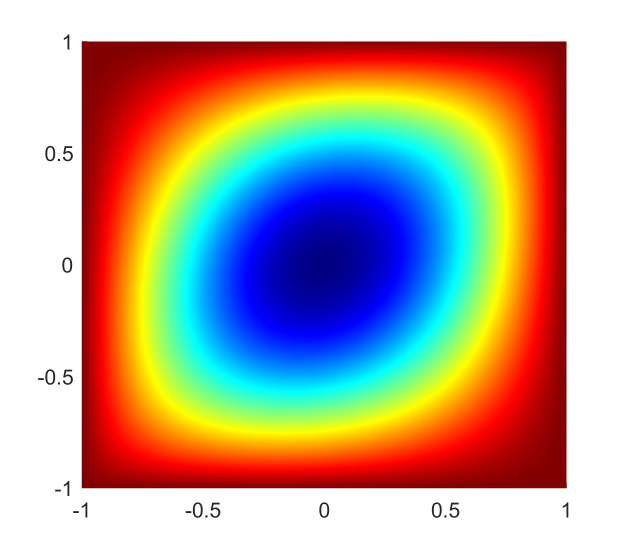}
         \caption{First EV ($\mu=0.5$)}
     \end{subfigure}
     \begin{subfigure}{0.32\textwidth}
         \includegraphics[height=4.5cm,width=5.5cm]{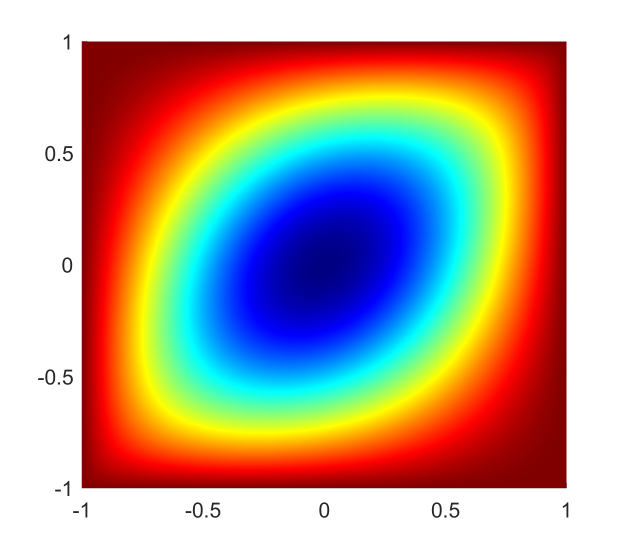}
         \caption{First EV ($\mu=0.75$)}
     \end{subfigure}
     \begin{subfigure}{0.32\textwidth}
         \includegraphics[height=4.5cm,width=5.5cm]{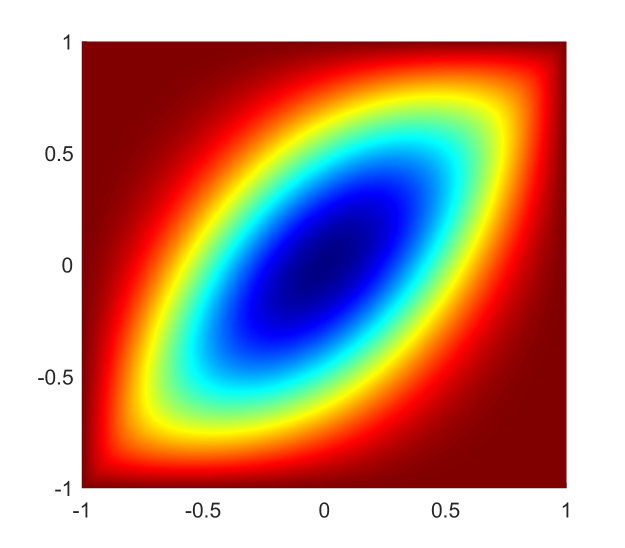}
          \caption{First EV ($\mu=1.25$)}
     \end{subfigure}
     \begin{subfigure}{0.32\textwidth}
         \includegraphics[height=4.5cm,width=5.5cm]{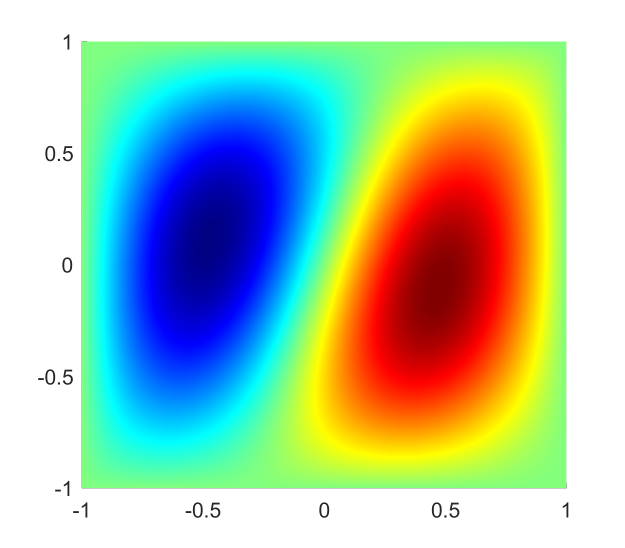}
         \caption{Second EV ($\mu=0.5$)}
     \end{subfigure}
     \begin{subfigure}{0.32\textwidth}
         \includegraphics[height=4.5cm,width=5.5cm]{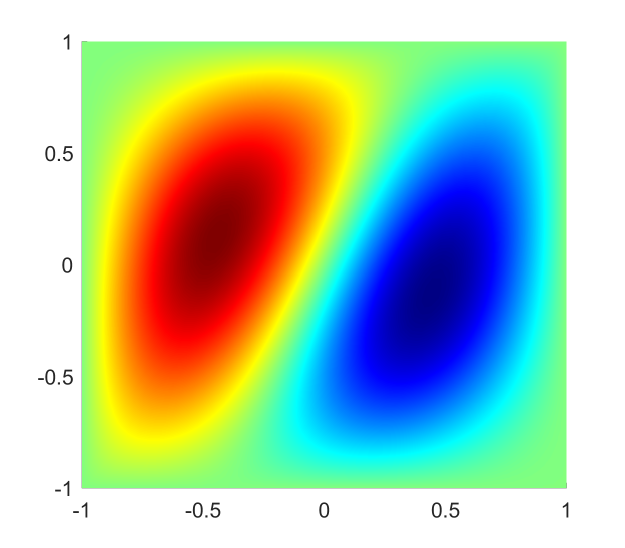}
         \caption{Second EV ($\mu=0.75$)}
     \end{subfigure}
     \begin{subfigure}{0.32\textwidth}
         \includegraphics[height=4.5cm,width=5.5cm]{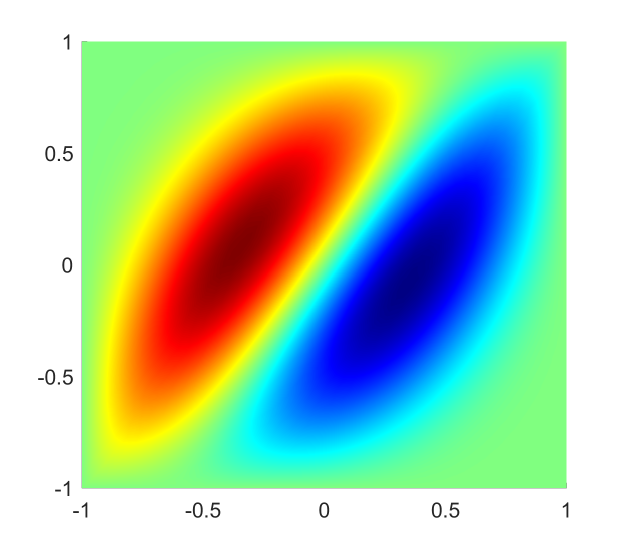}
          \caption{Second EV ($\mu=1.25$)}
     \end{subfigure}
     \begin{subfigure}{0.32\textwidth}
         \includegraphics[height=4.5cm,width=5.5cm]{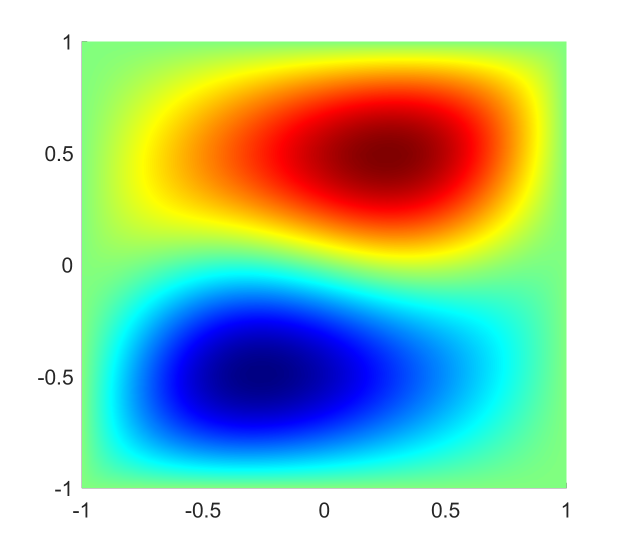}
         \caption{Third EV ($\mu=0.5$)}
     \end{subfigure}
     \begin{subfigure}{0.32\textwidth}
         \includegraphics[height=4.5cm,width=5.5cm]{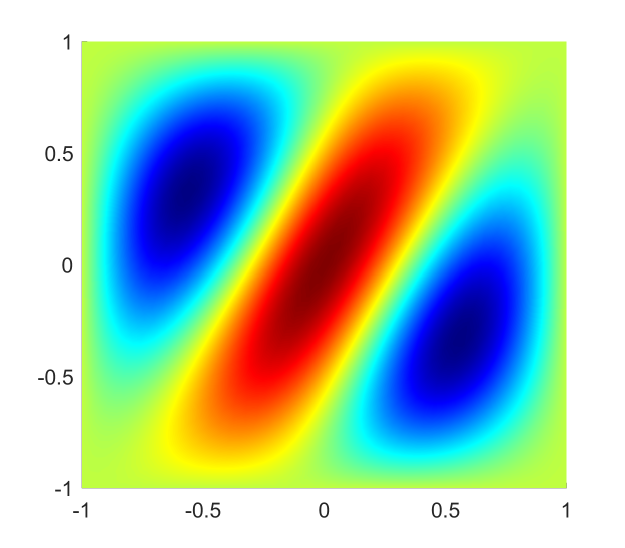}
         \caption{Third EV ($\mu=0.75$)}
     \end{subfigure}
     \begin{subfigure}{0.32\textwidth}
         \includegraphics[height=4.5cm,width=5.5cm]{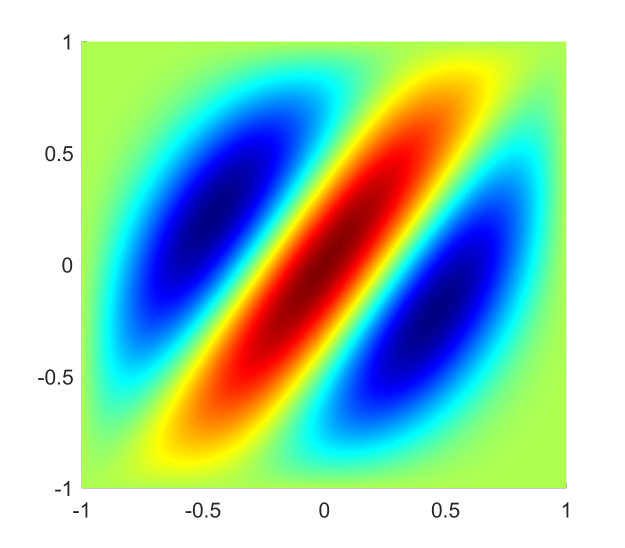}
          \caption{Third EV ($\mu=1.25$)}
     \end{subfigure}
     \begin{subfigure}{0.32\textwidth}
         \includegraphics[height=4.5cm,width=5.5cm]{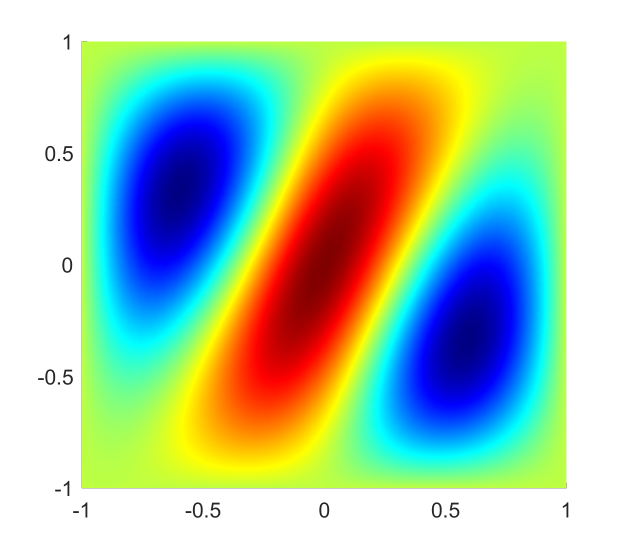}
         \caption{Fourth EV ($\mu=0.5$)}
     \end{subfigure}
     \begin{subfigure}{0.32\textwidth}
         \includegraphics[height=4.5cm,width=5.5cm]{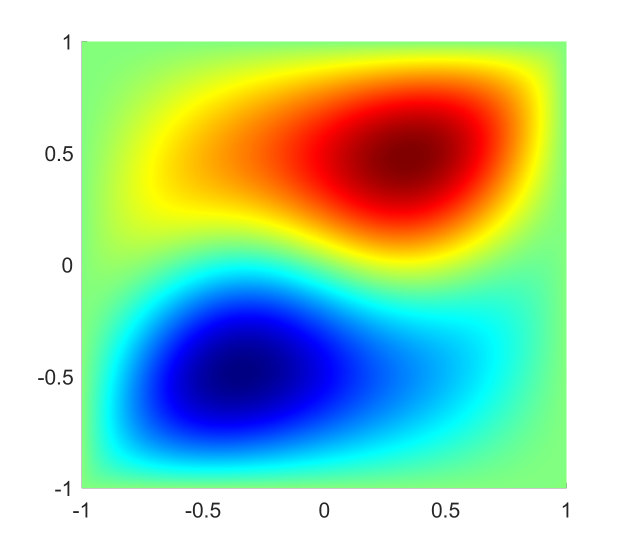}
         \caption{Fourth EV ($\mu=0.75$)}
     \end{subfigure}
     \begin{subfigure}{0.32\textwidth}
         \includegraphics[height=4.5cm,width=5.5cm]{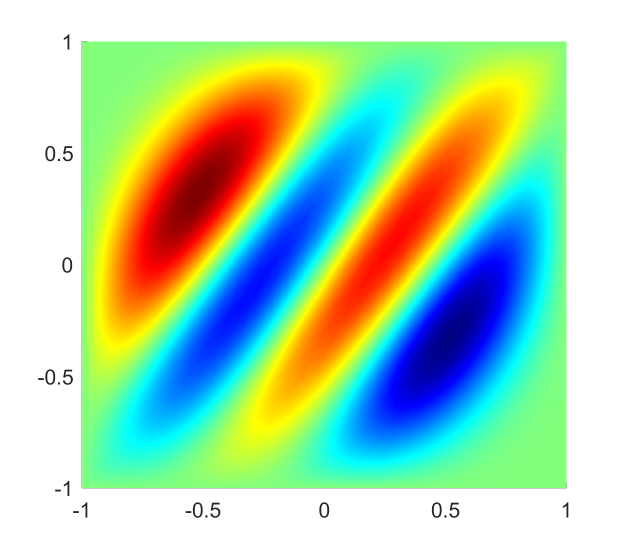}
          \caption{Fourth EV ($\mu=1.25$)}
     \end{subfigure}
        \caption{First four eigenfunctions at $\mu=0.5,0.75,1.25$ using FEM}
        \label{fig:3rdev_1par_pve}
  \end{figure}
\subsection{Reduced order method to obtain $\\lambda_1$ }
\begin{figure}
\centering
\includegraphics[scale=0.4]{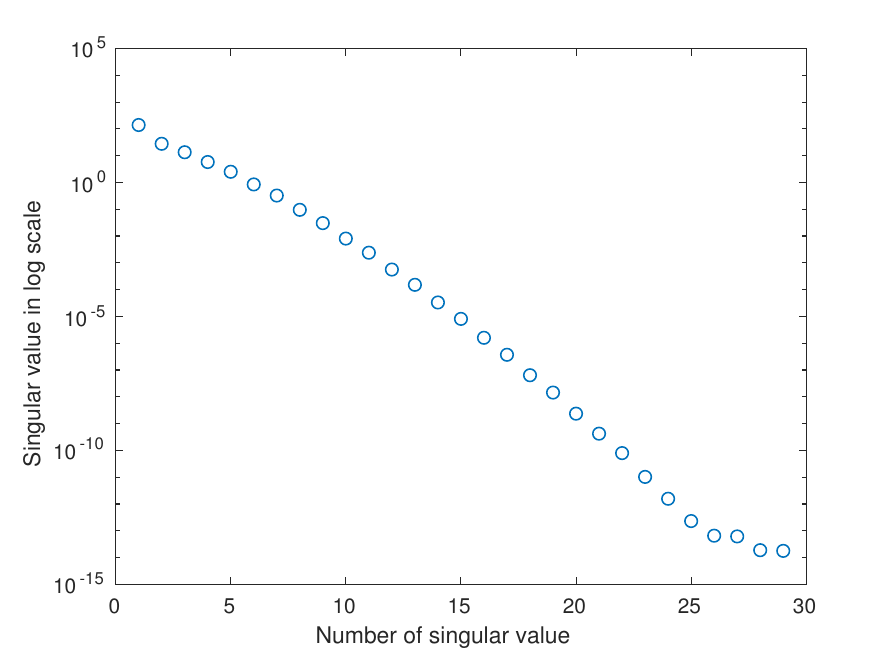}
\includegraphics[scale=0.4]{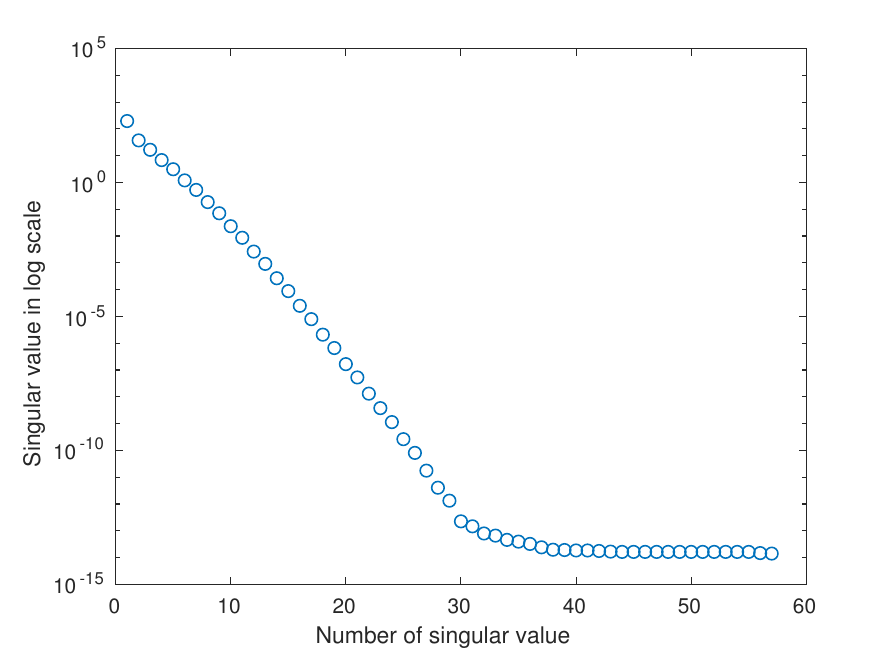}
\caption{Singular values corresponding to snapshot matrix based on $u_1$ when $\mu = -1.4:0.1:1.4$ and $\mu = -1.4:0.05:1.4.$}
\label{FIGURE:6}
\end{figure}

\begin{figure}
\centering
\includegraphics[scale=0.4]{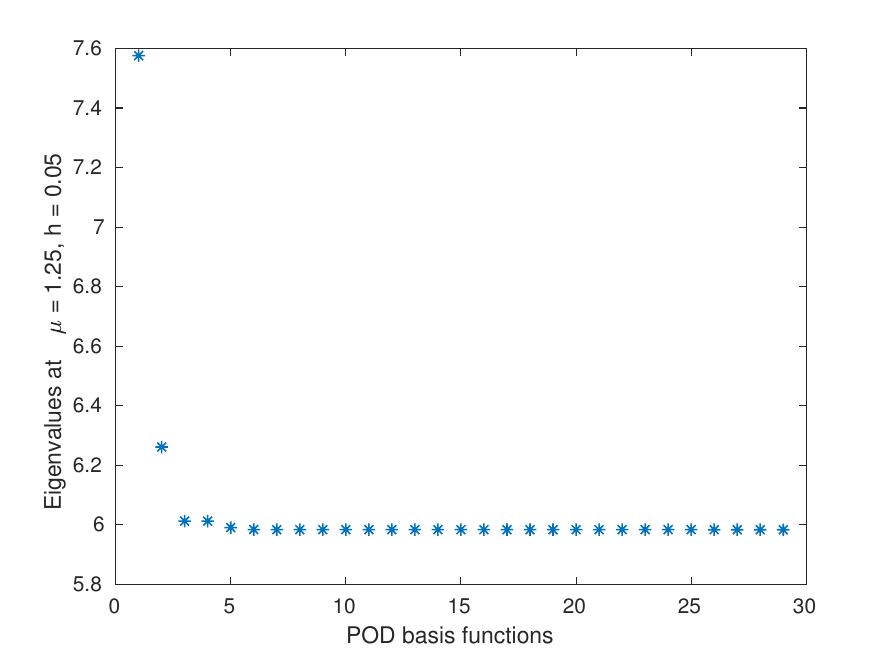}
\includegraphics[scale=0.4]{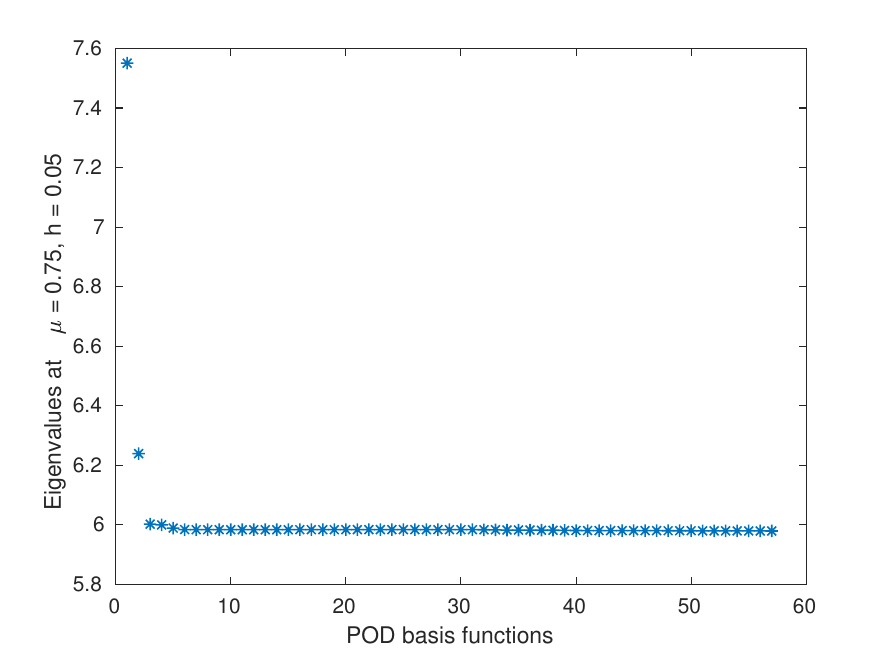}
\caption{Approximation of $\lambda_1$ with snapshot based on $u_1$: eigenvalues corresponding to different number of POD basis functions when $\mu = 1.25$ and $h = 0.05.$ }
\label{FIGURE:7}
\end{figure}

In order to obtain the first eigenvalue using reduced order method, we consider the snapshot matrix consisting of the first eigenvector $u_1$ columnwise at sample parameters. In Figure~\ref{FIGURE:6} we have presented the plot of the singular values of the snapshot matrix for different uniform partitions of the $\mu$ interval and observed that the singular values are decaying very fast. In Figure~\ref{FIGURE:7} we have shown the plot of first eigenvalue at different number of POD basis functions. It is evident from Figure~\ref{FIGURE:7} that the first eigenvalue obtained by ROM is converging to the first eigenvalue of FEM. The first eigenvalues obtained by FEM and ROM  and their relative errors at six sample points are reported in Table~\ref{TABLE:1}. The ROM dimensions mentioned in Table~\ref{TABLE:1} are obtained using the criterion \eqref{criterion} with tolerance $10^{-8}$ and it can be noted that considering very low ROM dimension, we have achieved an accuracy of order $10^{-7}$ -- $10^{-8}$.

\begin{table}[H]
\centering
\begin{tabular}{|c|c|c|c|c|c|c|c|} 
 \hline
{\begin{tabular}[c]{@{}c@{}} ROM dimension \end{tabular}} &
{\begin{tabular}[c]{@{}c@{}} $\mu$ \end{tabular}} &
 {\begin{tabular}[c]{@{}c@{}} 1st EV(FEM) \end{tabular}} & 
 {\begin{tabular}[c]{@{}c@{}}  1st EV(ROM)  \end{tabular}}&
{\begin{tabular}[c]{@{}c@{}}  Relative error \end{tabular}} \\
  \hline
   9 &-1.25&5.98379108& 5.98379387 & 4.7 $\times 10^{-7}$\\
    &-0.75&7.00305328&7.00305363 & 5.0 $\times 10^{-8}$\\
    &-0.50&7.23588322	&7.23588390	 & 9.3 $\times 10^{-8}$\\
    & 0.50&7.23585871 &7.23585938 & 9.2 $\times 10^{-8}$\\
    &0.75&7.00299299&7.00299335 & 5.1 $\times 10^{-8}$\\
    &1.25 & 5.98368037& 5.98368313 &  4.6 $\times 10^{-7}$\\ 
   \hline
   
  10 & -1.25&5.98379108& 5.98379164  & 9.4 $\times 10^{-8}$\\
   &-0.75&7.00305328&7.00305360 & 5.8 $\times 10^{-8}$\\
  &-0.50&7.23588322	&7.23588349	 & 3.6 $\times 10^{-8}$\\
  & 0.50&7.23585871 &	7.23585900 & 3.9 $\times 10^{-8}$\\
  & 0.75&7.00299299&7.00299332 & 5.4 $\times 10^{-8}$\\
  & 1.25 & 5.98368037& 5.98368097 & 1.0 $\times 10^{-7}$  \\ 
  \hline
\end{tabular}
   \caption{Approximation of $\lambda_1$ with snapshot based on $u_1$: comparison of FEM and ROM at $h = 0.05$}
 \label{TABLE:1}
 \end{table}

\subsection{Reduced order method to obtain $\lambda_{2}$}
In this subsection, we discuss numerical results for the second eigenvalue considering different combinations of eigenvectors in the snapshot matrix.
\subsubsection{Results of the EVP considering $u_2$ in the snapshot matrix}
\begin{figure}[H]
\centering
\includegraphics[scale=0.4]{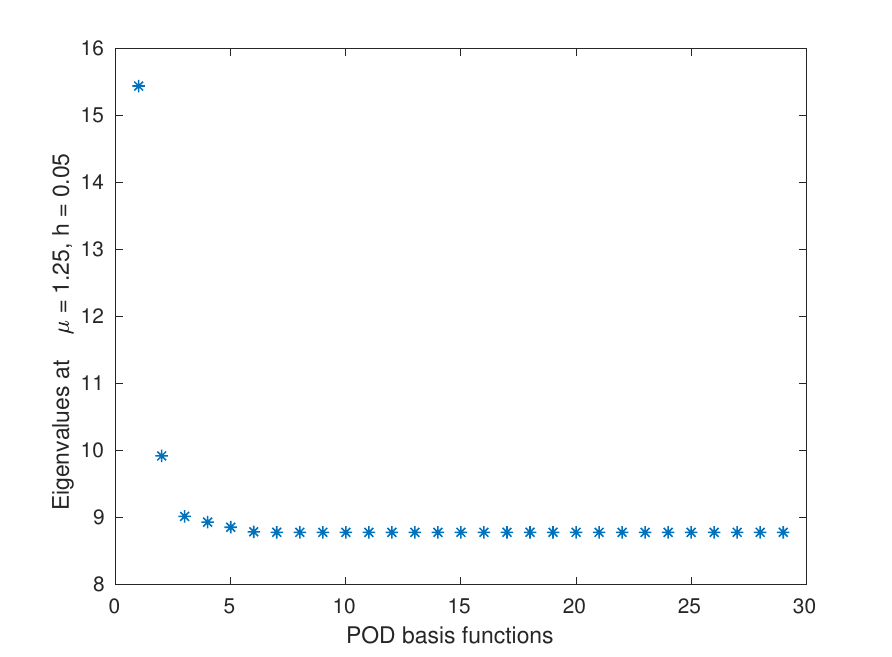}
\includegraphics[scale=0.4]{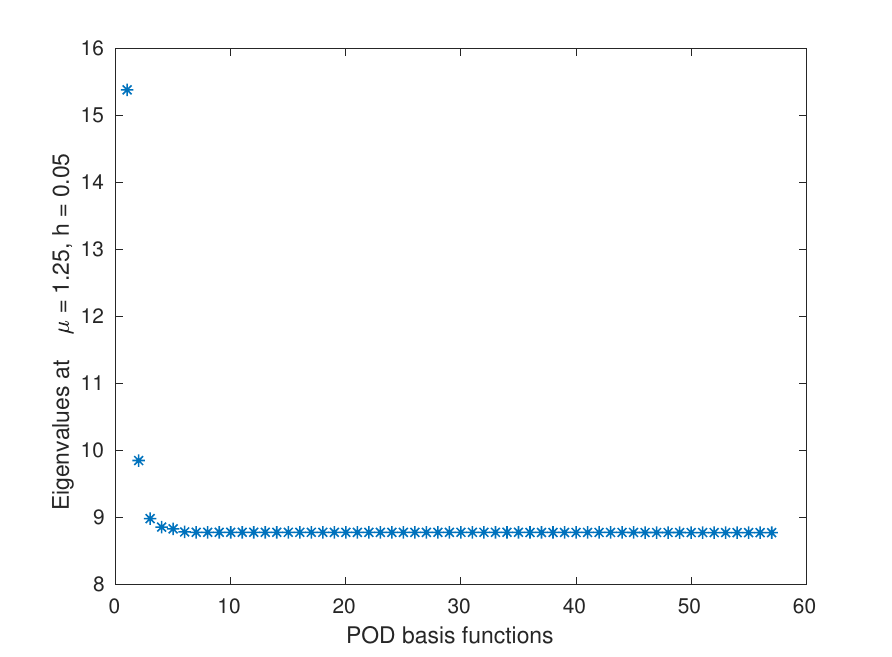}
\caption{Approximation of $\lambda_2$ with snapshot based on $u_2$: eigenvalues corresponding to different number of POD basis functions when $\mu = 1.25$ and $h = 0.05.$ }
\label{FIGURE:8}
\end{figure}

In order to obtain the second eigenvalue using reduced order method, we consider the snapshot matrix consisting of the second eigenvector $u_2$ columnwise at sample parameters. In Figure~\ref{FIGURE:8} we have shown the plot of the first eigenvalue of ROM for different number of POD basis functions. It is evident from Figure~\ref{FIGURE:8} that the first eigenvalue obtained by ROM is converging to the second eigenvalue of FEM as expected. 

The first  eigenvalues of ROM, the second eigenvalues of FEM and their relative errors are reported in Table~\ref{TABLE:2}.  It can be noted that considering very low ROM dimension the relative errors are of order $10^{-7}$ -- $10^{-8}$.

\begin{table}[H]
 	 	\centering
 	\begin{tabular}{|c|c|c|c|c|c|c|c|} 
 		\hline

 	        {\begin{tabular}[c]{@{}c@{}} ROM dimension \end{tabular}} &
		  {\begin{tabular}[c]{@{}c@{}} $\mu$ \end{tabular}} &
 		 {\begin{tabular}[c]{@{}c@{}} 2nd EV(FEM) \end{tabular}} & 
 		{\begin{tabular}[c]{@{}c@{}}  1st EV(ROM)  \end{tabular}}&
        {\begin{tabular}[c]{@{}c@{}} Relative error \end{tabular}}\\

  \hline

  11 &-1.25&8.77547249	& 8.77547548 & 3.4 $\times 10^{-7}$\\
    &-0.75&12.86380150&12.86380201 & 4.0 $\times 10^{-8}$\\
   &-0.50 & 13.95750227	&13.95750231 & 3.3 $\times 10^{-9}$\\
   &0.50 & 13.95736933 &	13.95736937 & 3.2 $\times 10^{-9}$\\
   &0.75 & 12.86364140&12.86364192 & 4.0 $\times 10^{-8}$	\\
   &1.25 & 8.77689931& 8.77690230  & 3.4 $\times 10^{-7}$\\ 
   \hline
  
   12 & -1.25&8.77547249 & 8.77547254 & 4.8 $\times 10^{-9}$	\\
   &-0.75&12.86380150&12.86380156 &  5.1 $\times 10^{-9}$ \\
  & -0.50 & 13.95750227	&13.95750256	&	2.0 $\times 10^{-8}$ \\
  & 0.50 & 13.95736933 &13.95736962 & 2.0 $\times 10^{-8}$\\
  & 0.75&12.86364140&12.86364147 & 5.0 $\times 10^{-9}$\\
  & 1.25 & 8.77689931&8.77689936 & 4.9 $\times 10^{-9}$\\ 
 \hline
   \end{tabular}
 	\caption{Approximation of $\lambda_2$ with snapshot based on $u_2$: comparison of FEM and ROM at $h = 0.05$}
  \label{TABLE:2}
 	\end{table}

\subsubsection{Results of the EVP considering $u_1 + u_2$  is in the snapshot matrix}
\begin{figure}
\centering
\includegraphics[scale=0.4]{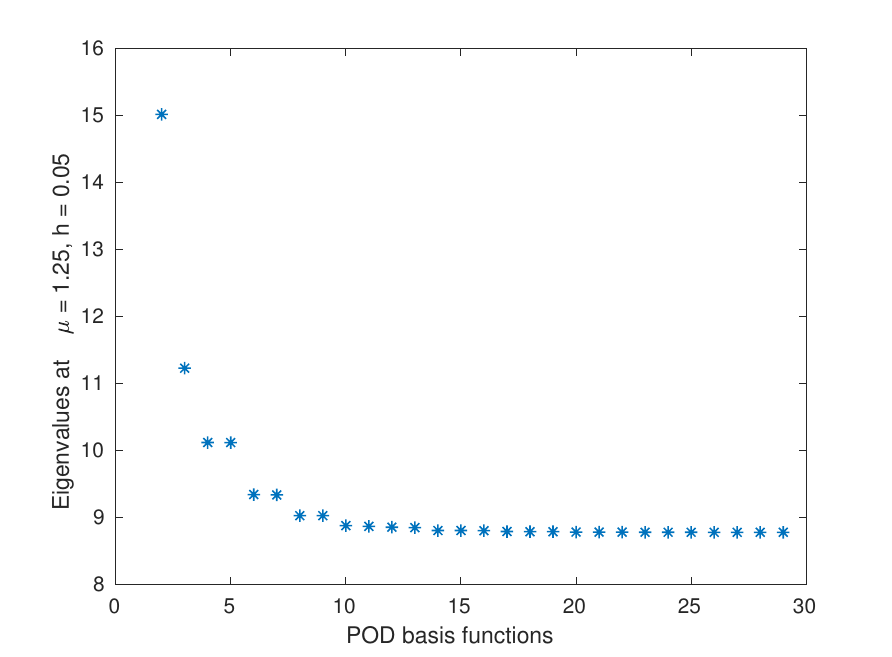}
\includegraphics[scale=0.4]{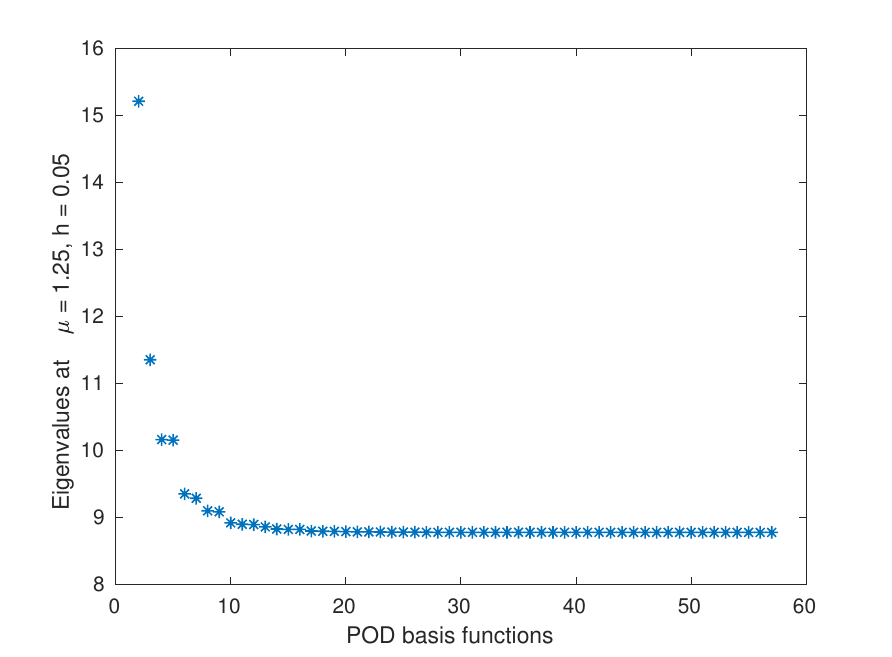}
\caption{Approximation of $\lambda_2$ with snapshot based on $u_1 + u_2$: eigenvalues corresponding to different number of POD basis functions when $\mu = 1.25$ and $h = 0.05.$  }
\label{FIGURE:9}
\end{figure}

We consider the snapshot matrix consisting of the combination of eigenvectors $u_1 + u_2$ column wise at sample parameters in order to compute the second eigenvalue using the reduced order method. In Figure~\ref{FIGURE:9} we have presented the plot of the second eigenvalue of ROM at different number of POD basis functions. The second eigenvalues obtained by FEM and ROM and their relative errors are reported in Table~\ref{TABLE:3}.

\begin{table}
 	 	\centering
 	\begin{tabular}{|c|c|c|c|c|c|c|c|} 
 		\hline

 	        {\begin{tabular}[c]{@{}c@{}} ROM dimension \end{tabular}} &
		  {\begin{tabular}[c]{@{}c@{}} $\mu$ \end{tabular}} &
 		 {\begin{tabular}[c]{@{}c@{}} 2nd EV(FEM) \end{tabular}} & 
 		{\begin{tabular}[c]{@{}c@{}}  2nd EV(ROM)  \end{tabular}}&
 		{\begin{tabular}[c]{@{}c@{}}  Relative error  \end{tabular}}\\
 	
  \hline

  14 &-1.25&8.77547249	& 8.78969339	&	$1.6 \times 10^{-3}$	 \\
    &-0.75&12.86380150 & 12.86380339 & $1.4 \times 10^{-7}$\\
   &-0.50&13.95750227	& 	13.95750436	& $1.5 \times 10^{-7}$\\
   &0.50&13.95736933 & 13.95737012 & $5.6 \times 10^{-8}$\\
    &0.75&12.86364140 & 	12.86379512 & $1.2 \times 10^{-5}$\\
   &1.25 & 	8.77689931 & 	8.81325897	&	$4.1 \times 10^{-3}$\\ 
   \hline
  
   15 & -1.25&8.77547249 	&8.79049305		&		$1.7 \times 10^{-3}$				\\
    &-0.75&12.86380150& 12.86380446 & $2.2 \times 10^{-7}$\\
    &-0.50&13.95750227	 & 	13.95750469	& $1.7 \times 10^{-7}$\\
   &0.50&13.95736933 & 	13.95737042  & $7.8 \times 10^{-8}$\\
    & 0.75&12.86364140 & 	12.86390127  & $2.0 \times 10^{-5}$\\
    & 1.25 & 8.77689931 & 	8.82153162 & $5.0 \times 10^{-3}$\\ 
  \hline
   \end{tabular}
   	\caption{Approximation of $\lambda_2$ with snapshot based on $u_1 + u_2$: comparison of FEM and ROM at $h = 0.05$}
    \label{TABLE:3}
 	\end{table}

It can be seen from the Table~\ref{TABLE:2} and Table~\ref{TABLE:3} that considering only $u_2$ in the snapshot matrix provides us slightly better results than considering $u_1 + u_2$ in the snapshot matrix. However, considering $u_1 + u_2$ in the snapshot matrix, both the eigenvalues $\lambda_1$ and $\lambda_2$ can be obtained simultaneously.

\subsection{Reduced order method to obtain $\lambda_{3}$}
In this subsection, we discuss numerical results for the third eigenvalue considering different combinations of eigenvectors in the snapshot matrix. It can be seen from Figure~\ref{fig3} that the third and the fourth eigenvalues are intersecting at some values of $\mu$. We observe many interesting phenomena in this case.
\subsubsection{Results of the EVP considering only $u_3$ in the snapshot matrix}
\begin{figure}
\centering
\includegraphics[scale=0.6]{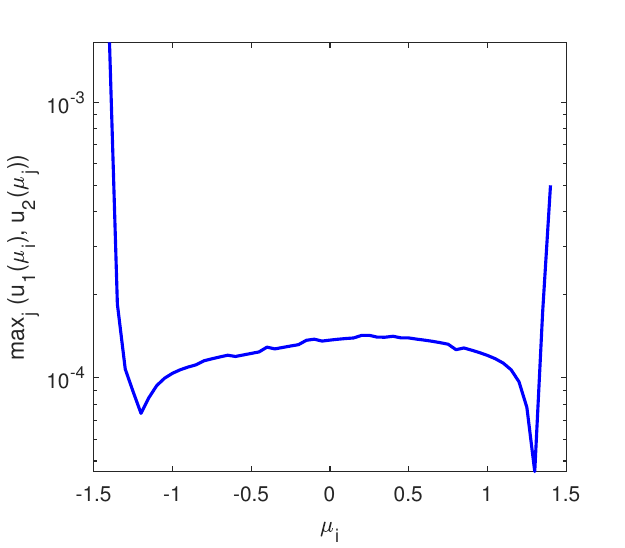}
\includegraphics[scale=0.6]{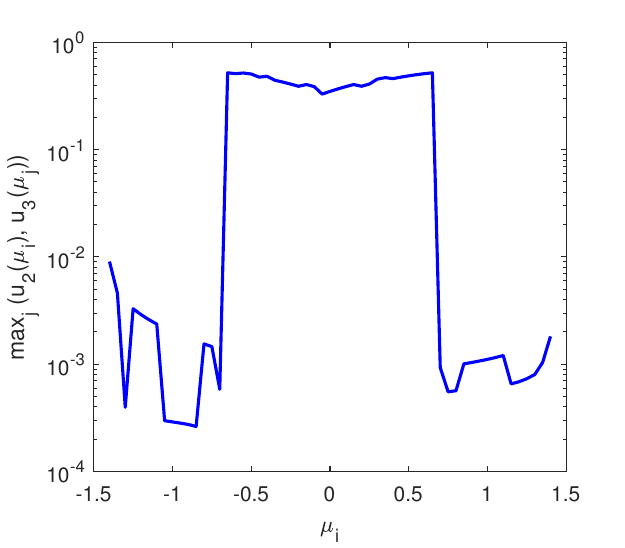}
\caption{\color{magenta}Scalar product of eigenfunctions for different parameters. Left: $\max_j(u_1(\mu_i),u_2(\mu_j))$ as a function of $\mu_i$; right: $\max_j(u_2(\mu_i),u_3(\mu_j)$ as a function of $\mu_i$. While $u_1(\mu_i)$ and $u_2(\mu_j)$ are almost orthogonal, $u_2(\mu_i)$ and $u_3(\mu_j)$ are not.}
\label{FIGURE:10a}
\end{figure}

\begin{figure}
\centering
\includegraphics[scale=0.4]{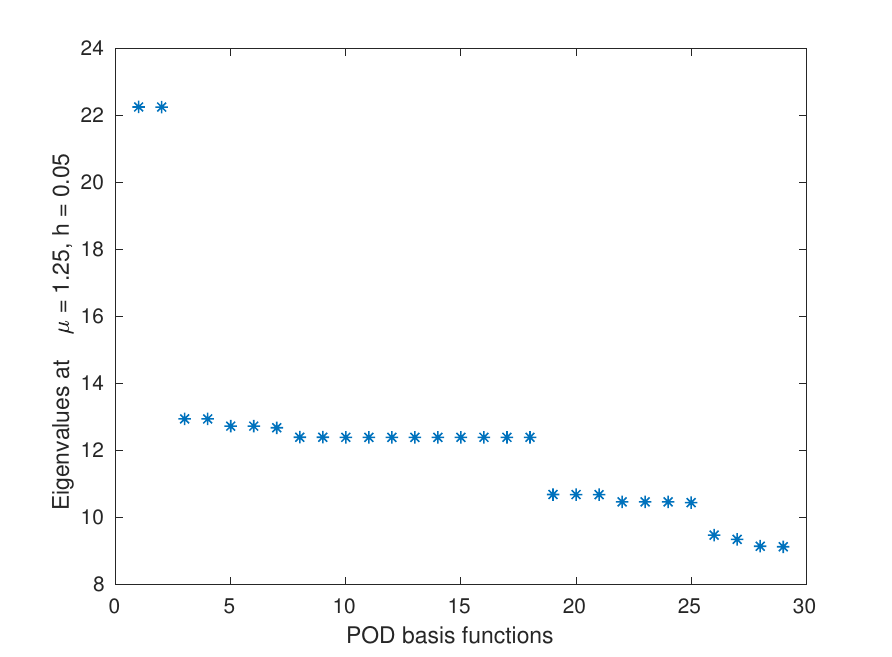}
\includegraphics[scale=0.4]{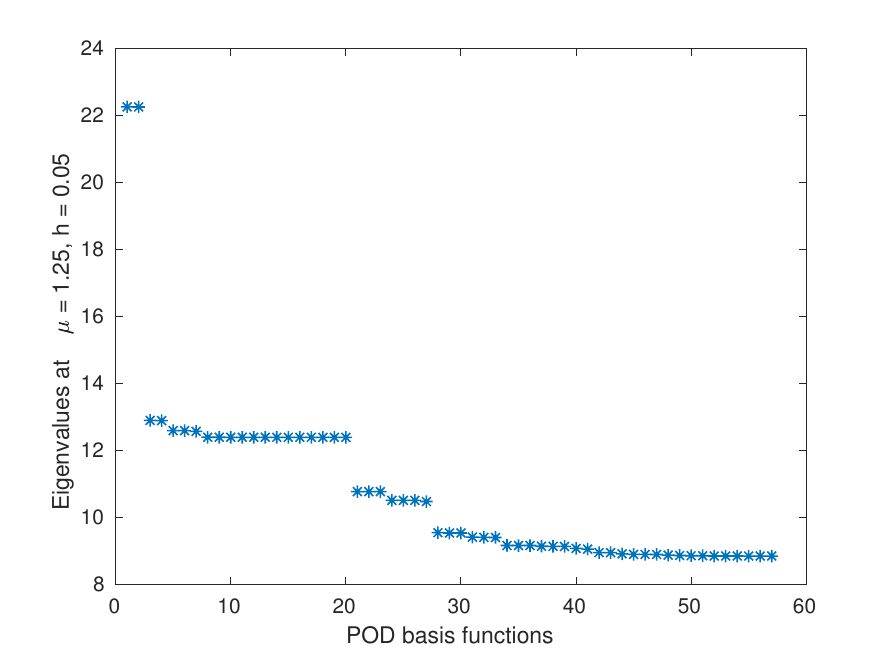}
\caption{Approximation of $\lambda_3$ with snapshot based on $u_3$: eigenvalues corresponding to different number of POD basis functions when $\mu = 1.25$ and $h = 0.05.$ }
\label{FIGURE:10}
\end{figure}

\begin{figure}
     \begin{subfigure}{0.32\textwidth}
         \includegraphics[height=4.5cm,width=5.5cm]{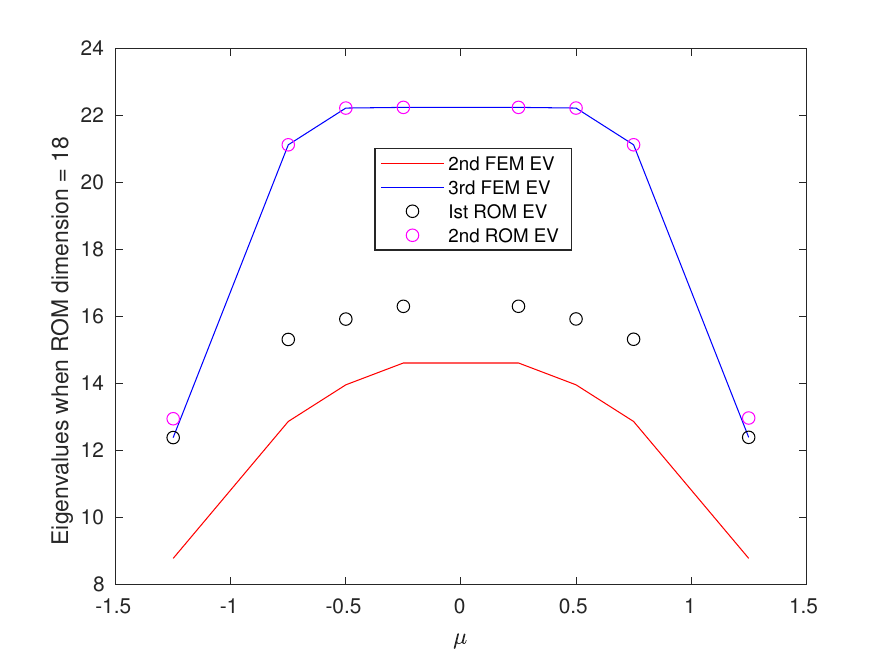}
         \caption{ROM dimension 18}
     \end{subfigure}
     \begin{subfigure}{0.32\textwidth}
         \includegraphics[height=4.5cm,width=5.5cm]{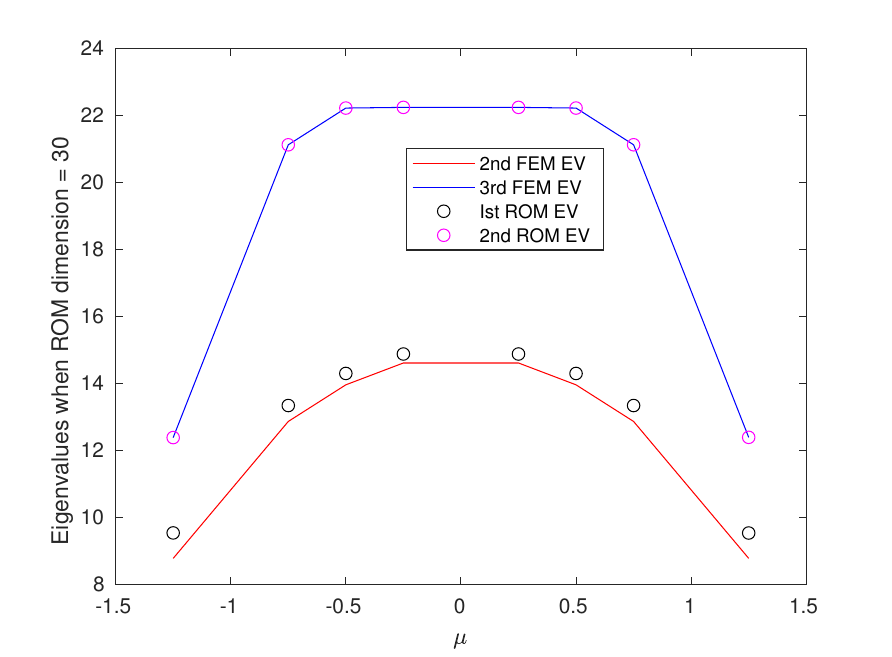}
         \caption{ROM dimension 30}
     \end{subfigure}
     \begin{subfigure}{0.32\textwidth}
         \includegraphics[height=4.5cm,width=5.5cm]{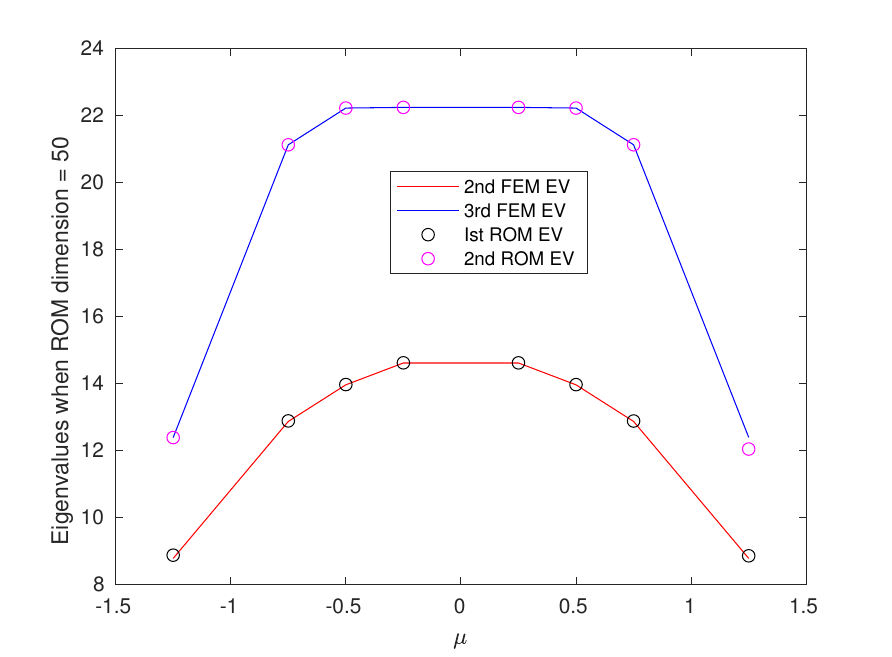}
          \caption{ROM dimension 50}
     \end{subfigure}
     \caption{FEM and ROM based eigenvalues at different ROM dimensions}
     \label{FIGURE:11}
\end{figure}

\begin{figure}
\centering
\includegraphics[scale=0.5]{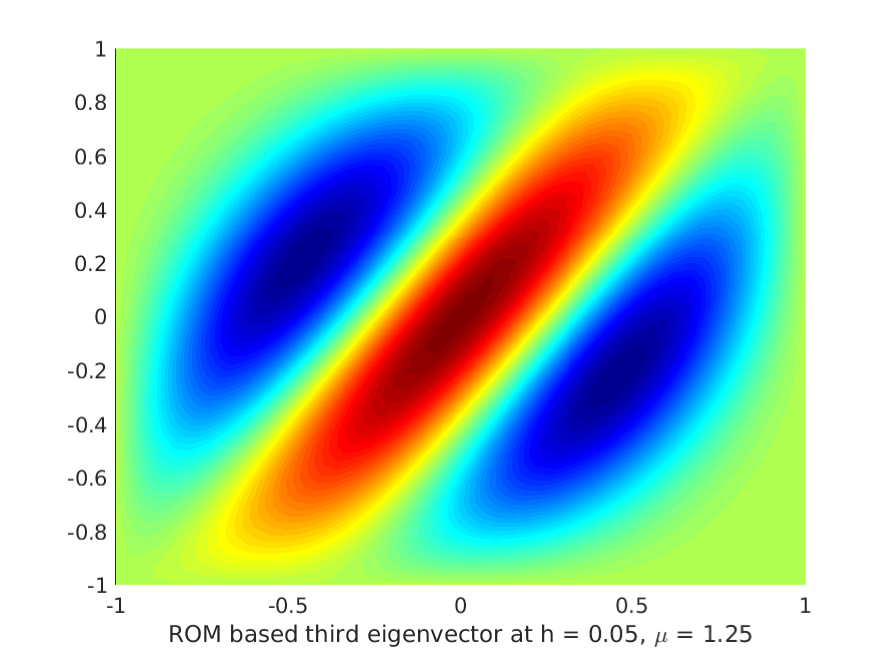}
\includegraphics[scale=0.5]{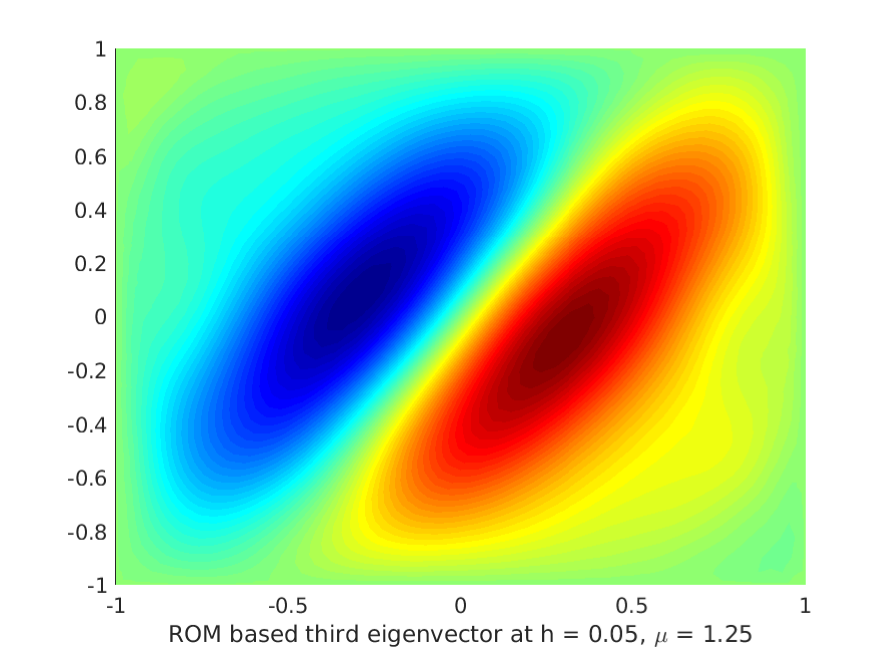}
\caption{Approximation of $u_3$ with snapshot based on $u_3$ when ROM dimension is 18 and 30, $\mu = 1.25$ and $h = 0.05.$ }
\label{FIGURE:12}
\end{figure}

Let us consider the snapshot matrix consisting only  third eigenvector at the sample points and compute the eigenvalues using reduced order method. Considering the criterion \eqref{criterion} with tolerance $10^{-8}$, the ROM dimension turns out to be 17 and 18 respectively. In this case, the first eigenvalue of ROM is not matching with the third eigenvalue of FEM, see Table~\ref{TABLE:4}. However, when we increase the number of POD basis functions, we notice that the
first eigenvalue of the ROM matches with the second eigenvalue of the FEM, whereas the second eigenvalue of the ROM matches with the third eigenvalue of the FEM, see Figure~\ref{FIGURE:11}. When we consider only the third eigenvector $u_{3}$ in the snapshot matrix, we were expecting that the first eigenvalue of the ROM would match with the third eigenvalue of the FEM, but this is not the case. The reason behind this is the fact that the $L^2$ inner product $(u_2(\mu_i),u_3(\mu_j))$ is not zero for $i\neq j$, so the snapshots contain some component of the second eigenvector. Note that the inner product is zero for $i=j$.
\rev{In Figure~\ref{FIGURE:10a}, we report the $L^2$ inner product of $(u_1(\mu_i),u_2(\mu_j))$ and $(u_2(\mu_i),u_3(\mu_j))$ as a function of $\mu_i$, respectively. The figure supports our claim that $u_1(\mu_i)$ and $u_2(\mu_j)$ are almost orthogonal, while $u_2(\mu_i)$ and $u_3(\mu_j)$ are not.} In Figure~\ref{FIGURE:10} we have presented the plot of first eigenvalue of the ROM at different number of POD basis functions. It can be easily seen from Figure~\ref{FIGURE:10} that the first eigenvalue of the ROM is not matching with the third eigenvalue of the FEM. It turns out that the third eigenvector obtained by the ROM is not the same using different number of POD basis functions, see Figure~\ref{FIGURE:12}. 


\begin{table}
\centering
\begin{tabular}{|c|c|c|c|c|c|c|c|} 
 \hline
 {\begin{tabular}[c]{@{}c@{}} ROM dimension \end{tabular}} &
 {\begin{tabular}[c]{@{}c@{}} $\mu$ \end{tabular}} &
 {\begin{tabular}[c]{@{}c@{}} 3rd EV(FEM) \end{tabular}} & 
 {\begin{tabular}[c]{@{}c@{}}  1st EV(ROM)  \end{tabular}}\\
  \hline
   17 &-1.25&12.38270171				& 12.38271385\\
    &-0.75&21.12292172&15.19805477\\
    &-0.50& 22.22010384  &15.81534625	\\
    &0.50&	22.22018744 &15.81572205\\
    &0.75&21.12370790&15.19851153\\
    &1.25 & 	12.38958778	& 12.38960108 \\ 
   \hline
  18 &-1.25&12.38270171	& 12.38270216\\
    &-0.75&21.12292172&15.31186360\\
    &-0.50& 22.22010384  &15.92038351		\\
    &0.50&	22.22018744 &15.92406607\\
    &0.75&21.12370790&15.31808566\\
    &1.25 & 	12.38958778	& 12.38958815 \\ 
  \hline
   \end{tabular}
   	\caption{Approximation of $\lambda_3$ with snapshot based on $u_3$ at $h = 0.05$}
	 	\label{TABLE:4}
 	\end{table}
 	
\begin{figure}
\centering
     \begin{subfigure}{0.32\textwidth}
         \includegraphics[height=4.5cm,width=5.5cm]{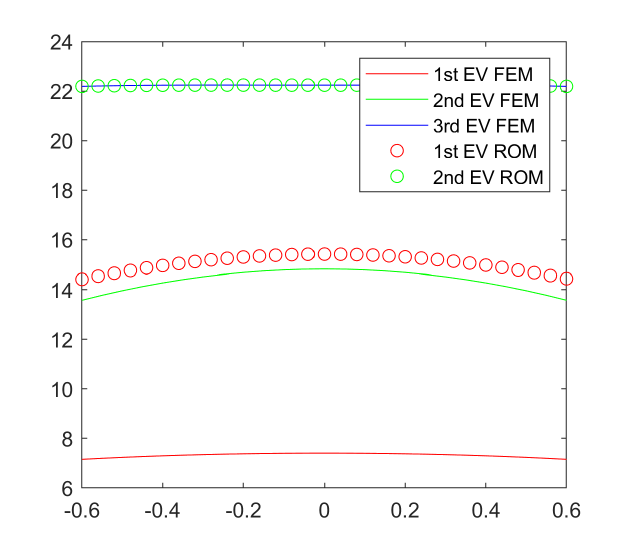}
         \caption{ROM dimension 10}
     \end{subfigure}
      \begin{subfigure}{0.32\textwidth}
         \includegraphics[height=4.5cm,width=5.5cm]{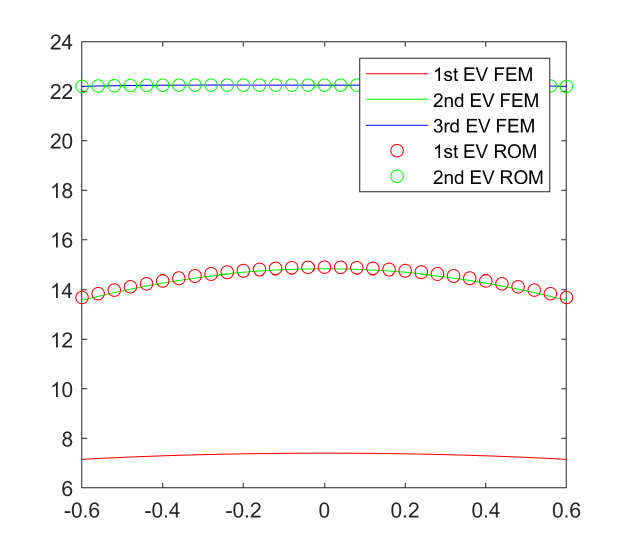}
         \caption{ROM dimension 15}
     \end{subfigure}
      \begin{subfigure}{0.32\textwidth}
         \includegraphics[height=4.5cm,width=5.5cm]{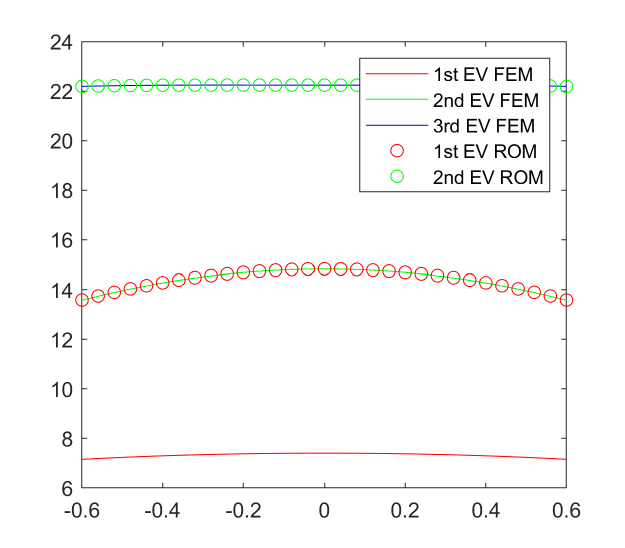}
         \caption{ROM dimension 20}
     \end{subfigure}
     \caption{Approximation of eigenvalues with snapshot based on $u_3$ and $\mathcal{M}=(-0.6,0.6)$ }
     \label{FIGURE:13}
\end{figure}
In Figure~\ref{FIGURE:13} we have shown the FEM and ROM eigenvalues for different number of POD basis functions for the parameter space $\mathcal{M}=(-0.6,0.6)$. Note that for this particular parameter interval, the two eigenvalues do not intersect (see Figure~\ref{fig3}). From the plot one can see that also in this case the situation is not as naively expected: the first eigenvalue of the ROM is matching with the second eigenvalue of the FEM and the second eigenvalue of the ROM is matching with the third eigenvalue of the FEM.

\subsubsection{Results of the EVP considering $u_3$ and $u_4$ in the snapshot matrix}
\begin{figure}[H]
\centering
\includegraphics[scale=0.35]{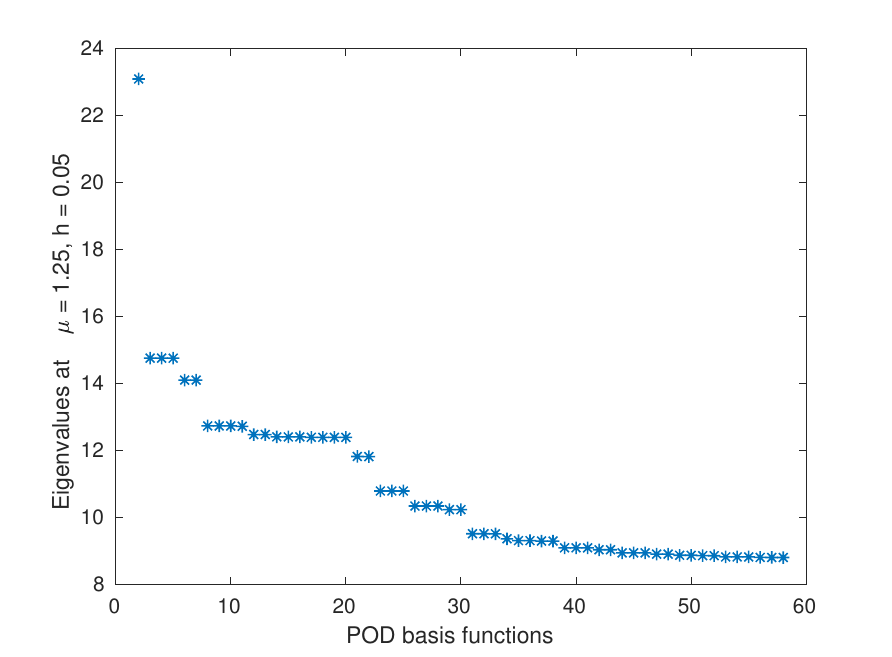}
\includegraphics[scale=0.35]{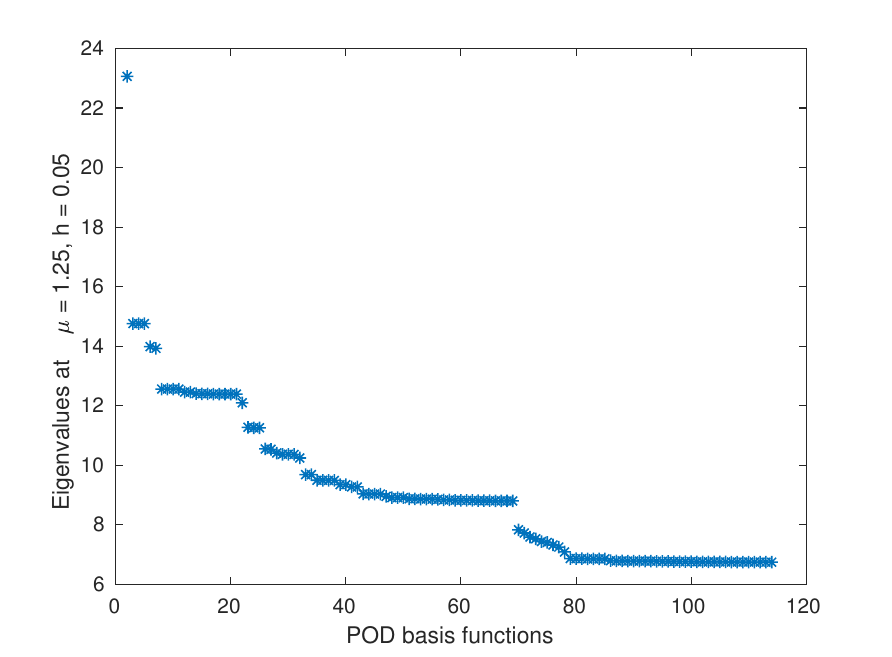}
\caption{Approximation of $\lambda_3$ with snapshot based on $u_3$ and $u_4$: eigenvalues corresponding to different number of POD basis functions when $\mu = 1.25$ and $h = 0.05.$ }
\label{FIGURE:14}
\end{figure}

Since the third and the fourth eigenvalues are intersecting, we consider $u_3$ and $u_4$ in the snapshot matrix and compute the eigenvalues using the reduced order method. In Figure~\ref{FIGURE:14} we have presented the plot of the first eigenvalue obtained by the ROM. It can be easily seen that the first eigenvalue of ROM is not converging to the third eigenvalue of the FEM. However, the second eigenvalue of the ROM is converging to the third eigenvalue of FEM. In Table~\ref{TABLE:5} we have reported the first and second eigenvalues of the ROM and the third eigenvalue of the FEM. \rev{Further, the third eigenvalue of the ROM converges to the fourth eigenvalue of the FEM. In Table~\ref{TABLE:5a}, we have reported the third eigenvalue of the ROM and the fourth eigenvalue of the FEM and the relative error between them.}

\begin{table}[H]
 	 	\centering
 	\begin{tabular}{|c|c|c|c|c|c|c|c|} 
 		\hline
 	        {\begin{tabular}[c]{@{}c@{}} ROM dimension \end{tabular}} &
		  {\begin{tabular}[c]{@{}c@{}} $\mu$ \end{tabular}} &
 		 {\begin{tabular}[c]{@{}c@{}} 3rd EV(FEM) \end{tabular}} & 
 		{\begin{tabular}[c]{@{}c@{}}  1st EV(ROM)  \end{tabular}}&
 		{\begin{tabular}[c]{@{}c@{}}  2nd EV(ROM)  \end{tabular}}\\
 	
  \hline
  26 &-1.25&12.38270171 & 10.32768230 & 12.38271820\\
    &-0.75&21.12292172&14.01748876&21.12292540\\
   &-0.50& 22.22010384&14.96276873	  &22.22010448		\\
    &0.50&	22.22018744&14.96700836 &22.22018808	\\
   &0.75&21.12370790&14.02396381&21.12371133\\
    &1.25 & 	12.38958778	& 10.33922731&12.38960275\\ 
   \hline
   30 &-1.25&12.38270171& 10.35774089& 12.38270198\\
    &-0.75&21.12292172& 14.02551777&           21.12292203\\
   &-0.50& 22.22010384 & 14.96893932			  &22.22010447	\\
    &0.50&	22.22018744&14.97067172 &	22.22018807		\\
    &0.75&21.12370790&  14.02852884&           21.12370817\\
    &1.25 & 	12.38958778	&  10.36928581 &      12.38958799 \\ 
  \hline
   \end{tabular}
	\caption{Approximation of $\lambda_3$ with snapshot based on $u_3$ and $u_4$ at $h = 0.05$}
	 	\label{TABLE:5}
 	\end{table}
 	
\begin{table}[H]
\centering
\begin{tabular}{|c|c|c|c|c|c|c|c|} 
 \hline
{\begin{tabular}[c]{@{}c@{}} ROM dimension \end{tabular}} &
{\begin{tabular}[c]{@{}c@{}} $\mu$ \end{tabular}} &
{\begin{tabular}[c]{@{}c@{}} 4th EV(FEM) \end{tabular}} & 
{\begin{tabular}[c]{@{}c@{}}  3rd EV(ROM)  \end{tabular}}&
{\begin{tabular}[c]{@{}c@{}}  Relative Error  \end{tabular}}\\
\hline
  26&-1.25 &16.13593156 & 16.13594080 & $5.7\times 10^{-7}$\\
    &-0.75 &22.08781810 & 22.08782217  & $1.8\times 10^{-7}$\\
    &-0.50 &24.02185801 & 24.02186412 & $2.5\times 10^{-7}$\\
    &0.50  &	24.02234525 & 24.02235057  & $2.2\times 10^{-7}$	\\
    &0.75  &22.08790319 & 22.08790743  & $1.9\times 10^{-7}$\\
    &1.25  &16.13984606 & 16.13985482 & $5.4\times 10^{-7}$\\ 
    \hline
   30&-1.25 &16.13593156 & 16.13593252  & $5.9\times 10^{-8}$\\
    &-0.75 &22.08781810 & 22.08781835  & $1.1\times 10^{-}$\\
    &-0.50 &24.02185801 & 24.02185811  & $4.1\times 10^{-9}$\\
    &0.50  &	24.02234525 & 24.02234536  & $4.5\times 10^{-9}$	\\
    &0.75  &22.08790319 & 22.08790344  & $1.1\times 10^{-8}$\\
    &1.25  &16.13984606 & 16.13984698 & $5.6\times 10^{-8}$\\ 
    \hline
   \end{tabular}
\caption{\rev{Approximation of $\lambda_4$ with snapshot based on $u_3$ and $u_4$ at $h = 0.05$}}
\label{TABLE:5a}
\end{table}

\subsubsection{Results of the EVP considering  $u_1$, $u_2$, and $u_3$ in the snapshot matrix}
\begin{figure}
\centering
\includegraphics[scale=0.4]{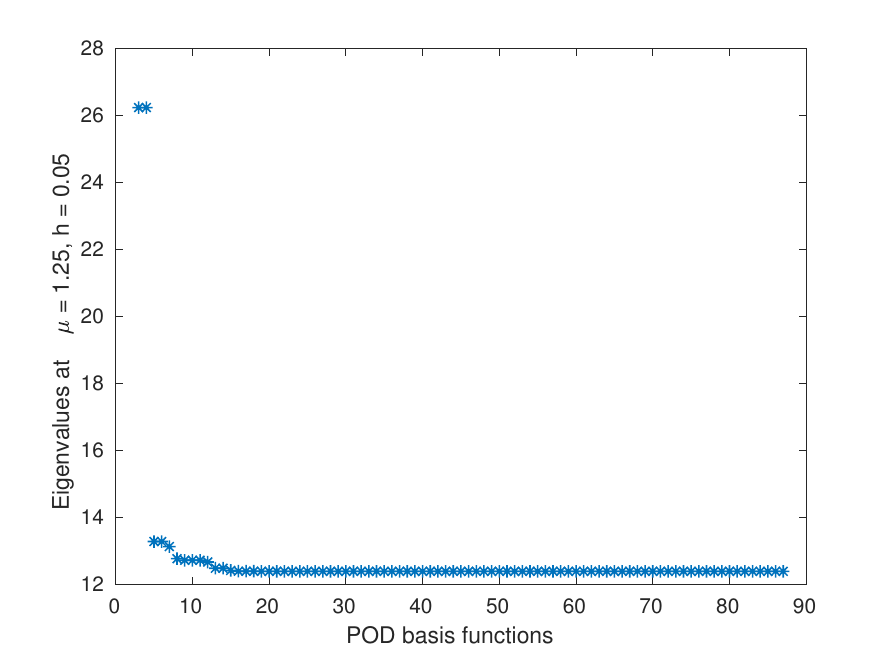}
\includegraphics[scale=0.4]{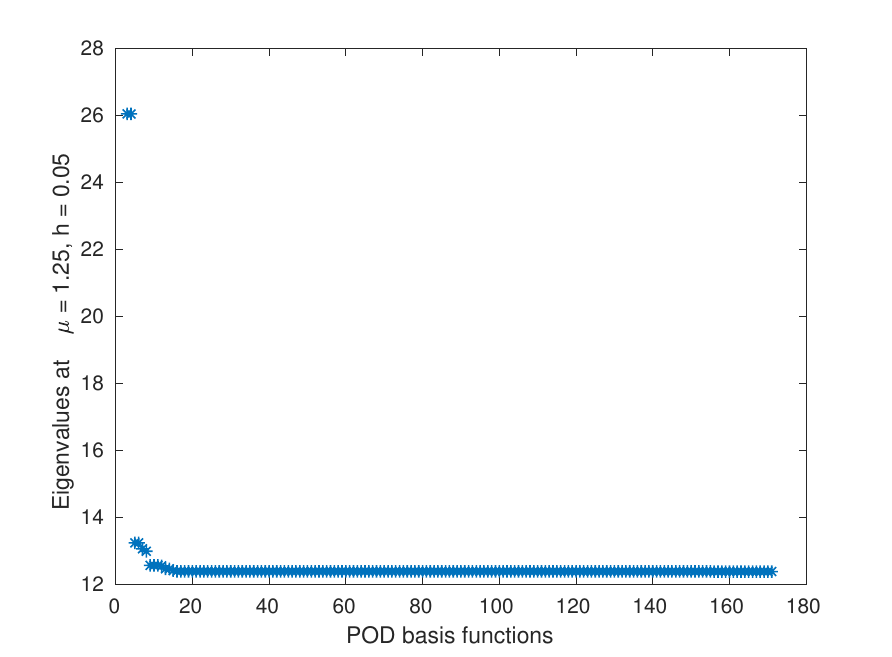}
\caption{Approximation of $\lambda_3$ with snapshot based on $u_1, u_2$ and $u_3$: eigenvalues corresponding to different number of POD basis functions when $\mu = 1.25$ and $h = 0.05.$ }
\label{FIGURE:15}
\end{figure}

\begin{figure}
\centering
\includegraphics[scale=0.4]{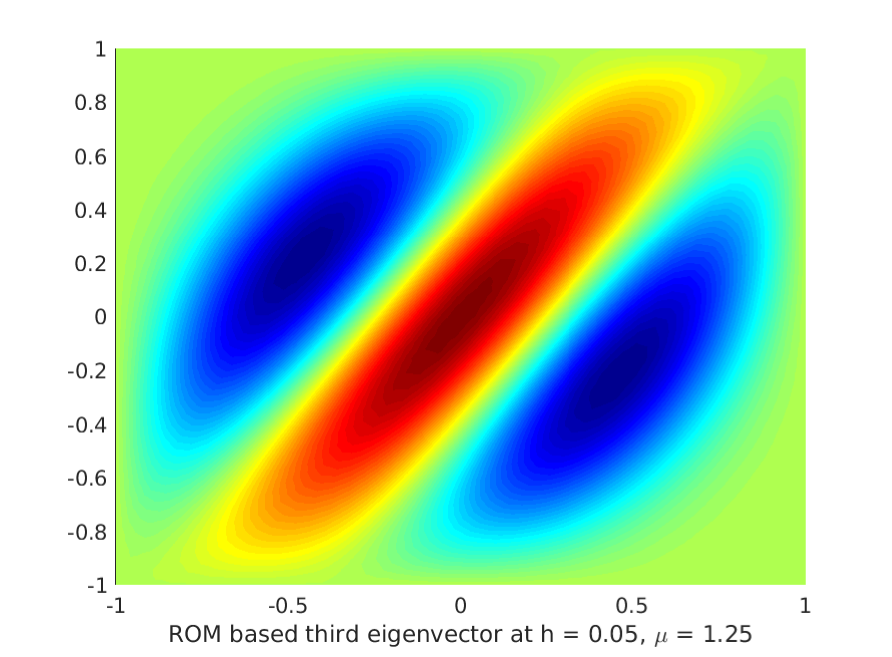}
\includegraphics[scale=0.4]{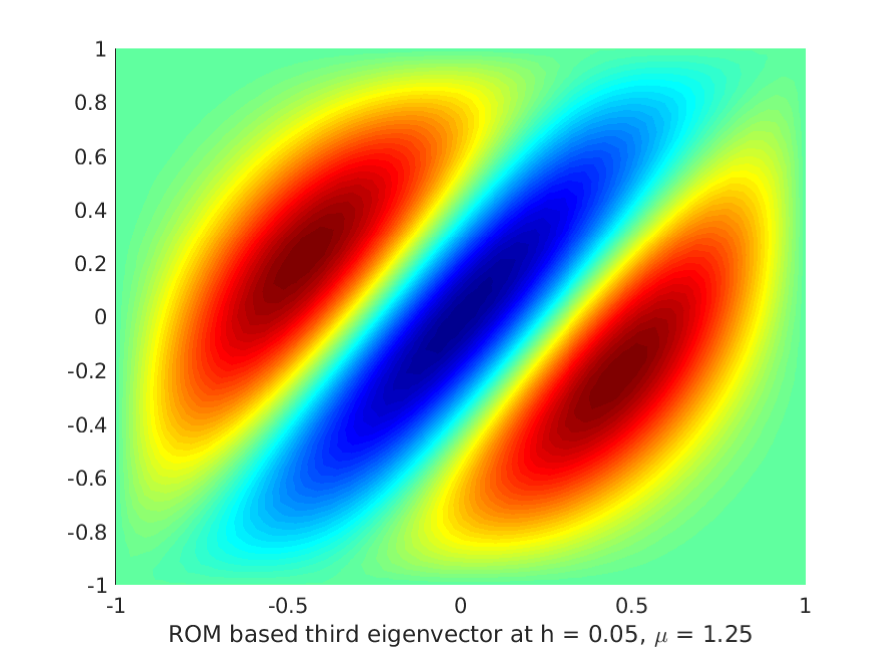}
\caption{Approximation of $u_3$ with snapshot based on $u_1, u_2$ and $u_3$ when ROM dimension is 34 and 80 respectively ($\mu = 1.25$ and $h = 0.05.$) }
\label{FIGURE:16}
\end{figure}

Let us consider the snapshot matrix consisting of $u_1$, $u_2$, and $u_3$ columnwise at the sample points and compute the first three smallest eigenvalues simultaneously using the reduced order method. Despite the intersection of the third and the fourth eigenvalues, the third eigenvalue of the ROM is converging to the third eigenvalue of the FEM in this case, see Figure~\ref{FIGURE:15}. In Table~\ref{TABLE:6} we have presented the first, second, and third eigenvalues of ROM at six sample points. It can be noted that we have obtained the first three smallest eigenvalues simultaneously. Eigenvalues obtained by the ROM are highly accurate. As expected, the third eigenvectors of the ROM are the same even when we increase the number of POD basis functions, see Figure~\ref{FIGURE:16}.



\begin{table}
 	 	\centering
 	\begin{tabular}{|c|c|c|c|c|c|c|c|} 
 		\hline

 	        {\begin{tabular}[c]{@{}c@{}} ROM dimension. \end{tabular}} &
		  {\begin{tabular}[c]{@{}c@{}} $\mu$ \end{tabular}} &
 		 {\begin{tabular}[c]{@{}c@{}} 3rd EV(FEM) \end{tabular}} & 
 		{\begin{tabular}[c]{@{}c@{}}  1st EV(ROM)  \end{tabular}}&
 		{\begin{tabular}[c]{@{}c@{}}  2nd EV(ROM)  \end{tabular}}&
 		{\begin{tabular}[c]{@{}c@{}}  3rd EV(ROM)  \end{tabular}}\\

  \hline

    32 &-1.25&12.38270171& 5.98379208 & 8.77547748 & 12.38270666\\
   &-0.75& 21.12292172&   7.00305355 & 12.86380188 & 21.12292235\\
   & -0.50& 22.22010384&7.23588577& 13.95750313			  &22.22010501		\\
   & 0.50&	22.22018744& 7.23586151 & 13.95737020	 &	22.22018859		\\
    &0.75&21.12370790&   7.00299337  & 12.86364179  & 21.12370853\\
    &1.25 & 	12.38958778	& 5.98368109 & 8.77690440  &12.38959316\\ 
   \hline
 
34 &-1.25&12.38270171&    5.98379248 & 8.77547319 & 12.38270251\\
    &-0.75&21.12292172&     7.00305409   &         12.86380239 & 21.12292251\\
    & -0.50& 22.22010384&    7.23588378	  &13.95750246 &22.22010682		\\
   & 0.50&	22.22018744&		7.23585920	 &	13.95736953	&22.22019046		\\
   &0.75&21.12370790&      7.00299381 &           12.86364223 & 21.12370869\\
   &1.25 & 	12.38958778	&   5.98368176  &      8.77689991 & 12.38958860\\ 
  \hline
   \end{tabular}
	\caption{Approximation of $\lambda_1, \lambda_2$ and $\lambda_3$ with snapshot based on $u_1, u_2$ and $u_3$ at $h = 0.05$}
	 	\label{TABLE:6}
 	\end{table}
 	
 \subsubsection{Results of the EVP considering  $u_1+u_2+u_3$  in the snapshot matrix}
\begin{figure}
\centering
\begin{subfigure}{0.34\textwidth}
\includegraphics[height=5cm,width=6cm]{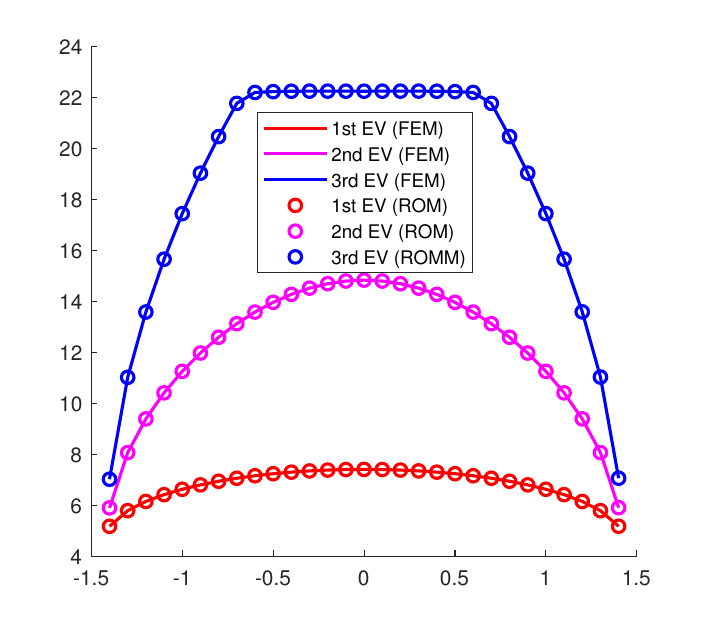}
 \caption{$u_1,u_2,u_3$}
 \end{subfigure}
 \begin{subfigure}{0.34\textwidth}
 \includegraphics[height=5cm,width=6cm]{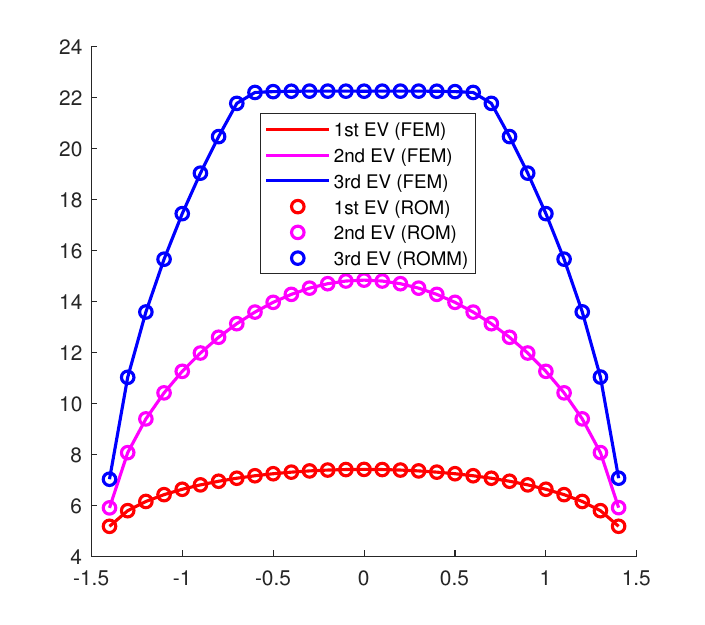}
 \caption{$u_1+u_2+u_3$}
 \end{subfigure}
 \caption{\rev{First three eigenvalues of FEM and ROM with snapshot based on $u_1,u_2,u_3$ and $u_1+u_2+u_3$}.}
     \label{FIGURE:17a}
\end{figure}

\begin{figure}[H]
\centering
\includegraphics[scale=0.4]{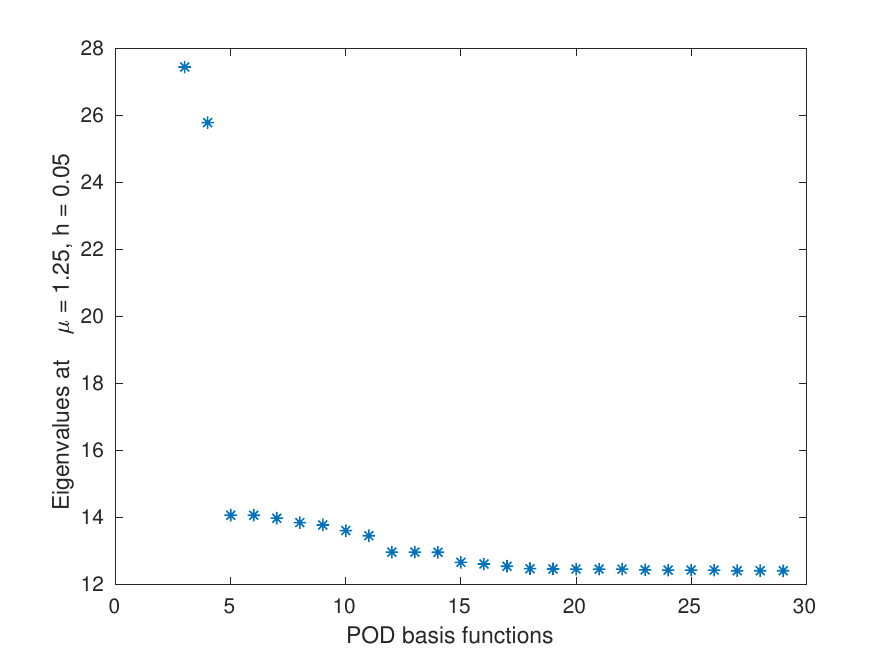}
\includegraphics[scale=0.4]{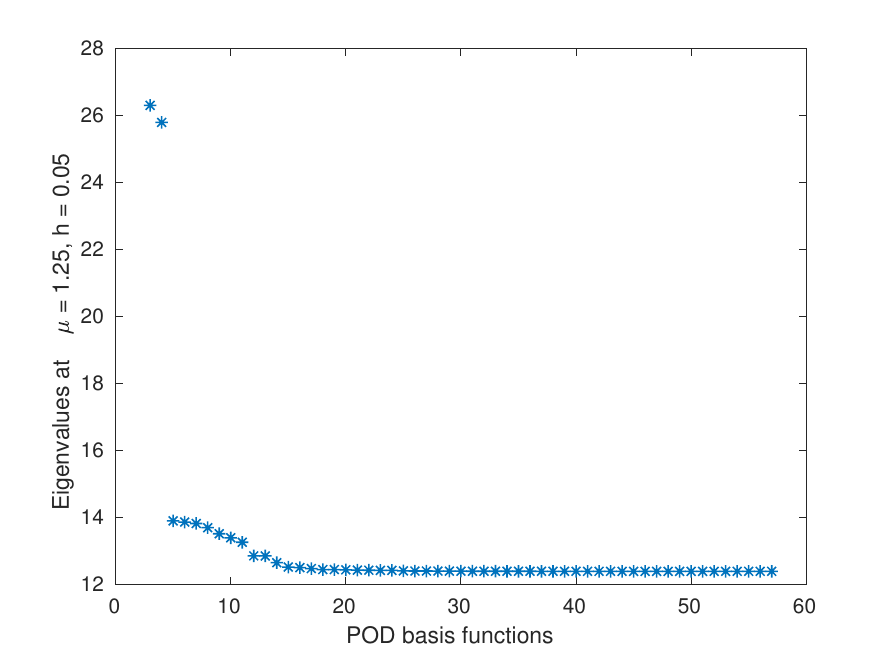}
\caption{Approximation of $\lambda_3$ with snapshot based on $u_1 + u_2 + u_3$: eigenvalues corresponding to different number of POD basis functions when $\mu = 1.25$ and $h = 0.05.$ }
\label{FIGURE:17}
\end{figure}

\begin{figure}[H]
\centering
\includegraphics[scale=0.4]{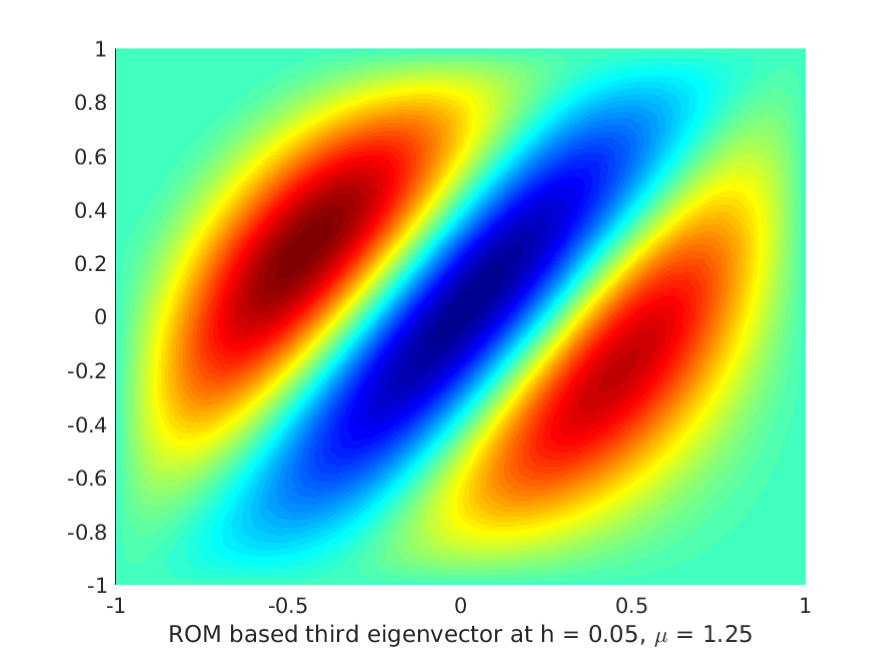}
\includegraphics[scale=0.4]{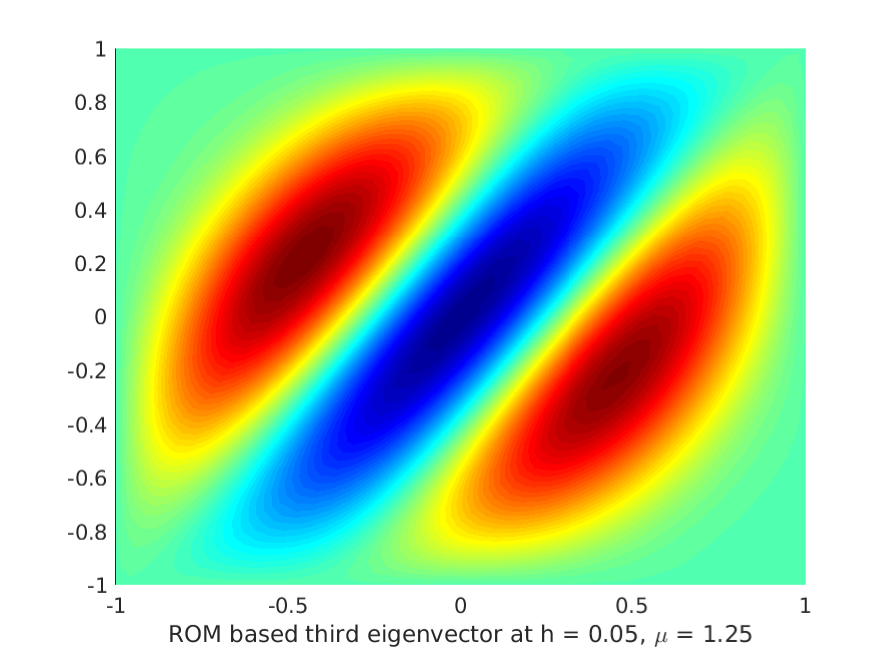}
\caption{Approximation of third eigenvector with snapshot based on $u_1 +u_2 + u_3$ when ROM dimension is 37 and 90 respectively at $\mu = 1.25$ and $h = 0.05.$}
\label{FIGURE:18}
\end{figure}

 We consider the snapshot matrix consisting $u_1 + u_2 + u_3$ column wise at the sample points and obtain first three smallest eigenvalues simultaneously using reduced order method.  In Figure~\ref{FIGURE:17} we have shown the plot of third eigenvalues of ROM. It is observed that third eigenvalues of ROM is converging to the third eigenvalue of FEM. In Table~\ref{TABLE:7} we have reported the first, second, and third eigenvalues of ROM at six sample points.
 \rev{We have shown only the third eigenvalue of FEM but all the first three eigenvalues of ROM converge to the respective eigenvalue of FEM, which is evident from Figure~\ref{FIGURE:17a}.}
 As expected, the third eigenvectors of ROM are the same despite different number of POD basis functions, see Figure~\ref{FIGURE:18}.

\begin{table}[H]
\centering
\begin{tabular}{|c|c|c|c|c|c|c|c|} 
\hline 
{\begin{tabular}[c]{@{}c@{}} ROM dimension \end{tabular}} &
{\begin{tabular}[c]{@{}c@{}} $\mu$ \end{tabular}} &
{\begin{tabular}[c]{@{}c@{}} 3rd EV(FEM) \end{tabular}} & 
{\begin{tabular}[c]{@{}c@{}}  1st EV(ROM)  \end{tabular}}&
{\begin{tabular}[c]{@{}c@{}}  2nd EV(ROM)  \end{tabular}}&
{\begin{tabular}[c]{@{}c@{}}  3rd EV(ROM)  \end{tabular}}\\
\hline
22 &-1.25 &12.38270171& 5.98859162& 8.77898776	& 12.38307378 \\
   &-0.75 &21.12292172&7.00306584 & 12.86381909 &21.12293201\\
   &-0.50& 22.22010384&7.23591373 &13.95751305&22.22012955		\\
   &0.50 &22.22018744& 7.23586309 &13.95747389&22.22026978	\\
    &0.75&21.12370790& 7.00308513 & 12.86383468 & 21.12375349\\
    &1.25 & 	12.38958778	& 5.99388504 & 8.86730411 &	12.45479398	\\ 
   \hline
25 &-1.25&12.38270171&  5.98673888& 8.77822960& 12.38272369 \\
    &-0.75&21.12292172&	7.00305515  & 12.86380311 & 21.12292566\\
    & -0.50& 22.22010384&7.23589141 &13.95750303	 &22.22011036\\
   & 0.50&	22.22018744&	7.23586289&13.95738837&22.22019762	\\
    &0.75&21.12370790& 7.00309327 & 12.86379347 &21.12372327\\
    &1.25 & 	12.38958778	&5.99440793 &8.79856593&12.40055649	\\ 
  \hline
   \end{tabular}
   \caption{Approximation of $\lambda_1, \lambda_2$ and $\lambda_3$ with snapshot based on $u_1 + u_2 + u_3$ at $h = 0.05$}
	 	\label{TABLE:7}
 	\end{table}


\begin{table}[H]
\centering
\begin{tabular}{|c|c|c|c|c|c|c|c|} 
\hline 
{\begin{tabular}[c]{@{}c@{}} $(c_1,c_2,c_3)$ \end{tabular}} &
{\begin{tabular}[c]{@{}c@{}} $\mu$ \end{tabular}} & 
{\begin{tabular}[c]{@{}c@{}} Rel. Error in 1st EV  \end{tabular}}&
{\begin{tabular}[c]{@{}c@{}} Rel. Error in 2nd EV  \end{tabular}}&
{\begin{tabular}[c]{@{}c@{}} Rel. Error in 3rd EV  \end{tabular}}\\
\hline
$(1,5,10)$ &-1.25 &$1.5 \times 10^{-3}$  &$2.3 \times 10^{-5}$  &$3.7 \times 10^{-6}$ \\ 
 &-.75 &$4.5 \times 10^{-5}$  &$4.5 \times 10^{-6}$  &$8.1 \times 10^{-8}$ \\ 
 &-0.50 &$3.3 \times 10^{-5}$  &$1.5 \times 10^{-6}$  &$1.6 \times 10^{-8}$ \\ 
 &0.50 &$1.2 \times 10^{-5}$  &$1.8 \times 10^{-7}$  &$6.8 \times 10^{-8}$ \\ 
 &0.75 &$3.7 \times 10^{-5}$  &$2.4 \times 10^{-7}$  &$1.3 \times 10^{-8}$ \\ 
  &1.25 &$2.6 \times 10^{-3}$  &$5.0 \times 10^{-5}$  &$1.3 \times 10^{-5}$ \\ 
  \hline
  $(5,10,1)$ &-1.25 &$4.6 \times 10^{-5}$  &$7.5 \times 10^{-6}$  &$8.3 \times 10^{-5}$ \\ 
 &-.75 &$2.2 \times 10^{-7}$  &$1.3\times 10^{-7}$  &$1.5\times 10^{-6}$ \\ 
 &-0.50 &$2.0 \times 10^{-7}$  &$3.2 \times 10^{-8}$  &$1.4 \times 10^{-6}$ \\ 
 &0.50 &$5.7\times 10^{-7}$  &$1.3 \times 10^{-7}$  &$1.0 \times 10^{-6}$ \\ 
 &0.75 &$3.1 \times 10^{-8}$  &$2.6 \times 10^{-8}$  &$3.1 \times 10^{-7}$ \\ 
  &1.25 &$4.1 \times 10^{-4}$  &$2.9 \times 10^{-6}$  &$4.5\times 10^{-3}$ \\ 
  \hline
  $(10,1,5)$ &-1.25 &$2.4 \times 10^{-5}$  &$4.8 \times 10^{-4}$  &$9.3 \times 10^{-6}$ \\ 
 &-.75 &$2.8 \times 10^{-7}$  &$5.4\times 10^{-5}$  &$1.6\times 10^{-7}$ \\ 
 &-0.50 &$1.9 \times 10^{-7}$  &$9.1 \times 10^{-6}$  &$1.2 \times 10^{-7}$ \\ 
 &0.50 &$1.5\times 10^{-7}$  &$1.3 \times 10^{-5}$  &$3.2 \times 10^{-7}$ \\ 
 &0.75 &$4.5 \times 10^{-8}$  &$7.5 \times 10^{-6}$  &$1.2 \times 10^{-7}$ \\ 
  &1.25 &$4.4 \times 10^{-4}$  &$2.2 \times 10^{-4}$  &$7.0\times 10^{-4}$ \\ 
  \hline
   \end{tabular}
   \caption{\rev{Relative error for the first three ROM eigenvalues computed with snapshot based on different linear combinations of $u_1,u_2,u_3$.}}
	 	\label{TABLE:7a}
 	\end{table}

\begin{table}
 \centering
 \begin{tabular}{|c|c|c|c|c|c|c|c|} 
 \hline
 {\begin{tabular}[c]{@{}c@{}} $\mu$ \end{tabular}} &
 {\begin{tabular}[c]{@{}c@{}} Relative error\\($u_1, u_2, u_3$) \end{tabular}} &
 {\begin{tabular}[c]{@{}c@{}}  Relative error\\($u_1 + u_2 + u_3$)  \end{tabular}}\\
 \hline
-1.25&   $4.0 \times 10^{-7}$   & $3.0 \times 10^{-5}$\\
 -0.75&  $3.0 \times 10^{-8}$   & $4.8 \times 10^{-7}$ \\
 -0.50&  $5.2 \times 10^{-8}$   & $1.2 \times 10^{-6}$ \\
  0.50&  $5.2 \times 10^{-8}$   & $3.7 \times 10^{-6}$ \\
 0.75&   $3.0 \times 10^{-8}$   & $2.2 \times 10^{-6}$\\
 1.25 &  $4.0 \times 10^{-7}$   & $5.2 \times 10^{-3}$ \\ 
  \hline
 -1.25& $6.4 \times 10^{-8}$ 	& $1.8 \times 10^{-6}$ \\
 -0.75& $3.7 \times 10^{-8}$     &  $1.9 \times 10^{-7}$ \\
 -0.50& $1.3 \times 10^{-7}$  	& $2.9 \times 10^{-7}$ \\
 0.50 & $1.4 \times 10^{-7}$   	& $4.6 \times 10^{-7}$ \\
 0.75 & $3.7 \times 10^{-8}$	    & $7.3 \times 10^{-7}$   \\
 1.25 & $6.5 \times 10^{-8}$ 	&  $8.8 \times 10^{-4}$   \\ 
  \hline
   \end{tabular}
	\caption{Relative error between FEM and ROM based third eigenvalue  with snapshot based on $u_1, u_2, u_3$ and $u_1 + u_2 + u_3$ at $h = 0.05$}
	 \label{TABLE:8}
 	\end{table}
From Table~\ref{TABLE:8} it is clear that relative error in the case when we considered $u_1$, $u_2$, and $u_3$ in the snapshot matrix, provides us with slightly better results than when considering $u_1 + u_2 + u_3$. However, the latter case is computationally cheaper.
\rev{As we mentioned in \eqref{com_snap}, if we take any non-zero linear combination of $u_1,u_2,\dots,u_k$ as snapshots then the first $n$ eigenvalues of ROM converge to the corresponding eigenvalue of the FEM. If we take different weights for the different eigenvectors then again the eigenvalues of ROM will converge to the corresponding eigenvalue of FEM, but the accuracy will not be the same for all eigenvalues. The most accurate eigenvalue will be that eigenvalue whose weight is maximum in the linear combination. In Table~\ref{TABLE:7a} we have reported the relative error for the first three eigenvalues while different combinations of $u_1,u_2,u_3$ are used as snapshot. Taking equal weights for all the eigenvectors gives the same results, no matter what is the value of the weight, since the corresponding snapshot matrices are multiple to each other. So the best choice is to take the sum or the average of all the eigenvectors $u_1,u_2,\dots,u_n$ in order to get eigenvalues converging in the right order.}

\subsection{Reduced order method to obtain $\lambda_4$}
In this subsection we discuss numerical results for the fourth eigenvalue considering different combinations of eigenvectors in the snapshot matrix. Since the third and the fourth eigenvalues are intersecting, the numerical results obtained in this subsection present several interesting features.
\subsubsection{Results of the EVP considering only $u_4$ is in the snapshot matrix}
\begin{figure}
\centering
\includegraphics[scale=0.35]{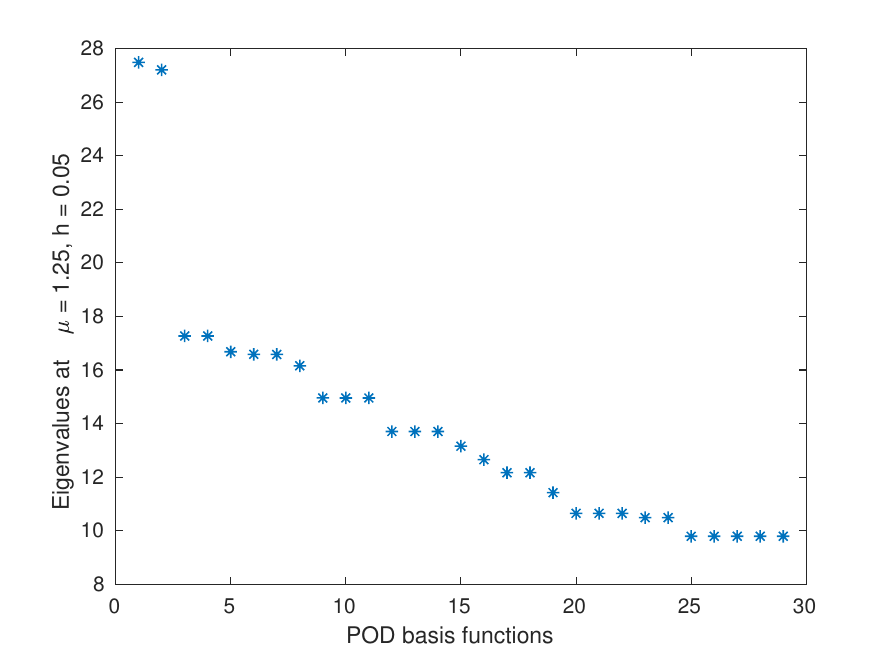}
\includegraphics[scale=0.35]{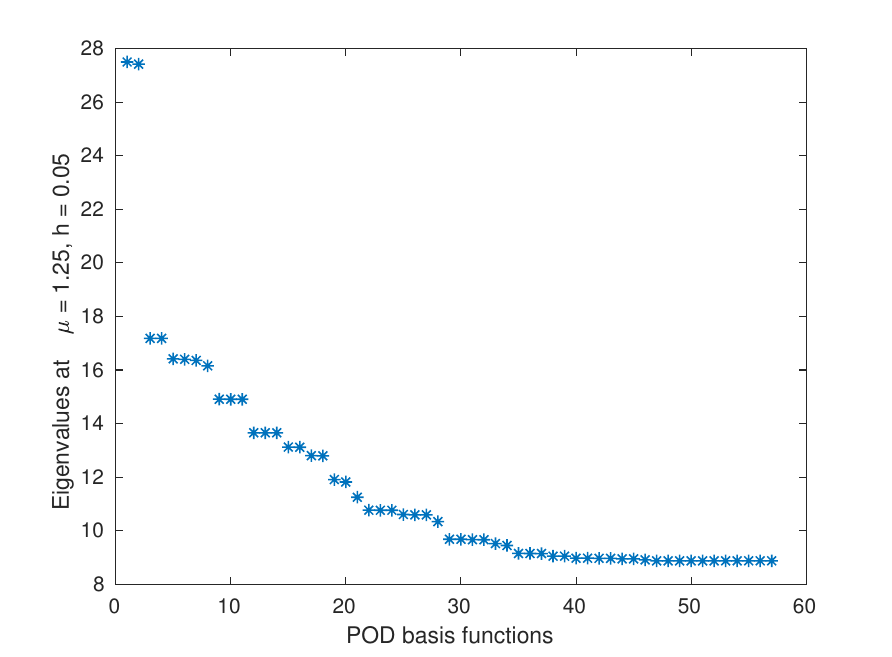}
\caption{Approximation of $\lambda_4$ with snapshot based on $u_4$: eigenvalues corresponding to different number of POD basis functions when $\mu = 1.25$ and $h = 0.05.$ }
\label{FIGURE:19}
\end{figure}

\begin{figure}
     \begin{subfigure}{0.32\textwidth}
         \includegraphics[height=4.5cm,width=5.5cm]{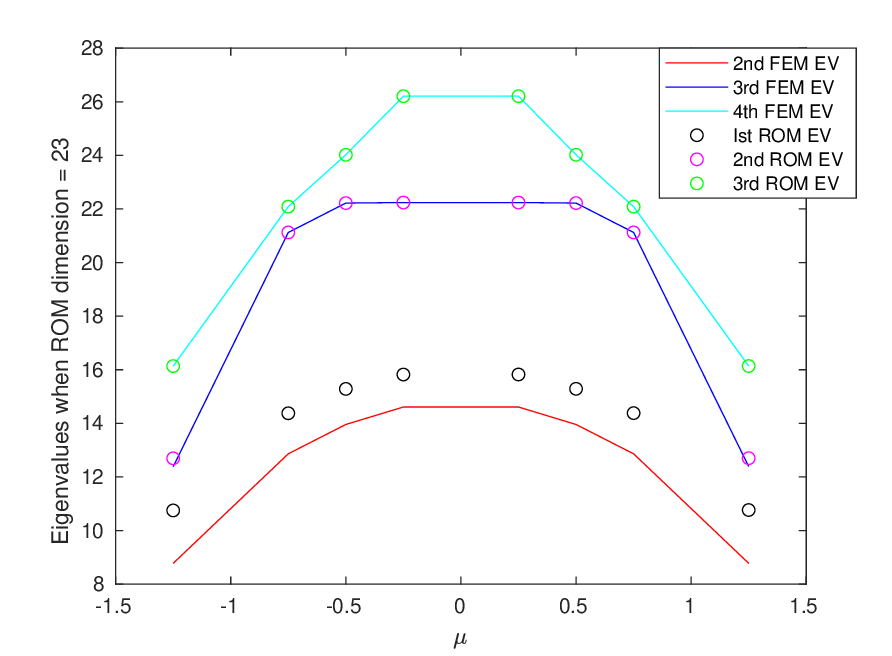}
     \end{subfigure}
     \begin{subfigure}{0.32\textwidth}
         \includegraphics[height=4.5cm,width=5.5cm]{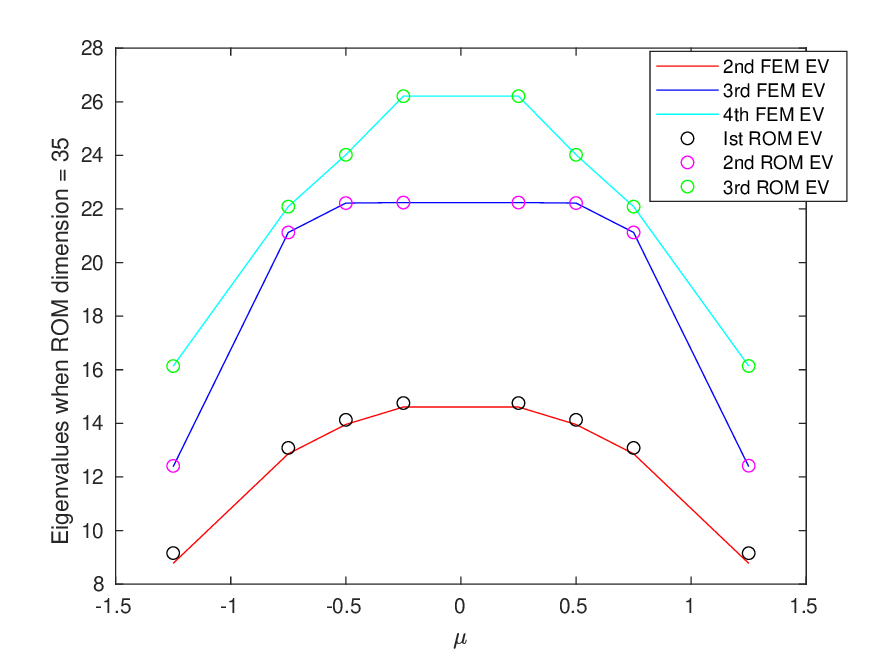}
     \end{subfigure}
     \begin{subfigure}{0.32\textwidth}
         \includegraphics[height=4.5cm,width=5.5cm]{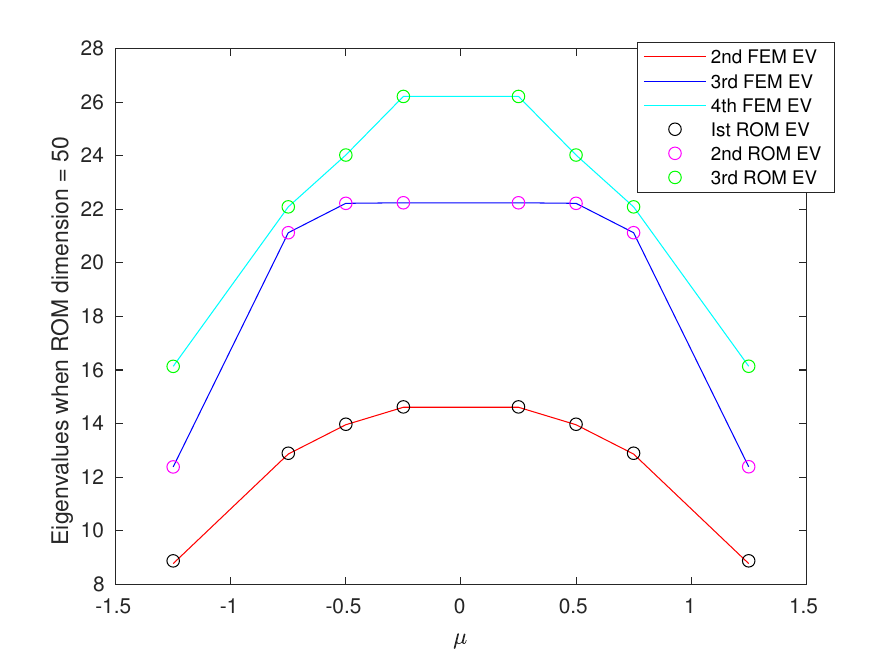}
     \end{subfigure}
     \caption{FEM and ROM based eigenvalues at different ROM dimensions}
     \label{FIGURE:20}
\end{figure}

Let us consider the snapshot matrix consisting only of the fourth eigenvector $u_4$ at the sample points and compute the eigenvalues using the reduced order method. In Table~\ref{TABLE:9} the fourth eigenvalue of the FEM and the first eigenvalue of the ROM are reported at few test points. It can be seen that the fourth eigenvalue of the FEM is not matching the first eigenvalue of the ROM, which is also evident from Figure~\ref{FIGURE:19}. Moreover, when we increase the number of POD basis functions, we notice that the first eigenvalue of the ROM matches with the second eigenvalues of the FEM, the second eigenvalue of the ROM matches with the third eigenvalue of the FEM, and the third eigenvalues of the ROM matches with the fourth eigenvalues of the FEM, see Figure~\ref{FIGURE:20}.


\begin{table}
\centering
\begin{tabular}{|c|c|c|c|c|c|c|c|} 
\hline
{\begin{tabular}[c]{@{}c@{}} ROM dimension \end{tabular}} &
{\begin{tabular}[c]{@{}c@{}} $\mu$ \end{tabular}} &
{\begin{tabular}[c]{@{}c@{}} 4th EV(FEM) \end{tabular}} & 
{\begin{tabular}[c]{@{}c@{}}  1st EV(ROM)  \end{tabular}}\\
\hline
 20 &-1.25&16.13593156		& 10.63650138			 \\
    &-0.75& 22.08781810& 14.31298332\\
    &-0.50& 24.02185801		& 15.23531319	\\
    &0.50& 24.02234525& 15.24102596\\
    &0.75&22.08790319	&   14.32162342 \\
    &1.25 & 16.13984606& 10.65123417	\\ 
   \hline
23 &-1.25&16.13593156		& 10.74968666		 \\
   &-0.75& 22.08781810&  	14.37615778 \\
   &-0.50& 24.02185801		& 15.28535039	\\
   &0.50& 24.02234525& 15.28887728	\\
   &0.75&22.08790319	&   14.38179829 \\
   &1.25 & 16.13984606& 10.76563540\\ 
   \hline
   \end{tabular}
   \caption{Approximation of $\lambda_4$ with snapshot based on $u_4$ at $h = 0.05$}
\label{TABLE:9}
\end{table}

\subsubsection{Results of the EVP considering $u_3$ and $u_4$ in the snapshot matrix}

\begin{figure}
\centering
\includegraphics[scale=0.4]{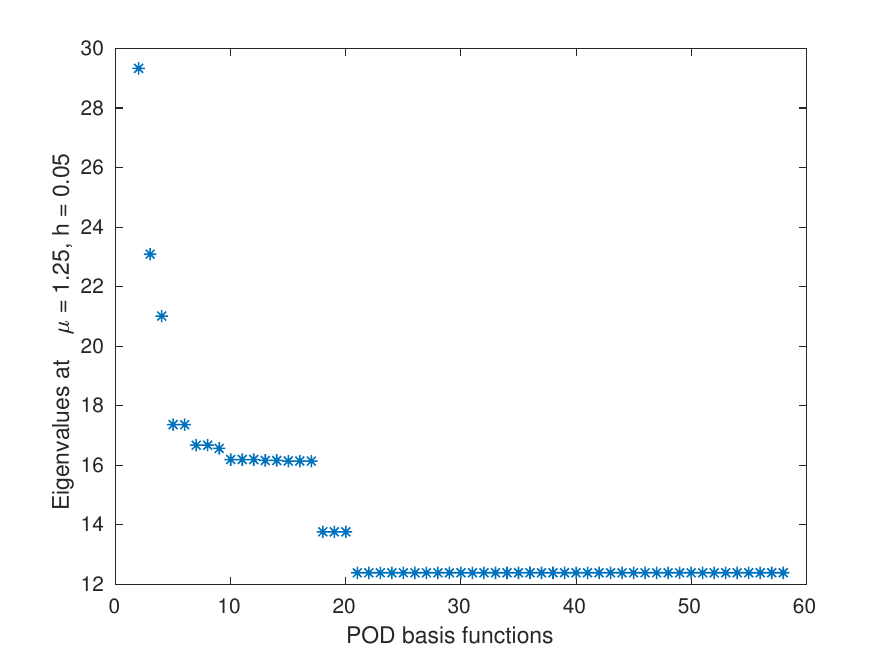}
\includegraphics[scale=0.4]{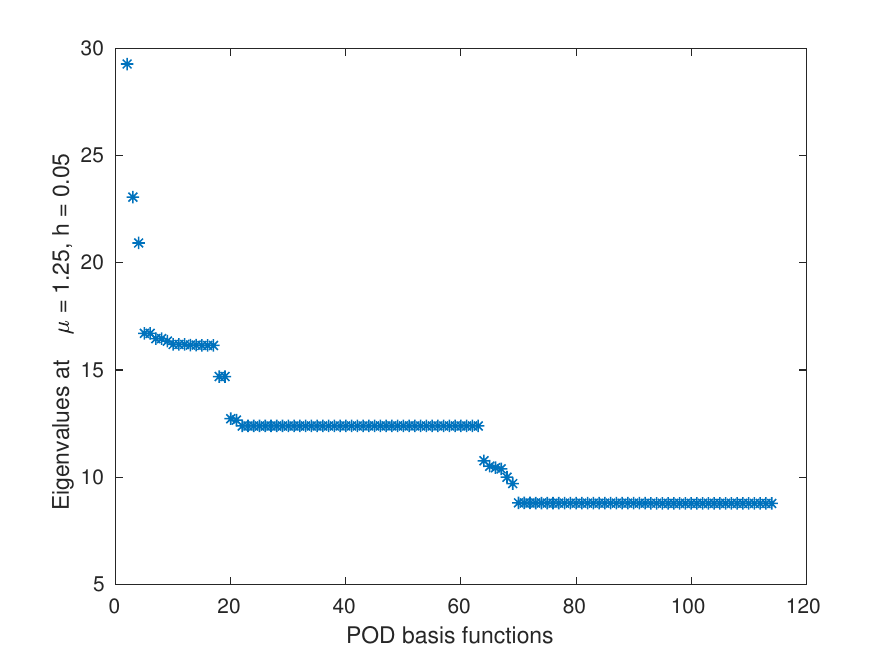}
\caption{Approximation of $\lambda_4$ with snapshot based on $u_3$ and $u_4$: eigenvalues corresponding to different number of POD basis functions when $\mu = 1.25$ and $h = 0.05.$  }
\label{FIGURE:22}
\end{figure}

Since $\lambda_3$ and $\lambda_4$ are intersecting, it is natural to consider both eigenvectors $u_3$ and $u_4$ columnwise in the snapshot matrix. In Table~\ref{TABLE:10} we have reported the first and second eigenvalues of the ROM and the fourth eigenvalue of the FEM. As we can see, the second eigenvalue of the ROM matches with the fourth eigenvalues of the FEM. However, as we increase the number of POD basis functions, the second eigenvalue of the ROM is not converging to the fourth eigenvalue of the FEM, see Figure~\ref{FIGURE:22}.

\begin{table}
 	 	\centering
 	\begin{tabular}{|c|c|c|c|c|c|c|c|} 
 		\hline

 	        {\begin{tabular}[c]{@{}c@{}} ROM dimension \end{tabular}} &
		  {\begin{tabular}[c]{@{}c@{}} $\mu$ \end{tabular}} &
 		 {\begin{tabular}[c]{@{}c@{}} 4th EV(FEM) \end{tabular}} & 
 		{\begin{tabular}[c]{@{}c@{}}  1st EV(ROM)  \end{tabular}}&
 		{\begin{tabular}[c]{@{}c@{}}  2nd EV(ROM)  \end{tabular}}\\
 	
  \hline

    26 &-1.25 & 16.13593156		& 10.32768230 & 12.38271820								 \\
   &-0.75 & 22.08781810 & 14.01748876 & 21.12292540\\
  &-0.50 & 24.02185801 & 14.96276873 & 22.22010448	\\
   & 0.50& 	24.02234525 & 	14.96700836 & 22.22018808\\
   &0.75& 22.08790319	&   14.02396381 & 21.12371133\\
   &1.25 & 16.13984606& 	10.33922731	& 12.38960275		\\ 
   \hline

    30 &-1.25&16.13593156		& 10.35774089& 	12.38270198					 \\
  &-0.75& 22.08781810&  14.02551777& 21.12292203\\
  &-0.50 & 24.02185801 & 14.96893932	 &22.22010447	\\
  & 0.50& 	24.02234525 & 	14.97067172 & 22.22018807\\
   &0.75&22.08790319	&   	14.02852884&	21.12370817 \\
   &1.25 & 16.13984606& 10.36928581&	12.38958799		\\ 
   \hline
   \end{tabular}
   \caption{Approximation of $\lambda_3$ and $\lambda_4$ with snapshot based on $u_3$ and $u_4$ at $h = 0.05$}
	 	\label{TABLE:10}
 	\end{table}

\subsubsection{Results of the EVP considering  $u_1$, $u_2$, $u_3$, and $u_4$ in the snapshot matrix}
\begin{figure}
\centering
\includegraphics[scale=0.4]{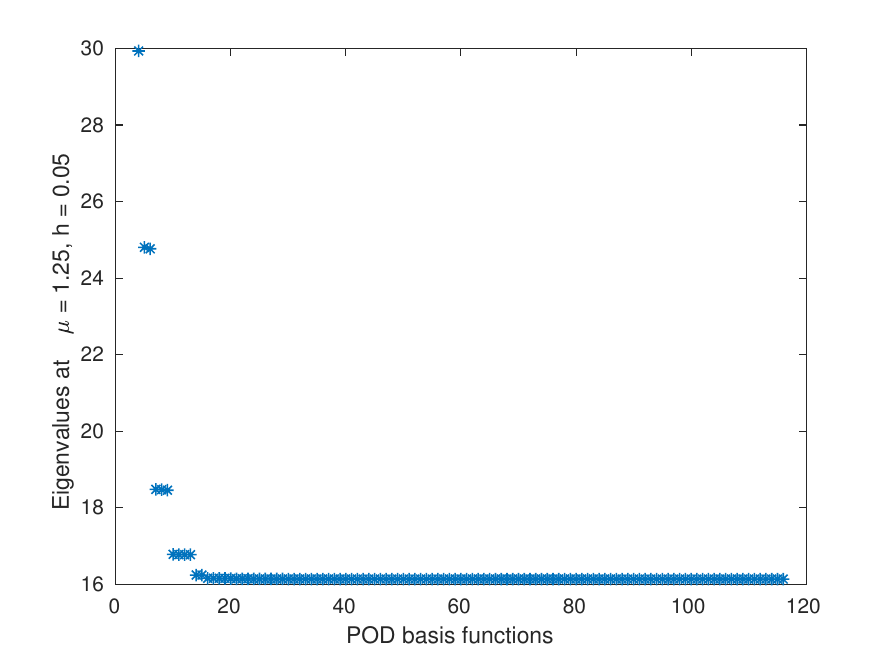}
\includegraphics[scale=0.4]{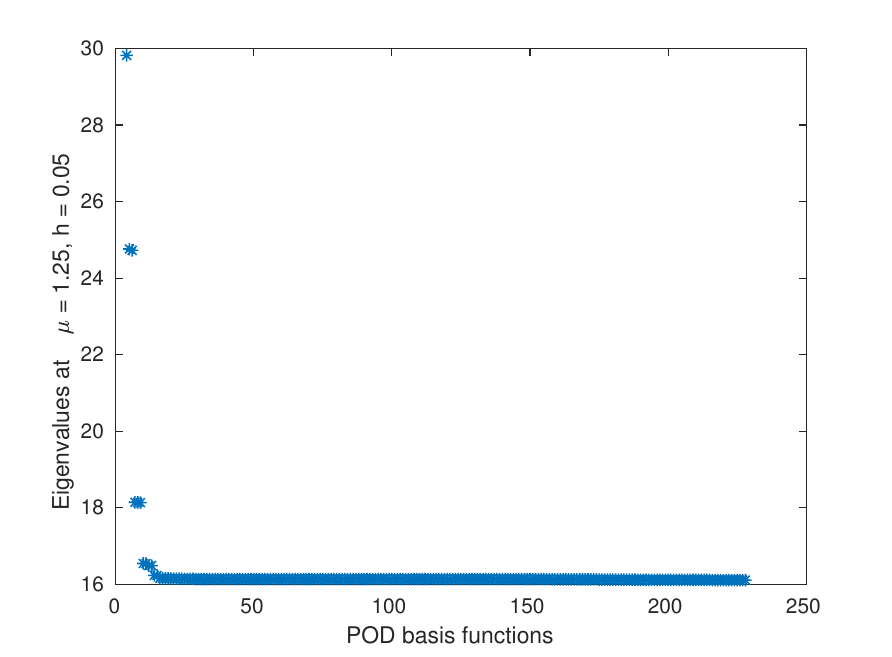}
\caption{Approximation of $\lambda_4$ with snapshot based on $u_1, u_2, u_3$ and $u_4$: eigenvalues corresponding to different number of POD basis functions when $\mu = 1.25$ and $h = 0.05.$  }
\label{FIGURE:24}
\end{figure}

\begin{figure}
\centering
\includegraphics[scale=0.5]{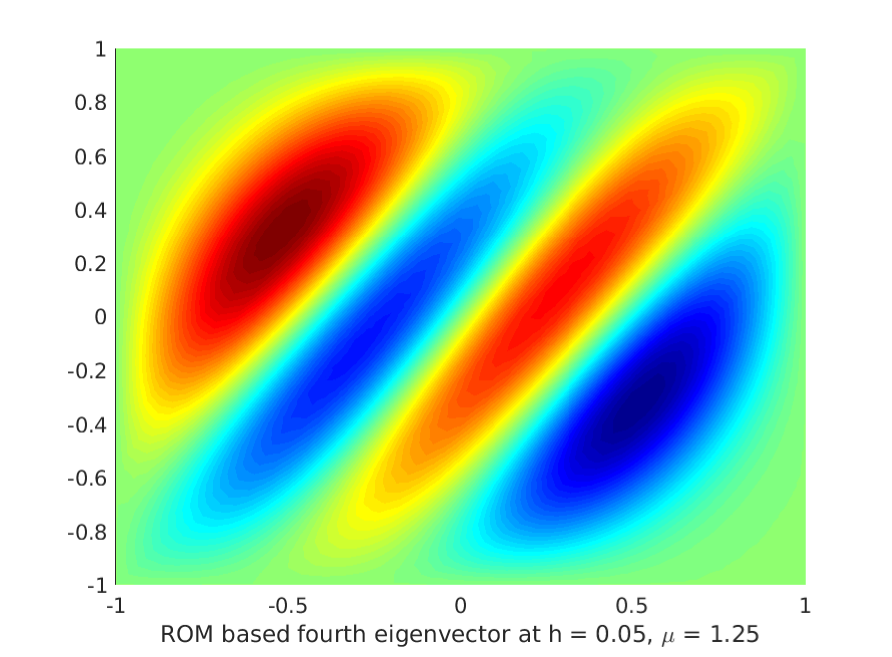}
\includegraphics[scale=0.5]{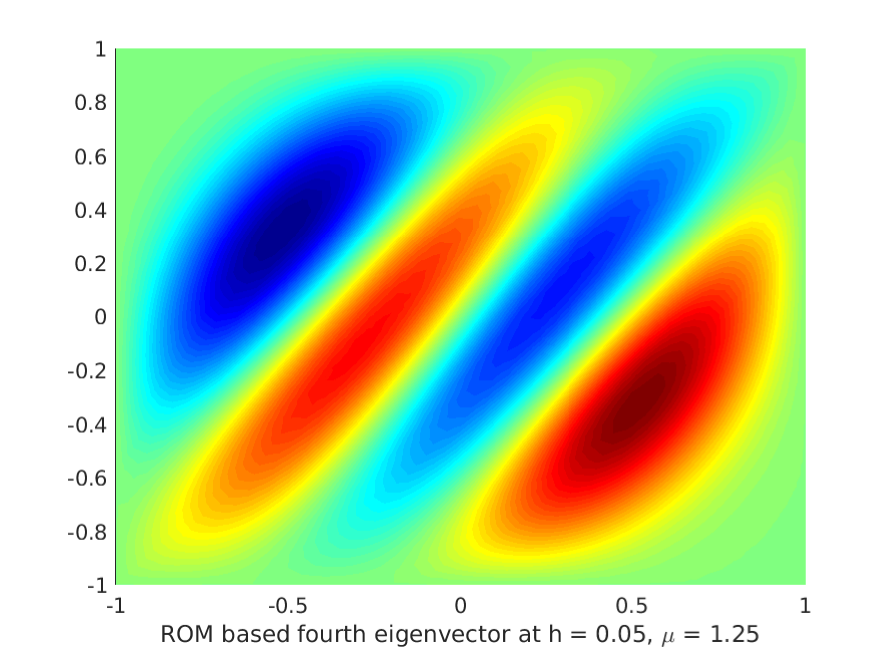}
\caption{Approximation of $u_4$ with snapshot based on $u_1, u_2, u_3$ and $u_4$ when ROM dimension is 44 and 90 respectively ($\mu = 1.25$ and $h = 0.05.$ }
\label{FIGURE:25}
\end{figure}

Let us consider $u_1$, $u_2$, $u_3$, and $u_4$ columnwise in the snapshot matrix and compute the first four smallest eigenvalues simultaneously using the reduced order method. In Figure~\ref{FIGURE:24} we have presented the plot of the fourth eigenvalue of the ROM. It can be observed that the fourth eigenvalue of the ROM is converging to the fourth eigenvalue of the FEM. In Table~\ref{TABLE:11} the first four smallest eigenvalues obtained by the ROM and the fourth eigenvalue of the FEM are reported. The obtained numerical results are highly accurate. As expected, the fourth eigenvector of ROM remains the same at different number of POD basis functions, see Figure~\ref{FIGURE:25}.

\begin{table}
 	 	\centering
 	\begin{tabular}{|c|c|c|c|c|c|c|c|} 
 		\hline

 	        {\begin{tabular}[c]{@{}c@{}} ROM dim. \end{tabular}} &
		  {\begin{tabular}[c]{@{}c@{}} $\mu$ \end{tabular}} &
 		 {\begin{tabular}[c]{@{}c@{}} 4th EV(FEM) \end{tabular}} & 
 		{\begin{tabular}[c]{@{}c@{}}  1st EV(ROM)  \end{tabular}}&
 		{\begin{tabular}[c]{@{}c@{}}  2nd EV(ROM)  \end{tabular}}&
 		{\begin{tabular}[c]{@{}c@{}}  3rd EV(ROM)  \end{tabular}}&
 		{\begin{tabular}[c]{@{}c@{}}  4th EV(ROM)  \end{tabular}}\\
 	
\hline

    41 &-1.25&16.13593156& 5.98379230		 & 8.77547626		 & 12.38271182&	16.13594735			\\
   &-0.75& 22.08781810&  	7.00305433 & 	12.86380238 & 	21.12292388& 22.08782059 \\
    &-0.50 & 24.02185801 & 7.23588570	 &13.95750374		 & 22.22010425	&24.02185908			\\
   & 0.50& 	24.02234525 & 		7.23586132 & 13.95737083 & 22.22018787 & 24.02234624	 \\
    &0.75&22.08790319&  	7.00299407 & 12.86364225 &	21.12370996& 	22.08790566\\
   &1.25 & 16.13984606	& 5.98368138 & 	8.77690306 &12.38959840& 16.13986166	\\ 
   \hline
 
44 &-1.25&16.13593156& 5.98379293 & 8.77547307	& 12.38270419 & 16.13593241	\\
  &-0.75& 22.08781810&  	7.00305414 & 12.86380198 & 21.12292253 & 22.08781883\\
   &-0.50 & 24.02185801 & 7.23588682			 &13.95750267 & 22.22010421		&24.02185980			\\
   & 0.50& 	24.02234525 & 	7.23586277 & 	13.95736973 & 22.22018781 &	24.02234706	 \\
  &0.75&22.08790319&  	7.00299410 & 12.86364189 & 21.12370888 & 	22.08790390\\
  &1.25 & 16.13984606	& 	5.98368190	 & 8.77689982	 &	12.38959032 & 	16.13984691	\\ 
   \hline
   \end{tabular}
   \caption{Approximation of $\lambda_1, \lambda_2, \lambda_3$ and $\lambda_4$ with snapshot based on $u_1, u_2, u_3$ and $u_4$ at $h = 0.05$}
	 	\label{TABLE:11}
 	\end{table}
 	
\subsubsection{Results of the EVP considering $u_1 + u_2 + u_3 + u_4$ in the snapshot matrix}
\begin{figure}[H]
\centering
\includegraphics[scale=0.4]{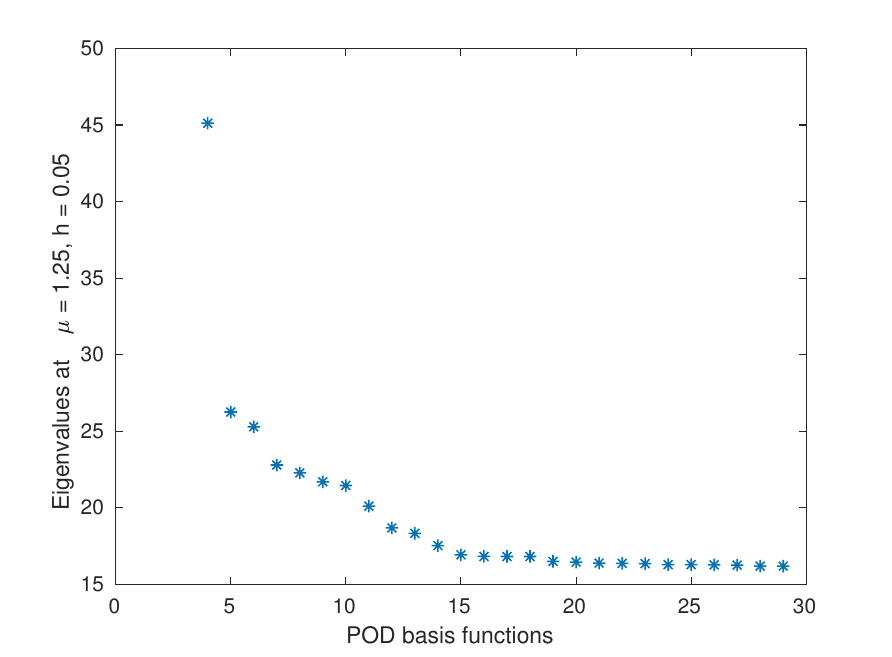}
\includegraphics[scale=0.4]{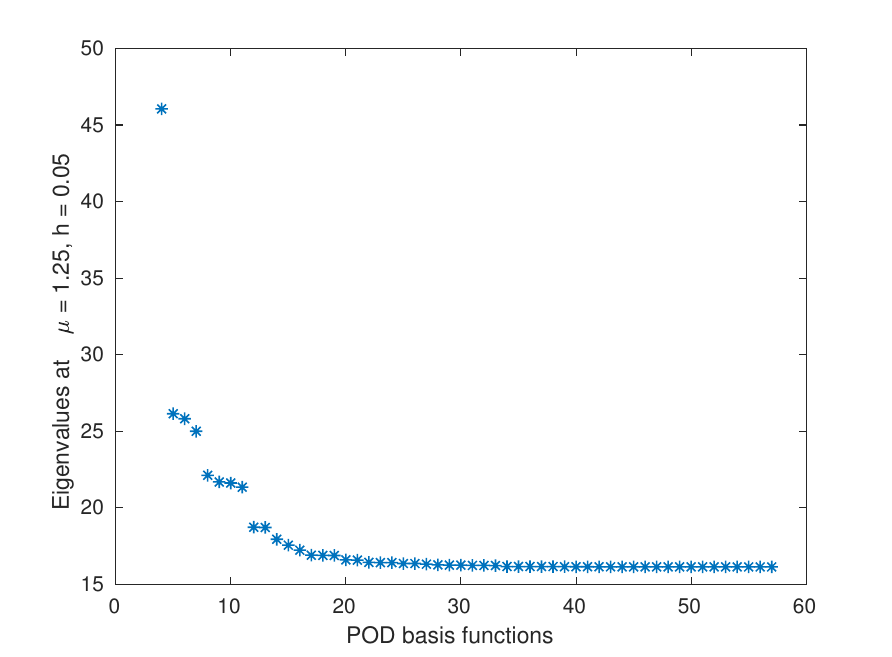}
\caption{Approximation of $\lambda_4$ with snapshot based on $u_1 + u_2 + u_3 + u_4$: eigenvalues corresponding to different number of POD basis functions when $\mu = 1.25$ and $h = 0.05.$  }
\label{FIGURE:26}
\end{figure}

\begin{figure}
\centering
\includegraphics[scale=0.4]{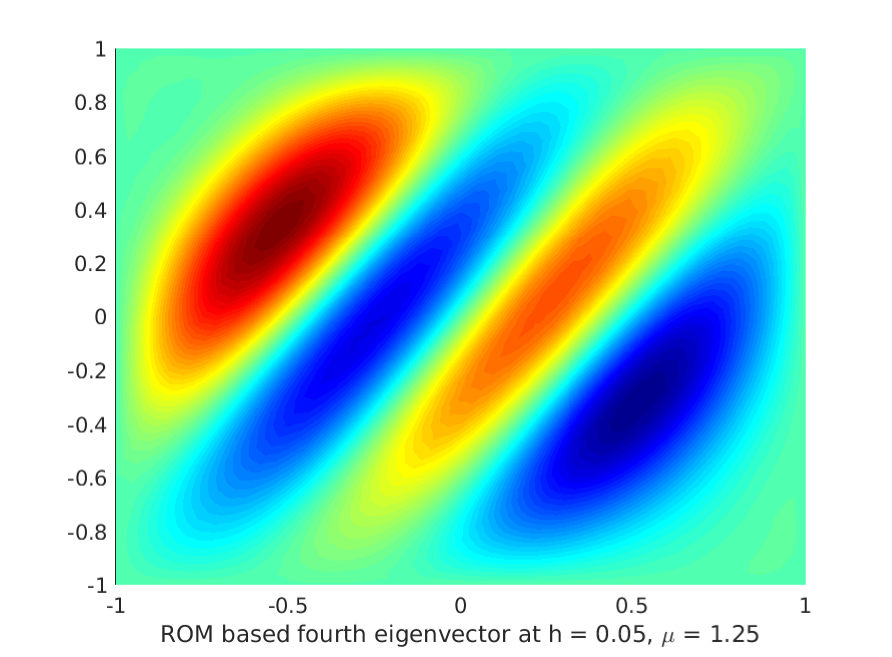}
\includegraphics[scale=0.4]{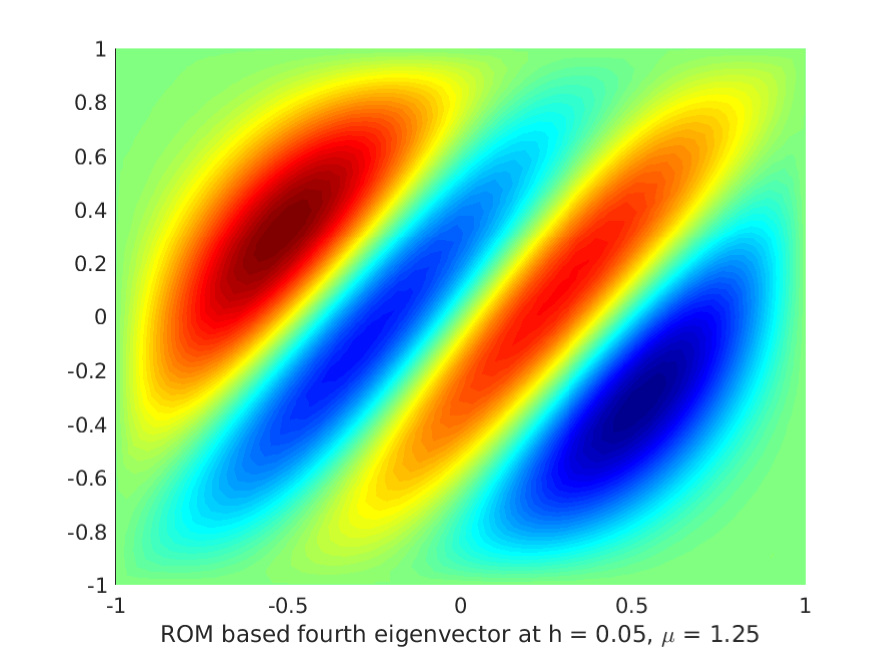}
\caption{Approximation of $u_4$ with snapshot based on $u_1 +  u_2 + u_3 + u_4$ when ROM dimension is 30 and 50 respectively ($\mu = 1.25$ and $h = 0.05.$}
\label{FIGURE:27}
\end{figure}

Let us consider the snapshot matrix consisting of the combination $u_1 + u_2 + u_3 + u_4$ columnwise at the sample points and compute the first four smallest eigenvalues simultaneously using the ROM. In Figure~\ref{FIGURE:26} we have shown the plot of the fourth eigenvalue of the ROM. It is observed that the fourth eigenvalue of the ROM is converging to the fourth eigenvalue of the FEM even when the number of POD basis functions increases. In Table~\ref{TABLE:12} we have reported the first four smallest eigenvalues of ROM at six sample points.  As expected, the fourth eigenvector of the ROM remains the same despite different number of POD basis functions, see Figure~\ref{FIGURE:27}.

\begin{table}[H]
\centering
\begin{tabular}{|c|c|c|c|c|c|c|c|} 
 \hline
 {\begin{tabular}[c]{@{}c@{}} ROM dim. \end{tabular}} &
 {\begin{tabular}[c]{@{}c@{}} $\mu$ \end{tabular}} &
 {\begin{tabular}[c]{@{}c@{}} 4th EV(FEM) \end{tabular}} & 
 {\begin{tabular}[c]{@{}c@{}}  1st EV(ROM)  \end{tabular}}&
 {\begin{tabular}[c]{@{}c@{}}  2nd EV(ROM)  \end{tabular}}&
 {\begin{tabular}[c]{@{}c@{}}  3rd EV(ROM)  \end{tabular}}&
 {\begin{tabular}[c]{@{}c@{}}  4th EV(ROM)  \end{tabular}}\\
\hline
25 &-1.25&16.13593156 &6.00975075 & 8.79767218 & 12.50952660 &16.23425928	\\
   &-0.75&22.08781810 &7.00397060 & 12.86485702& 21.12369222 & 22.08841968 \\
  &-0.50 &24.02185801 &7.23599290	&13.95753367 & 22.22012647	&24.02187149 \\		
   & 0.50& 	24.02234525 & 		7.23652419 & 	13.95782932 & 22.22070774 &	24.02277614 \\
    &0.75&22.08790319&  		7.00543420& 12.86604820 &		21.12640999 & 		22.09043982\\
   &1.25 & 16.13984606	& 		6.01170085 & 			8.81168198 &		12.56293955	 & 16.28632874		\\ 
   \hline
 
28 &-1.25&16.13593156& 5.99264622 &8.78419766	& 12.51450348						 & 16.24263319						\\
  &-0.75& 22.08781810&  	7.00338482 & 12.86406343	& 21.12305747 & 22.08797149\\
   &-0.50 & 24.02185801 & 	7.23594905	 &	13.95750914		 & 	22.22013876	&	24.02186464			\\
   & 0.50& 	24.02234525 & 			7.23637667 & 13.95781101 & 22.22063923 & 	24.02275865	 \\
  &0.75&22.08790319&  	7.00524937 & 12.86578107	 & 	21.12624130 & 22.09037931\\
  &1.25 & 16.13984606	& 	6.01088776	 & 		8.80962709 &		12.52964177	 & 		16.27949336	\\ 
   \hline
   \end{tabular}
   \caption{Approximation of $\lambda_1, \lambda_2, \lambda_3$ and $\lambda_4$ with snapshot based on $u_1 + u_2 + u_3 +u_4$ at $h = 0.05$}
	 	\label{TABLE:12}
 	\end{table}
 	
 	
 \begin{table}
 	 	\centering
 \begin{tabular}{|c|c|c|c|c|c|c|c|} 
 		\hline
        	
 	        {\begin{tabular}[c]{@{}c@{}} $\mu$ \end{tabular}} &
 	        {\begin{tabular}[c]{@{}c@{}} Relative error\\$u_1, u_2, u_3, u_4$ \end{tabular}} &
 		 {\begin{tabular}[c]{@{}c@{}}  Relative error\\$u_1 + u_2 + u_3 + u_4$  \end{tabular}}\\
 \hline
   -1.25& $9.8 \times 10^{-7}$     & $6.1 \times 10^{-3}$ 				\\
  -0.75& $1.1 \times 10^{-7}$ & $2.7 \times 10^{-5}$   \\
  -0.50& $4.4 \times 10^{-8}$ & $5.6 \times 10^{-7}$   \\
  0.50& $4.1 \times 10^{-8}$ & $1.8 \times 10^{-5}$   \\
  0.75&$1.1 \times 10^{-7}$ &   $1.1 \times 10^{-4}$ \\
  1.25 & $9.7 \times 10^{-7}$  &	$9.0 \times 10^{-3}$	 \\ 
  \hline
 
 -1.25 &$5.3 \times 10^{-8}$ &$6.6 \times 10^{-3}$ 	\\
 -0.75 &$3.3 \times 10^{-8}$ &$6.9 \times 10^{-6}$   \\
 -0.50 &$7.5 \times 10^{-8}$ &$2.7 \times 10^{-7}$   \\
 0.50  &$7.5 \times 10^{-8}$ &$1.7 \times 10^{-5}$   \\
 0.75  &$3.2 \times 10^{-8}$ &$1.1 \times 10^{-4}$    \\
  \hline
   \end{tabular}
	\caption{Relative error between FEM and ROM based fourth eigenvalue with snapshot based on $u_1, u_2, u_3, u_4$ and $u_1 + u_2 + u_3+u_4$ at $h = 0.05$}
	 	\label{TABLE:13}
 	\end{table}
It is apparent from Table~\ref{TABLE:13} that the relative error in the case when we consider $u_1$, $u_2$, $u_3$, and $u_4$ in the snapshot matrix is slightly smaller than when considering the combination $u_1 + u_2 + u_3 + u_4$ in the snapshot matrix. However, the latter case is computationally cheaper.

\section{Numerical results for eigenvalue problem with multiple parameters}\label{se:2D}
\begin{figure}[H]
\centering
     \begin{subfigure}{0.4\textwidth}
        \includegraphics[height=5cm,width=7cm]{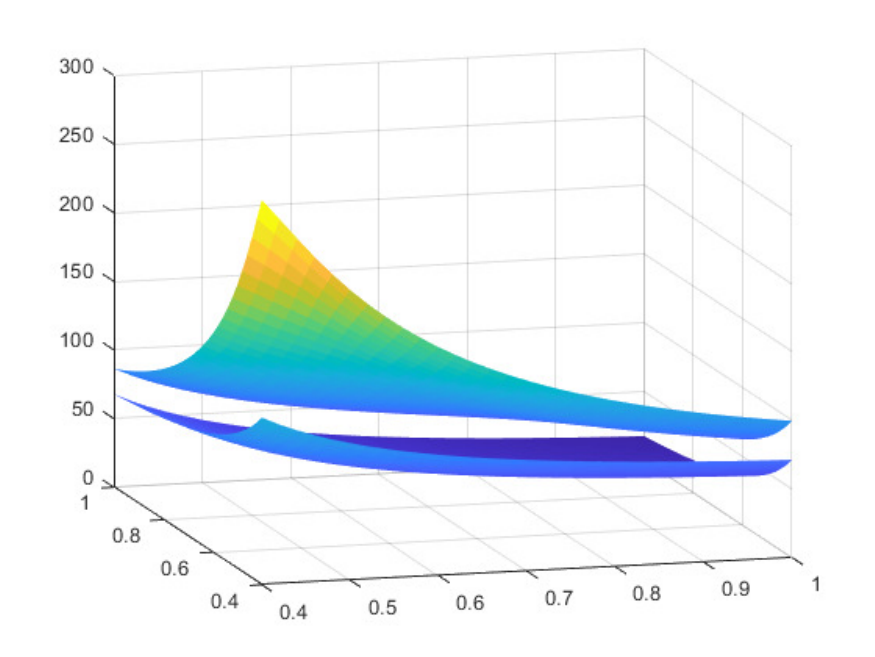}
       \caption{First and second}
     \end{subfigure}
     \begin{subfigure}{0.4\textwidth}
        \includegraphics[height=5cm,width=7cm]{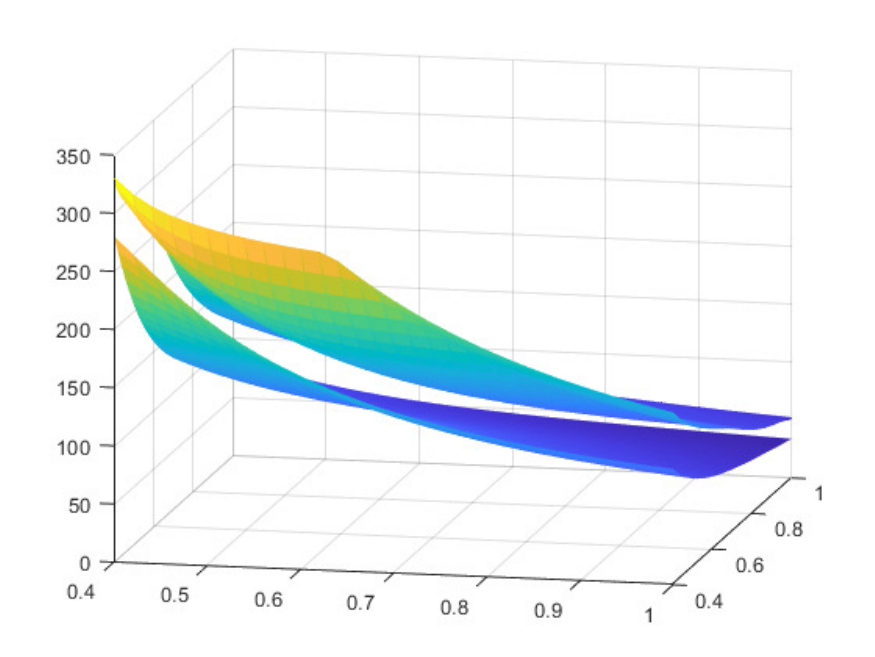}
         \caption{2nd and 3rd}
     \end{subfigure}\\
    \begin{subfigure}{0.4\textwidth}
 \includegraphics[height=5cm,width=7cm]{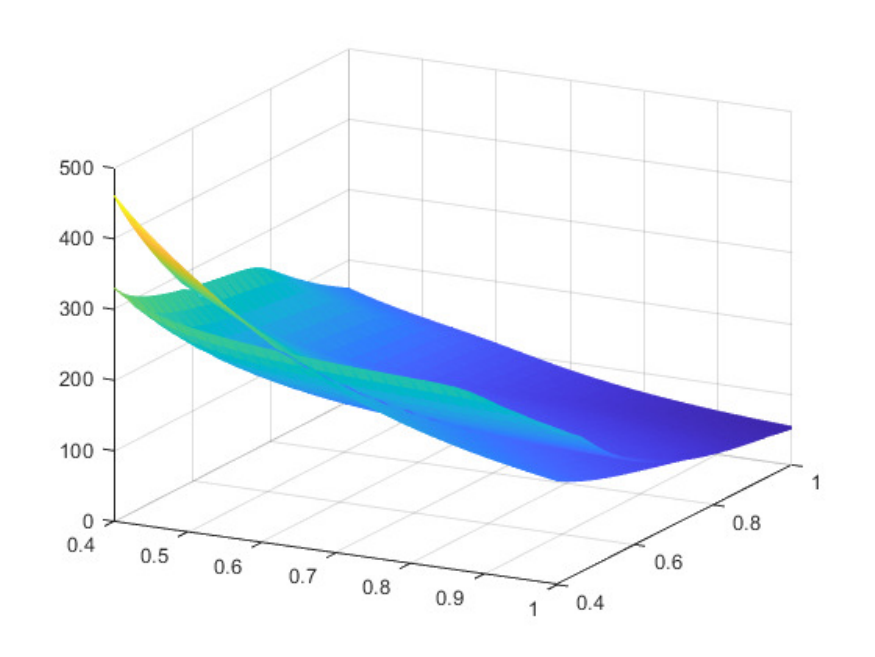}
         \caption{3rd and 4th}
     \end{subfigure}
      \begin{subfigure}{0.32\textwidth}
 \includegraphics[height=5cm,width=7cm]{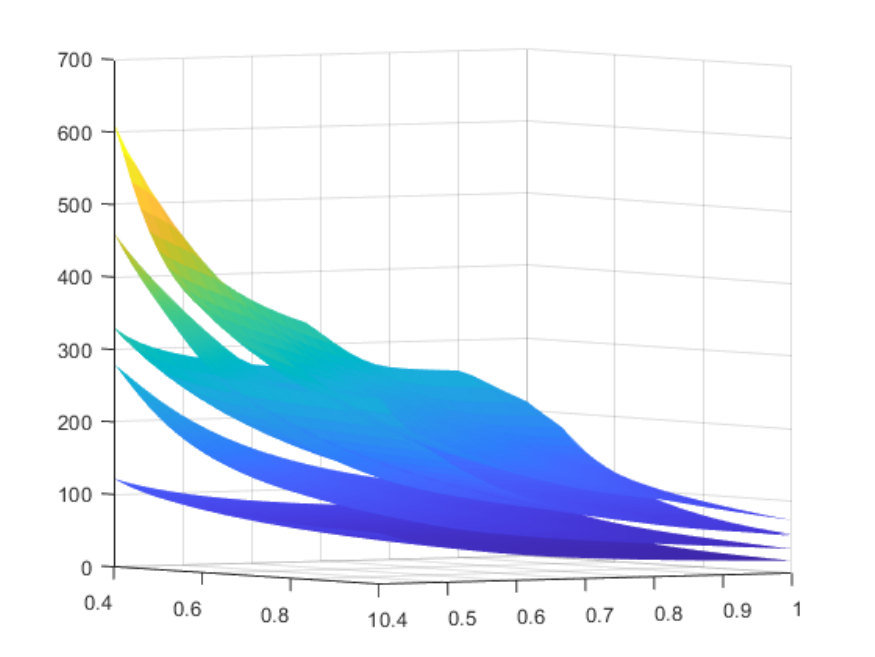}
         \caption{First five}
     \end{subfigure}
  \caption{Surface plot for the eigenvalues of the EVP \eqref{eq:PDE_example}.}
        \label{evplt2}
  \end{figure}

 In this section we investigate the behavior of the reduced eigenvalues and eigenvectors for different choices of snapshot matrix for the following eigenvalue problem with two parameters:
 \begin{equation}
\label{eq:PDE_example}
\left\{
\begin{array}{ll}
     -\operatorname{div}(A(\pmb{\mu})\nabla u(\pmb{\mu}))=\lambda(\pmb{\mu})u(\pmb{\mu})&  \textrm{ in }\Omega=(0,1)^2\\
     u(\pmb{\mu})=0& \textrm{ on }\partial\Omega 
\end{array}
\right.
\end{equation}
where the diffusion matrix $A(\pmb{\mu})\in\mathbb R^{2\times 2}$ is given by 
\begin{equation*}
    A(\pmb{\mu})= 
    \begin{pmatrix}
        \frac{1}{\mu_1^2} &\frac{0.7}{\mu_2}\\
        \frac{0.7}{\mu_2} &\frac{1}{\mu_2^2}
    \end{pmatrix},
\end{equation*}
with $\pmb{\mu}=(\mu_1,\mu_2)\in \mathcal{M}\subset \mathbb{R}^2.$
The problem is symmetric and the parameter space $\mathcal{M}$ is chosen in such a way 
that the matrix is positive definite. It can be easily checked that the matrix is positive definite, for instance, for any nonzero value of $\mu_2$ and $\mu_1\in (-1.42,1.42)\setminus\{0\}$. For our numerical tests, we choose the parameter space to be $\mathcal{M}=[0.4,1]^2$. In Figure~\ref{evplt2} we reported the surface plot for the eigenvalues of the eigenvalue problem (EVP) \eqref{evplt2}.
Numerical results with $25$ and $49$ uniform sample points are presented here. We reported the plot of the singular values of the snapshot matrix, plot for the FEM and the ROM eigenvalues for two sets of parameters, where $\mu_1$ varies and $\mu_2$ is fixed. Also, we have reported the values obtained using the FEM and the ROM at four  points $(0.5,0.6)$, $(0.5,0.8)$, $(0.8,0.6)$, and $(0.8,0.8)$. The number of POD basis is chosen using the criterion \eqref{criterion} with tolerance $10^{-8}$. To test the behavior of the eigenvalues, we have plotted the eigenvalues of the ROM by varying the number of POD basis at the point $\pmb{\mu}_{1,{\ast}}=(0.5,0.6)$. For all the experiments the mesh size is $h=0.05$.

\subsection{Reduced order method to obtain $\lambda_{1}$}
\begin{figure}
  \centering
     \begin{subfigure}{0.4\textwidth}
          \includegraphics[height=5cm,width=5.7cm]{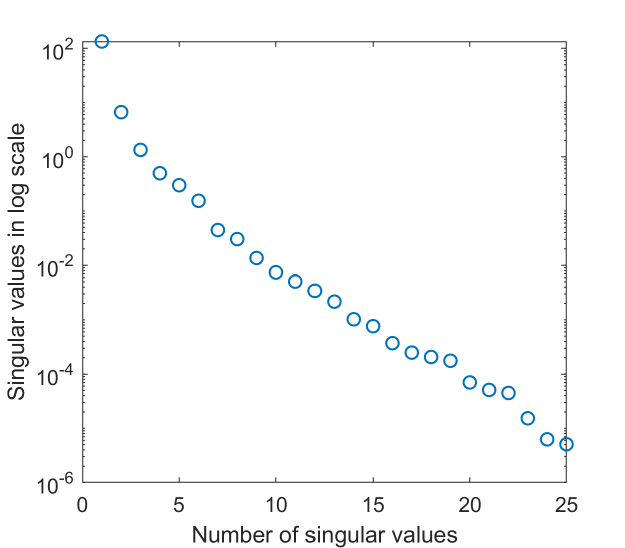}
         \caption{25 sampling}
     \end{subfigure}
     \begin{subfigure}{0.4\textwidth}
    \includegraphics[height=5cm,width=5.7cm]{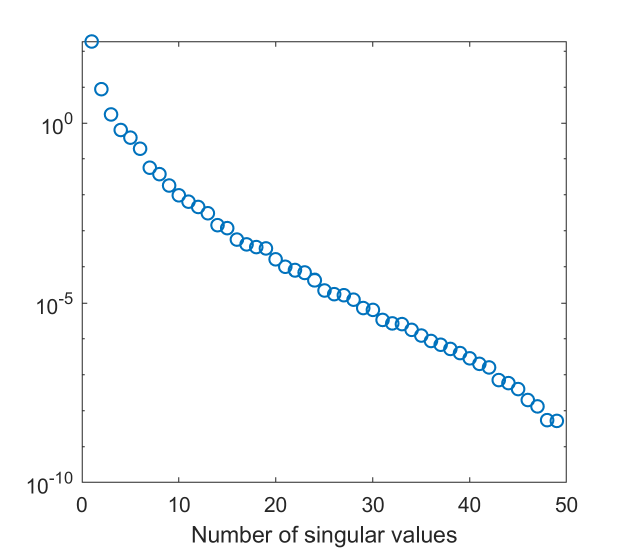}
         \caption{49 sampling}
     \end{subfigure}
  \caption{Singular values corresponding to snapshot matrix based on $u_1$.}
        \label{svd_u1}
  \end{figure}

\begin{figure}
     \centering
     \begin{subfigure}{0.4\textwidth}
         \includegraphics[height=5cm,width=5.7cm]{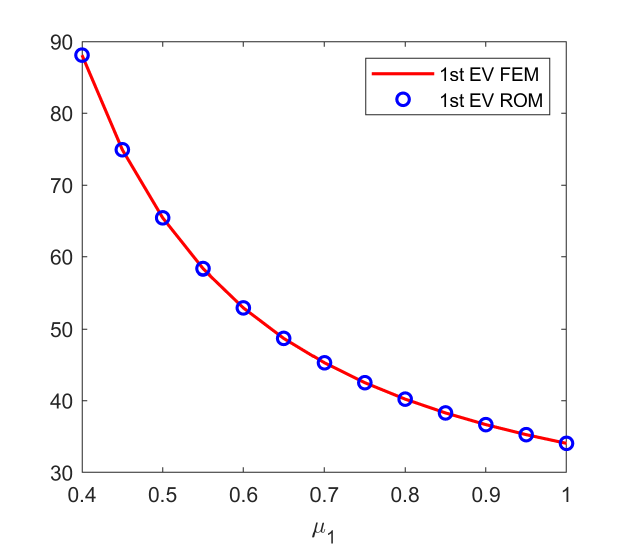}
         \caption{$\mu_2=0.6$}
     \end{subfigure}
      \begin{subfigure}{0.4\textwidth}
          \includegraphics[height=5cm,width=5.7cm]{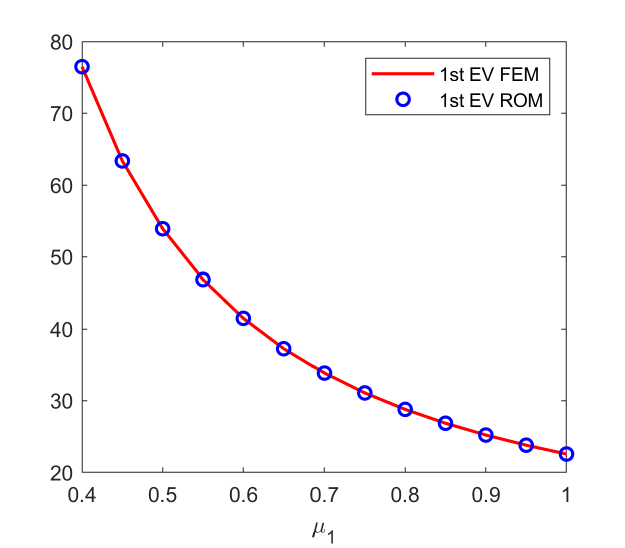}
          \caption{$\mu_2=0.8$}
     \end{subfigure}
  \caption{Approximation of $\lambda_1$ with snapshot based on $u_1$: comparison of FEM and ROM eigenvalues with varying $\mu_1$ and fixed $\mu_2$, and 49 sample points.}
        \label{evct_u1}
  \end{figure}

We consider the snapshot matrix consisting of the first eigenvectors column-wise at the sample parameters. In Figure~\ref{svd_u1} we have shown the plot for the singular values of the snapshot matrix and noticed that the singular values are decaying very fast. As expected, the first eigenvalues of the ROM match with the first eigenvalues of the FEM, which is evident in Figure~\ref{evct_u1}.

\begin{figure}
     \centering
     \begin{subfigure}{0.4\textwidth}
          \includegraphics[height=5cm,width=5.7cm]{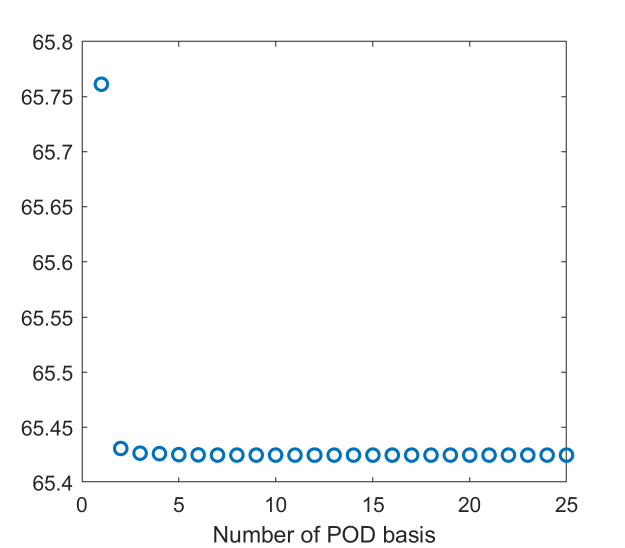}
         \caption{25 sampling points}
     \end{subfigure}
      \begin{subfigure}{0.4\textwidth}
          \includegraphics[height=5cm,width=5.7cm]{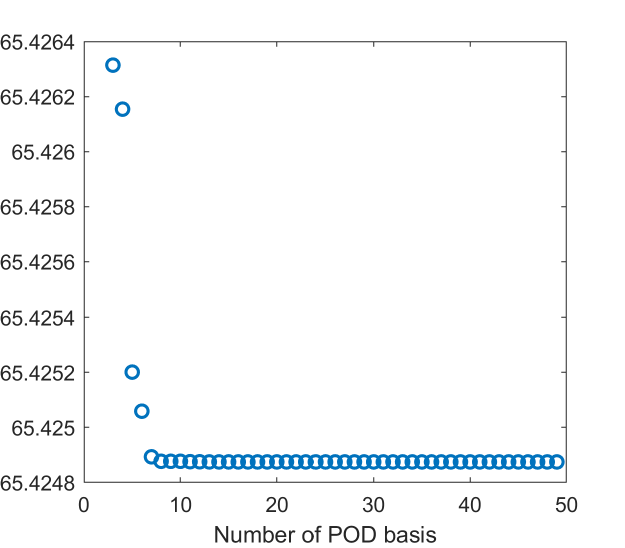}
         \caption{49 sampling points}
     \end{subfigure}
  \caption{Approximation of $\lambda_1$ at $\mu=(0.5,0.6)$ with snapshot based on $u_1$: varying the number of POD basis.}
        \label{evpod_u1}
  \end{figure}

The results of FEM and ROM eigenvalues and their relative errors at four sample points are reported in Table~\ref{table_u1} for all the sampling. The relative error are of order $10^{-7}$ -- $10^{-9}$. The number of POD dimensions is mentioned in Table~\ref{table_u1}, which is obtained using the criterion \eqref{criterion} with tolerance $10^{-8}$. In Figure~\ref{evpod_u1} we have plotted the first eigenvalue at the point $(0.5,0.6)$ using a different number of POD basis. The eigenvalues are converging to the exact eigenvalue with an increasing number of POD basis.

  \begin{table}
 	 	\centering
 	\begin{tabular}{|c|c|c|c|c|c|c|c|} 
 		\hline
 	 Sampling points&
 	 {\begin{tabular}[c]{@{}c@{}} ROM dim. \end{tabular}} &
		  {\begin{tabular}[c]{@{}c@{}} $\mu$ \end{tabular}} &
 		 {\begin{tabular}[c]{@{}c@{}} 1st EV(FEM) \end{tabular}} & 
 		{\begin{tabular}[c]{@{}c@{}}  1st EV(ROM)  \end{tabular}}&
 			{\begin{tabular}[c]{@{}c@{}} Rel. Error   \end{tabular}}\\
  \hline
25 & 8 &(0.5,0.6)&65.42487472& 65.42487526&$8.3\times 10^{-9}$ \\
&  &(0.5,0.8)&53.90654385&53.90654541 & $2.8\times 10^{-8}$\\
&  &(0.8,0.6)&40.22444914&40.22445184 &$6.7\times 10^{-8}$\\
&  &(0.8,0.8)&28.79395153& 28.79395783 & $2.1\times 10^{-7}$ \\ 
   \hline
  40 & 9 &(0.5,0.6)&65.42487472& 65.42487518&$7.0\times 10^{-9}$ \\
&  &(0.5,0.8)&53.90654385&53.90654555 & $3.1\times 10^{-8}$\\
&  &(0.8,0.6)&40.22444914&40.22444771 &$1.4\times 10^{-8}$\\
&  &(0.8,0.8)&28.79395186& 28.79395783 & $1.1\times 10^{-8}$ \\ 
  \hline
 \end{tabular}
\caption{Approximation of $\lambda_1$ with snapshot based on $u_1$: comparison of FEM and ROM.}
\label{table_u1}
 \end{table}
 
 \subsection{Reduced order method to obtain $\lambda_{2}$}
\begin{figure}
     \centering
     \begin{subfigure}{0.4\textwidth}
         \includegraphics[height=5cm,width=5.7cm]{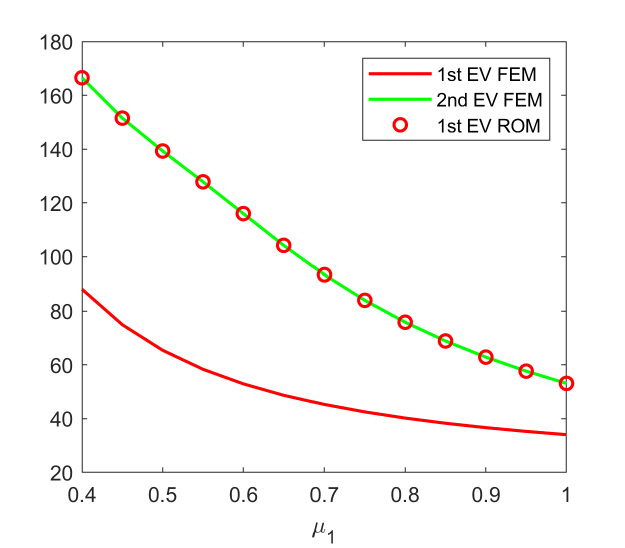}
         \caption{$\mu_2=0.6$}
     \end{subfigure}
      \begin{subfigure}{0.4\textwidth}
          \includegraphics[height=5cm,width=5.7cm]{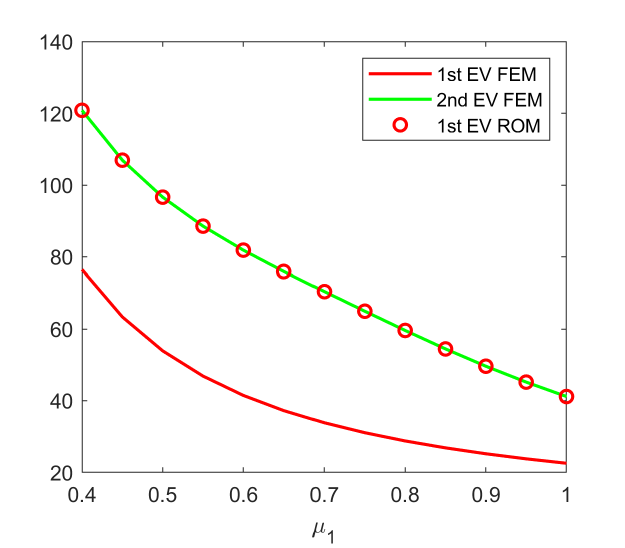}
          \caption{$\mu_2=0.8$}
     \end{subfigure}
  \caption{Approximation of $\lambda_2$ with snapshot based on $u_2$: comparison of FEM and ROM eigenvalues with varying $\mu_1$ and fixed $\mu_2$, and 49 sample points}
        \label{evct_u2}
  \end{figure}

We consider the snapshot matrix consisting of the second eigenvectors columnwise at the sample parameters. As expected, the first eigenvalue of the ROM matches with the second eigenvalue of the FEM, which is evident in Figure~\ref{evct_u2}.

\begin{figure}
     \centering
     \begin{subfigure}{0.4\textwidth}
          \includegraphics[height=5cm,width=5.7cm]{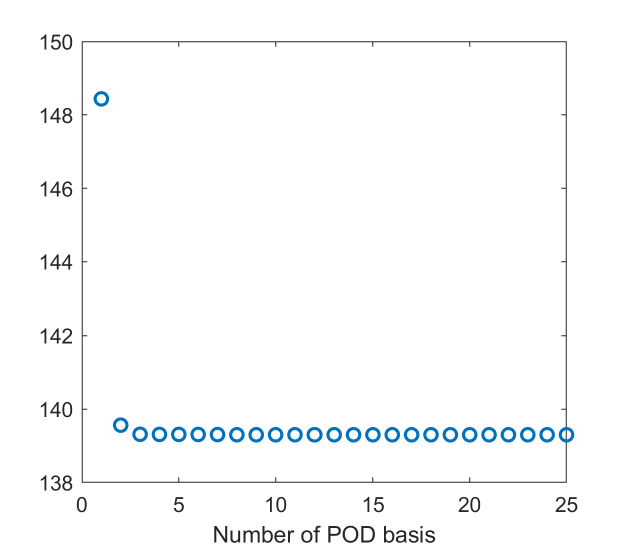}
         \caption{25 sampling points}
     \end{subfigure}
      \begin{subfigure}{0.4\textwidth}
          \includegraphics[height=5cm,width=5.7cm]{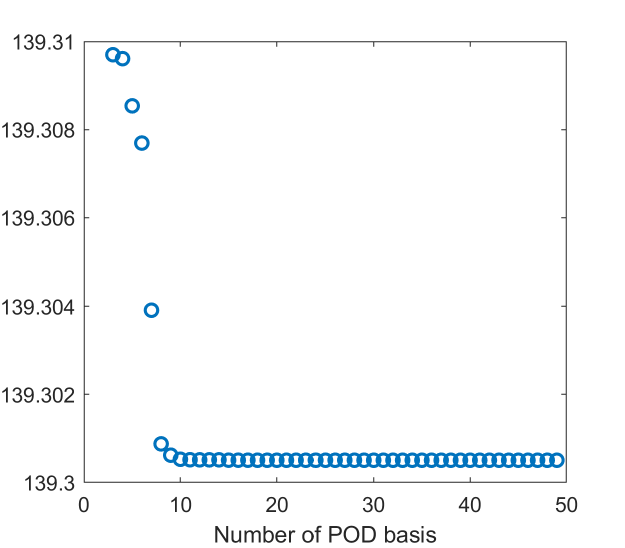}
         \caption{49 sampling points}
     \end{subfigure}
 \caption{Approximation of $\lambda_2$ at $\mu=(0.5,0.6)$ with snapshot based on $u_2$: varying the number of POD basis.}
        \label{evpod_u2}
  \end{figure}

The results of FEM and ROM eigenvalues and their relative errors at four sample points are reported in Table~\ref{table_u2} for all the sampling. The relative errors are of order $10^{-7}$ -- $10^{-8}$. The number of POD dimensions is mentioned in Table~\ref{table_u2}, which is obtained using the criterion \eqref{criterion} with tolerance $10^{-8}$. In Figure~\ref{evpod_u2} we have plotted the first eigenvalue at the point $(0.5,0.6)$ using a different number of POD basis. The eigenvalue is converging to the exact eigenvalue with an increasing number of POD basis.
  
\begin{table}
 	 	\centering
 	\begin{tabular}{|c|c|c|c|c|c|c|c|} 
 		\hline
 	Sampling& {\begin{tabular}[c]{@{}c@{}} ROM dim \end{tabular}} &
		  {\begin{tabular}[c]{@{}c@{}} $\mu$ \end{tabular}} &
 		 {\begin{tabular}[c]{@{}c@{}} 2nd EV(FEM) \end{tabular}} & 
 		{\begin{tabular}[c]{@{}c@{}}  1st EV(ROM)  \end{tabular}}&
 	 	{\begin{tabular}[c]{@{}c@{}}  Rel. Error  \end{tabular}}	\\
  \hline
25
 &12&(0.5,0.6)&139.30050529&139.30052177 & $1.2\times 10^{-7}$\\
&  &(0.5,0.8)& 96.61215367 &	96.61216806	& $1.5\times 10^{-7}$\\
&  &(0.8,0.6)&75.84831152&	75.84831902&$9.8\times 10^{-8}$	\\
&  &(0.8,0.8)& 59.58973796&59.58974344&$1.1\times 10^{-7}$	\\ 
  \hline
49
 &12&(0.5,0.6)&139.30050529&139.30051498 &$6.9\times 10^{-8}$\\
&  &(0.5,0.8)& 96.61215367 &96.61216466	&$1.1\times 10^{-7}$ \\
&  &(0.8,0.6)&75.84831152&75.84831665&$6.7\times 10^{-8}$	\\
&  &(0.8,0.8)& 59.58973796&59.58974473&$9.1\times 10^{-8}$	\\ 
  \hline
\end{tabular}
\caption{Approximation of $\lambda_2$ with snapshot based on $u_2$: comparison of FEM and ROM.}
\label{table_u2}
 \end{table}

\subsection{Reduced order method to obtain $\lambda_{3}$}
In this subsection, we discuss the results for third eigenvalues considering different combinations of eigenvectors in the snapshot matrix.
\subsubsection{Results of the EVP considering only $u_3$ in the snapshot matrix}
\begin{figure}
     \centering
     \begin{subfigure}{0.4\textwidth}
         \includegraphics[height=5cm,width=5.7cm]{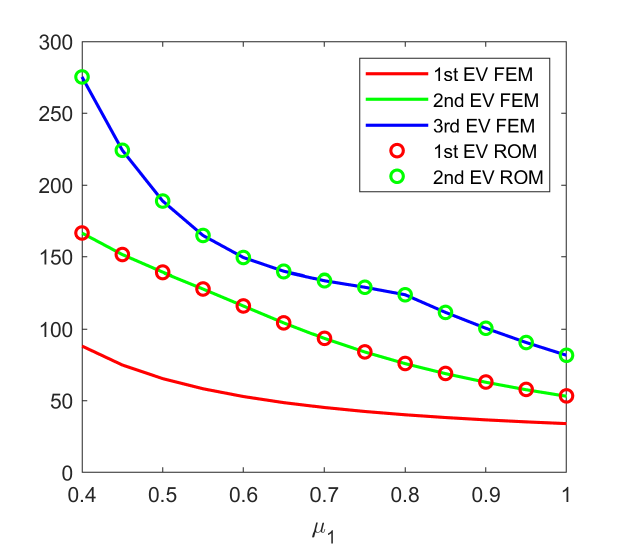}
         \caption{$\mu_2=0.6$}
     \end{subfigure}
      \begin{subfigure}{0.4\textwidth}
          \includegraphics[height=5cm,width=5.7cm]{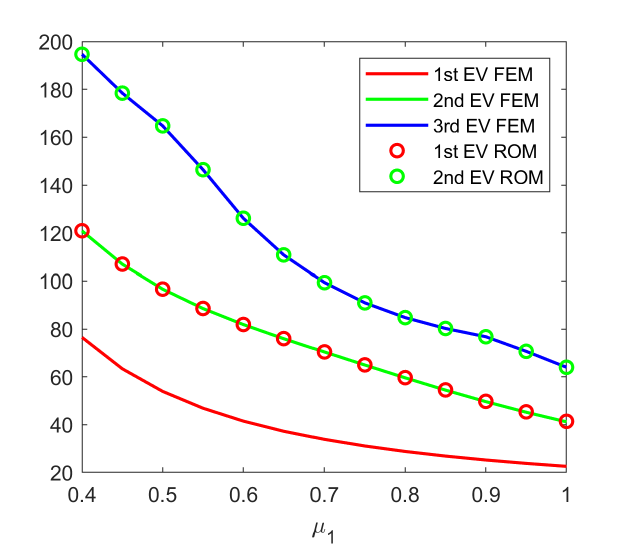}
          \caption{$\mu_2=0.8$}
     \end{subfigure}
  \caption{Approximation of $\lambda_3$ with snapshot based on $u_3$: comparison of FEM and ROM eigenvalues with varying $\mu_1$ and fixed $\mu_2$, and 49 sample points}
        \label{evct_u3}
  \end{figure}

\begin{figure}
     \centering
     \begin{subfigure}{0.4\textwidth}
          \includegraphics[height=5cm,width=5.7cm]{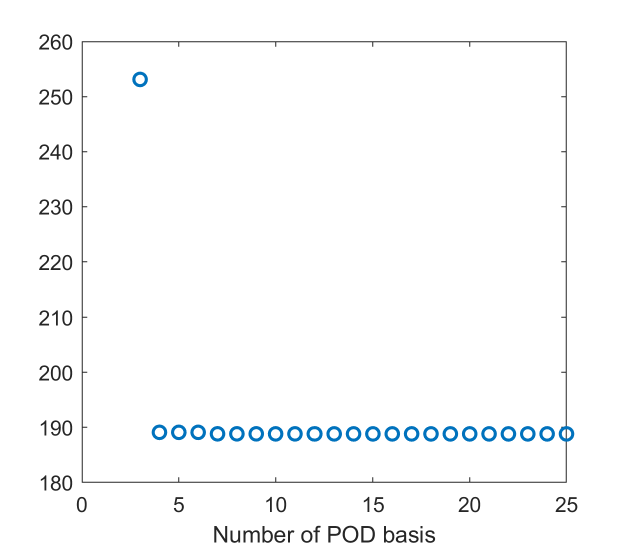}
         \caption{25 sampling points}
     \end{subfigure}
      \begin{subfigure}{0.4\textwidth}
          \includegraphics[height=5cm,width=5.7cm]{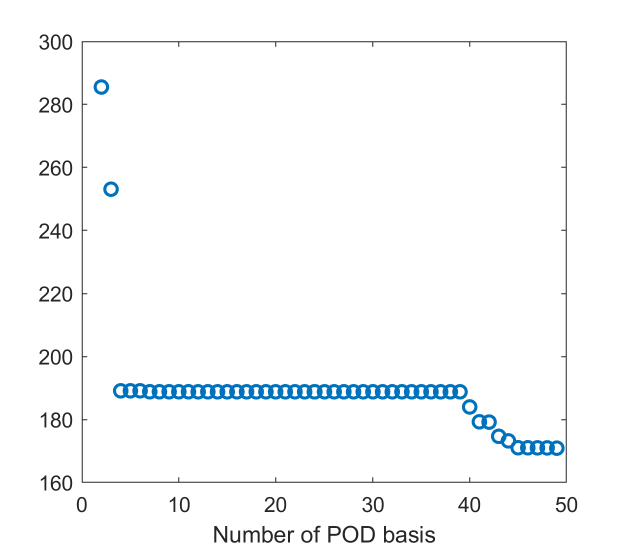}
         \caption{49 sampling points}
     \end{subfigure}
\caption{Approximation of $\lambda_3$ at $\mu=(0.5,0.6)$ with snapshot based on $u_3$: varying the number of POD basis.}
        \label{evpod_u3}
  \end{figure}

Let us consider the snapshot matrix containing only the third eigenvectors at the sample points.
In this case, from Figure~\ref{evct_u3} it can be seen that the first eigenvalue of the ROM matches the second eigenvalue of the FEM, and the second eigenvalue of the ROM matches the third eigenvalue of the FEM. But, as we consider $u_3$ in the snapshot matrix, one should expect that the first eigenvalue of the ROM should match with the third eigenvalue of the FEM. The reason behind this is the fact that the $L^2$ inner product $(u_2(\mu_i),u_3(\mu_j)$ is not zero for $i\neq j$, so the snapshots contain some component of the second eigenvector. If we calculate the inner product $(u_1(\mu_i),u_2(\mu_j)$ then we see that these values are small for $i\neq j$, that is why we are not getting the first eigenvalue. Note that the inner product is zero for $i=j$.  From Table~\ref{table_u3} it is observed that the maximum relative error among the four test points is $10^{-6}$. In Figure~\ref{evpod_u3} we have plotted the second eigenvalue of the ROM at the point $(0.5,0.6)$ with an increasing number of ROM dimensions. From Figure~\ref{evpod_u3} one can see that if the number of POD basis is up to 38 then the second ROM eigenvalue is converging to the exact one, while then it starts to decrease.

\begin{figure}
     \centering
     \begin{subfigure}{0.32\textwidth}
         \includegraphics[height=5cm,width=5.5cm]{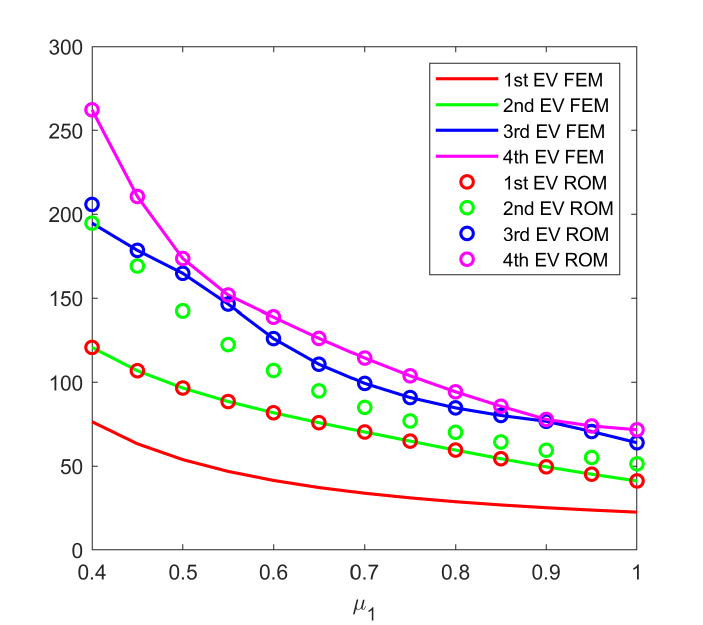}
         \caption{49 samples}
     \end{subfigure}
      \begin{subfigure}{0.32\textwidth}
          \includegraphics[height=5cm,width=5.5cm]{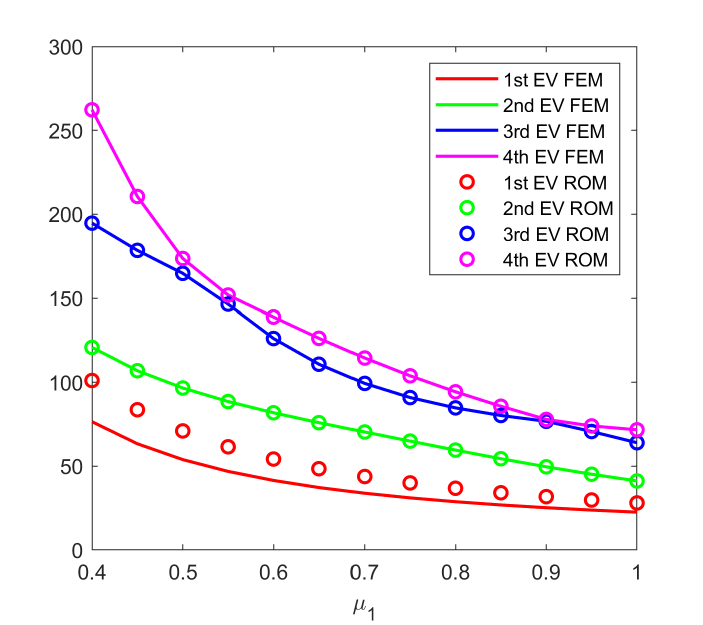}
          \caption{81 samples}
     \end{subfigure}
     \begin{subfigure}{0.32\textwidth}
          \includegraphics[height=5cm,width=5.5cm]{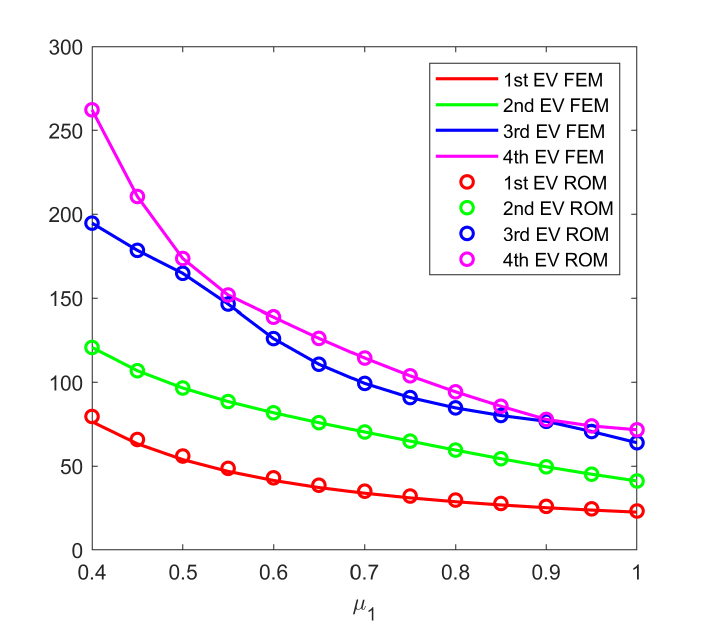}
          \caption{121 samples}
     \end{subfigure}
  \caption{Approximation of $\lambda_3$ with snapshot based on $u_3$: comparison of FEM and ROM eigenvalues with varying $\mu_1$ and  $\mu_2=0.8$, and with different number of sample points}
        \label{evct_u3_exp}
  \end{figure}

In order to investigate this issue in more detail, we use more sample points and consider all left singular vectors as ROM basis and plot the first four eigenvalues of the FEM and the ROM for $\mu_2=0.8$ in Figure~\ref{evct_u3_exp}; we can observe that with the increasing number of sample points all the first four eigenvalues of the ROM converge to the corresponding eigenvalues of the FEM.
  
 \begin{table}
 	 	\centering
 	\begin{tabular}{|c|c|c|c|c|c|c|c|} 
 		\hline
 	Sampling&{\begin{tabular}[c]{@{}c@{}} ROM dim \end{tabular}} &
		  {\begin{tabular}[c]{@{}c@{}} $\mu$ \end{tabular}} &
 		 {\begin{tabular}[c]{@{}c@{}} 3rd EV(FEM) \end{tabular}} &
 		{\begin{tabular}[c]{@{}c@{}}  2nd EV(ROM) \end{tabular}} &
 		{\begin{tabular}[c]{@{}c@{}}  Rel. Error \end{tabular}}\\
  \hline
25
&21 &(0.5,0.6)&188.82589386&188.82716018 &$6.7\times 10^{-6}$\\
&  &(0.5,0.8)&164.71911477&164.71911669 &$1.1\times 10^{-8}$\\
&  &(0.8,0.6)&123.85385935&123.85419564 &$2.7\times 10^{-6}$\\
&  &(0.8,0.8)&84.76248191&84.76248514 &$3.8 \times 10^{-8}$ \\
  \hline
 49
&24 &(0.5,0.6)&188.82589386 &188.82590770 &$7.3\times 10^{-8}$\\
&  &(0.5,0.8)&164.71911477&164.71911878 &$2.4\times 10^{-8}$\\
&  &(0.8,0.6)&123.85385935&123.85387416 &$1.1\times 10^{-7}$\\
&  &(0.8,0.8)&84.76248191&84.76248574 &$4.5 \times 10^{-8}$ \\
  \hline
 \end{tabular}
\caption{Approximation of $\lambda_3$ with snapshot based on $u_3$: comparison of FEM and ROM.}
\label{table_u3}
 \end{table} 
 
\subsubsection{Results of the EVP considering $u_1$, $u_2$, and $u_3$ in the snapshot matrix}
\begin{figure}
     \centering
     \begin{subfigure}{0.4\textwidth}
         \includegraphics[height=5cm,width=5.7cm]{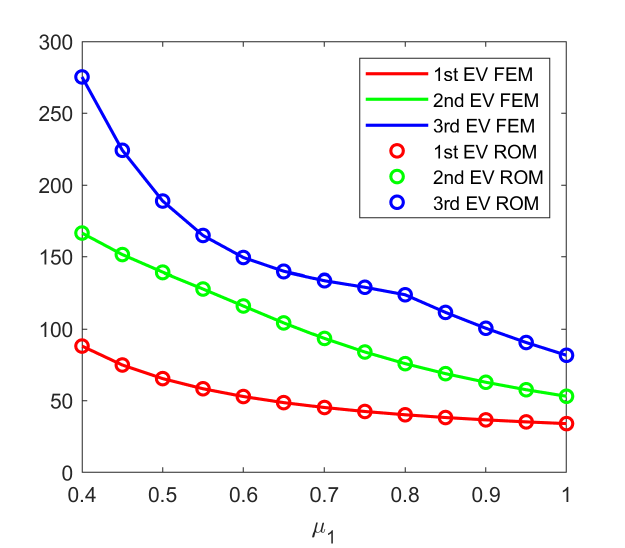}
         \caption{$\mu_2=0.6$}
     \end{subfigure}
      \begin{subfigure}{0.4\textwidth}
          \includegraphics[height=5cm,width=5.7cm]{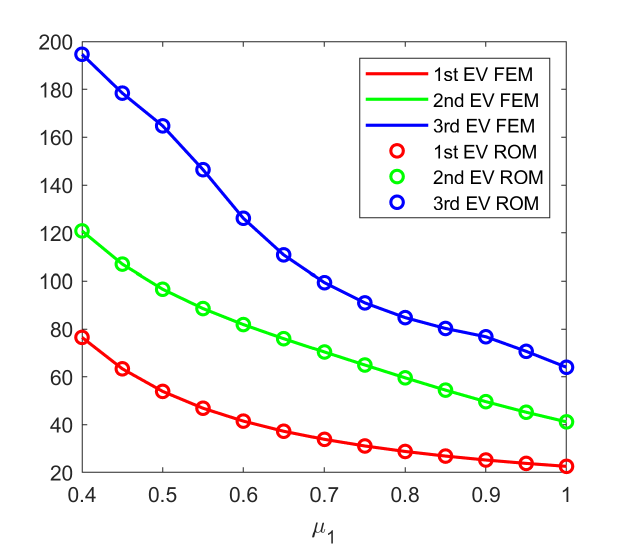}
          \caption{$\mu_2=0.8$}
     \end{subfigure}
  \caption{Approximation of $\lambda_3$ with snapshot based on $u_1.u_2,u_3$: comparison of FEM and ROM eigenvalues with varying $\mu_1$ and fixed $\mu_2$, and 49 sample points.}
        \label{evct_u13}
  \end{figure}

We consider the snapshot matrix containing the first three eigenvectors at the sample points. In this case, all three ROM eigenvalues coincide with the corresponding FEM eigenvalues, as it can be seen in Figure~\ref{evct_u13}.

\begin{figure}
     \centering
     \begin{subfigure}{0.4\textwidth}
          \includegraphics[height=5cm,width=5.7cm]{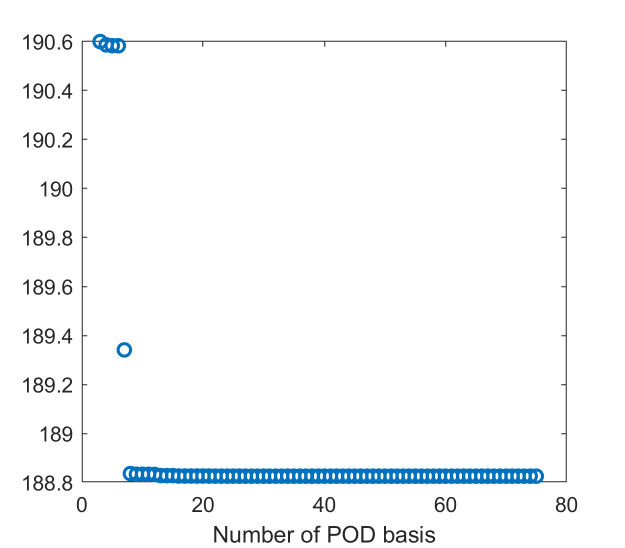}
         \caption{25 sampling points}
     \end{subfigure}
      \begin{subfigure}{0.4\textwidth}
          \includegraphics[height=5cm,width=5.7cm]{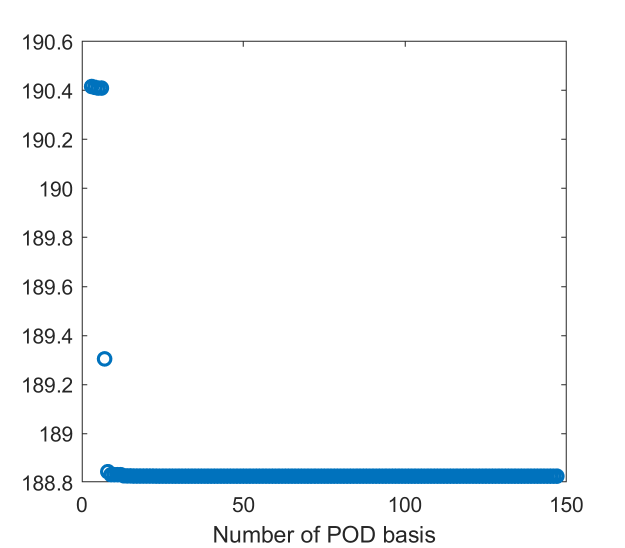}
         \caption{49 sampling points}
     \end{subfigure}
 \caption{Approximation of $\lambda_3$ at $\mu=(0.5,0.6)$ with snapshot based on $u_1,u_2,u_3$: varying the number of POD basis.}
        \label{evpod_u13}
  \end{figure}

The eigenvalues of the FEM and the ROM at the four sample points are reported in Table~\ref{table_u13}. The relative error between the FEM and ROM eigenvalues are also reported and the maximum error is $10^{-6}$. Even if we are increasing the number of POD basis then also the 3rd eigenvalue corresponding to the ROM converges to the 3rd eigenvalue of the FEM as it is shown in Figure~\ref{evpod_u13}.

  \begin{table}
 	 	\centering
 	\begin{tabular}{|c|c|c|c|c|c|c|c|} 
 		\hline
 	Sampling&{\begin{tabular}[c]{@{}c@{}} dim ROM \end{tabular}} & 
		  {\begin{tabular}[c]{@{}c@{}} $\mu$ \end{tabular}} &
 		{\begin{tabular}[c]{@{}c@{}} 3rd EV(FEM)  \end{tabular}}&
 		{\begin{tabular}[c]{@{}c@{}} 3rd EV(ROM)  \end{tabular}}&
         {\begin{tabular}[c]{@{}c@{}}  Rel. Error ($\lambda_3$) \end{tabular}}\\
  \hline
25
&32 &(0.5,0.6)&188.82589386 &188.82590469&$5.7\times 10^{-8}$  \\
 & &(0.5,0.8)&164.71911477 &164.71917647 &$3.7\times 10^{-7}$ \\
 & &(0.8,0.6)&123.85385935 &123.85438852 &$4.2\times 10^{-6}$ \\
&  &(0.8,0.8)&84.76248191 &84.76249669 &$1.7\times 10^{-7}$  \\ 
  \hline
 45
&33 &(0.5,0.6)&188.82589386  &188.82591285 &$1.0\times 10^{-7}$   \\
 & &(0.5,0.8)&164.71911477  &164.71913934&$1.4\times 10^{-7}$  \\
 & &(0.8,0.6)&123.85385935   &123.85389384&$2.7\times 10^{-7}$   \\
&  &(0.8,0.8)&84.76248191&84.76249350 &$1.3\times 10^{-7}$   \\
  \hline
 \end{tabular}
\caption{Approximation of $\lambda_3$ with snapshot based on $u_1,u_2,u_3$: comparison of FEM and ROM.}
\label{table_u13}
 \end{table}

\subsubsection{Results of the EVP considering $u_1+u_2+u_3$ in the snapshot matrix}
\begin{figure}
     \centering
     \begin{subfigure}{0.4\textwidth}
         \includegraphics[height=5cm,width=5.7cm]{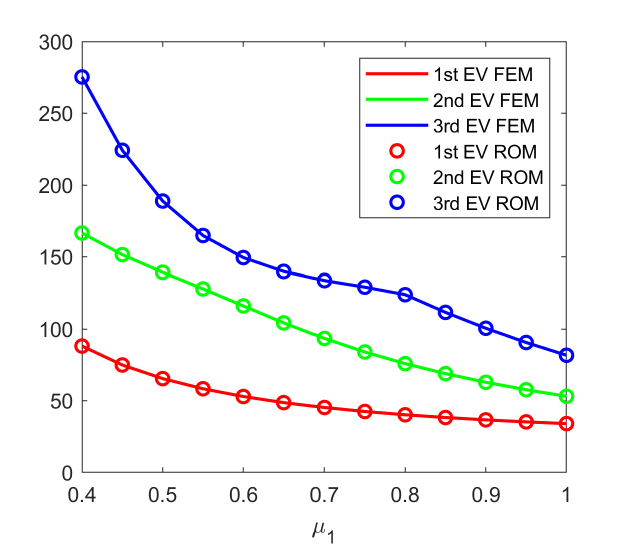}
         \caption{$\mu_2=0.6$}
     \end{subfigure}
      \begin{subfigure}{0.4\textwidth}
          \includegraphics[height=5cm,width=5.7cm]{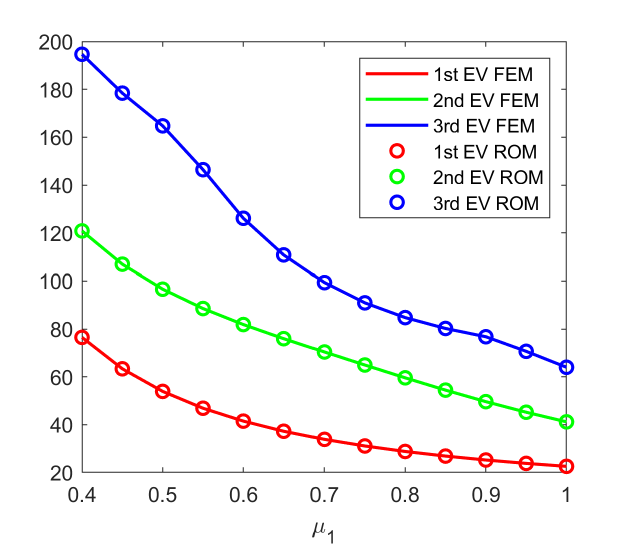}
          \caption{$\mu_2=0.8$}
     \end{subfigure}
  \caption{Approximation of $\lambda_1,\lambda_2,\lambda_3$ with snapshot based on $u_1+u_2+u_3$: comparison of FEM and ROM eigenvalues with varying $\mu_1$ and fixed $\mu_2$, and 49 sample points.}
        \label{evct_au3}
  \end{figure}

The results of the ROM are good when we consider the first three eigenvectors in the snapshot matrix and the order of the eigenvalues is also preserved. However, the number of snapshots is three times the number of sample points. In order to try to reduce the computational cost, we consider the snapshot matrix containing the sum of the first three eigenvectors at the sample points as the column. Then the number of columns of the snapshot matrix is equal to the number of snapshots and the snapshots contain some components of all three eigenvectors. In this case, also, the first three eigenvalues of the ROM coincide with the first eigenvalues of the FEM and preserve the order, see Figure~\ref{evct_au3}.

The 3rd eigenvalues of the ROM and the FEM at the four test points and the corresponding relative errors are reported in Table~\ref{table_au3}. The maximum relative error is $10^{-4}$. Note that by taking the combination of the eigenvectors we reduce the number of snapshots at the price of having results which are not as good as in the former case when we considered all three eigenvectors. In Figure~\ref{fig:3rdev_moataj} we reported the third eigenvectors obtained using the FEM and the ROM at the four sample points and they are matching up to the sign.

 	\begin{table}
 	 	\centering
 	\begin{tabular}{|c|c|c|c|c|c|c|c|} 
 		\hline
 	Sampling&{\begin{tabular}[c]{@{}c@{}} ROM dim \end{tabular}} &
		  {\begin{tabular}[c]{@{}c@{}} $\mu$ \end{tabular}} &
 	 {\begin{tabular}[c]{@{}c@{}} 3rd EV(FEM) \end{tabular}} & 
 		{\begin{tabular}[c]{@{}c@{}} 3rd EV(ROM)  \end{tabular}}&
    {\begin{tabular}[c]{@{}c@{}}  Rel. Error($\lambda_3$) \end{tabular}}\\
  \hline
 25
& 22 &(0.5,0.6)& 188.82589386 &188.82837895 &$1.3\times 10^{-5}$  \\
&  &(0.5,0.8)&164.71911477 & 164.72021897&$6.7\times 10^{-6}$ \\
&  &(0.8,0.6)& 123.85385935  & 123.87322460&$1.5\times 10^{-4}$  	\\
&  &(0.8,0.8)&84.76248191 & 84.76361731&$1.3\times 10^{-5}$ \\ 
  \hline
 49
& 28 &(0.5,0.6)&188.82589386&188.82593692&$2.2 \times 10^{-7}$   \\
&  &(0.5,0.8)&164.71911477& 164.71919081&$4.6\times 10^{-7}$ \\
&  &(0.8,0.6)&123.85385935  & 123.85406453 &$1.6\times 10^{-6}$ \\
&  &(0.8,0.8)&84.76248191 & 84.76249965&$2.0\times 10^{-7}$  \\ 
  \hline
 \end{tabular}
\caption{Approximation of $\lambda_3$ with snapshot based on $u_1+u_2+u_3$: comparison of FEM and ROM.}
\label{table_au3}
 \end{table}

 \begin{figure}
     \centering
     \begin{subfigure}{0.4\textwidth}
          \includegraphics[height=5cm,width=5.7cm]{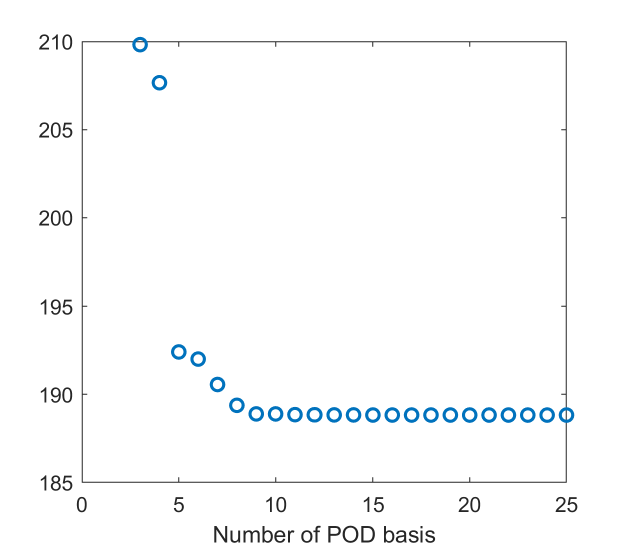}
         \caption{25 sampling points}
     \end{subfigure}
      \begin{subfigure}{0.4\textwidth}
          \includegraphics[height=5cm,width=5.7cm]{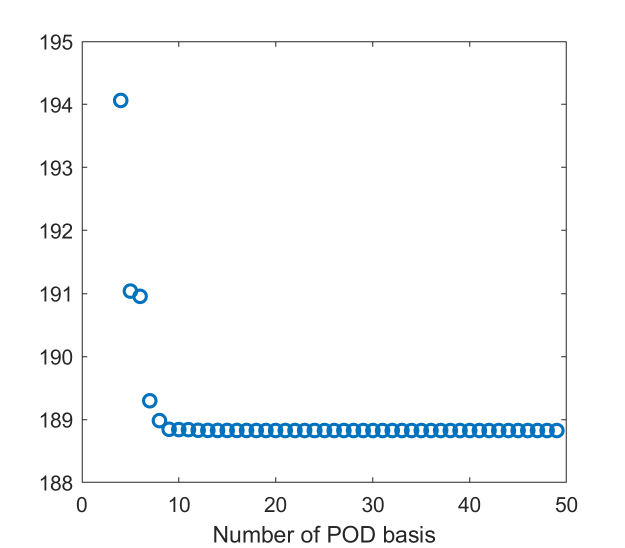}
         \caption{49 sampling points}
     \end{subfigure}
 \caption{Approximation of $\lambda_3$ at $\mu=(0.5,0.6)$ with snapshot based on $u_1+u_2+u_3$: varying the number of POD basis.}
        \label{evpod_au3}
  \end{figure}
 \begin{figure}
     \centering
     \begin{subfigure}{0.22\textwidth}
         \centering
         \includegraphics[height=4cm,width=4cm]{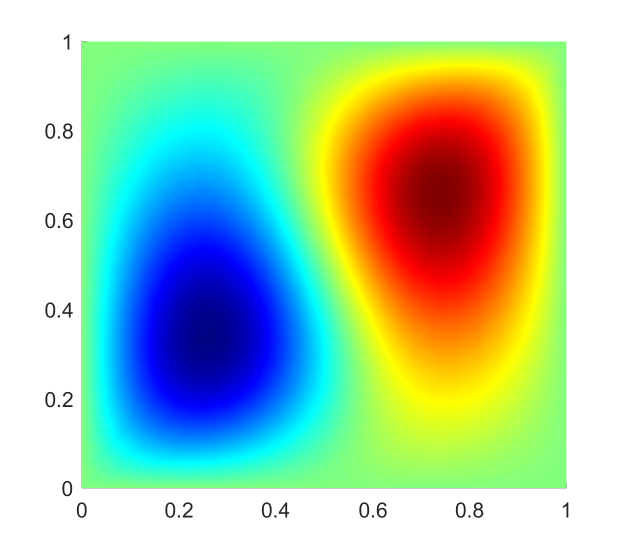}
     \end{subfigure}
     \begin{subfigure}{0.22\textwidth}
         \centering
         \includegraphics[height=4cm,width=4cm]{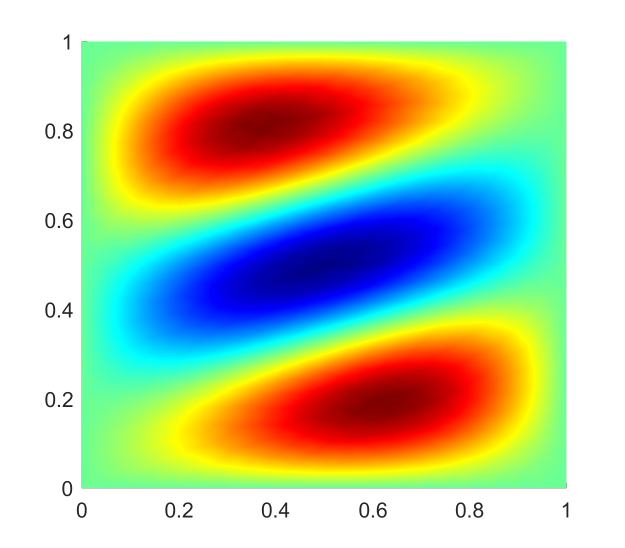}
     \end{subfigure}
     \begin{subfigure}{0.22\textwidth}
         \centering
         \includegraphics[height=4cm,width=4cm]{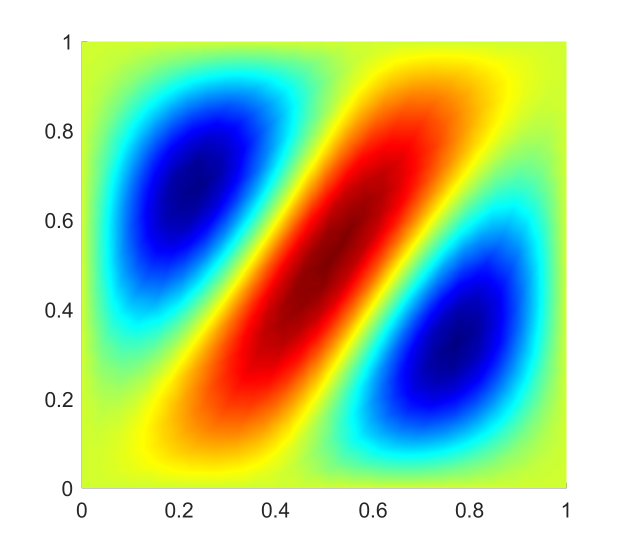}
     \end{subfigure}
     \begin{subfigure}{0.22\textwidth}
         \centering
         \includegraphics[height=4cm,width=4cm]{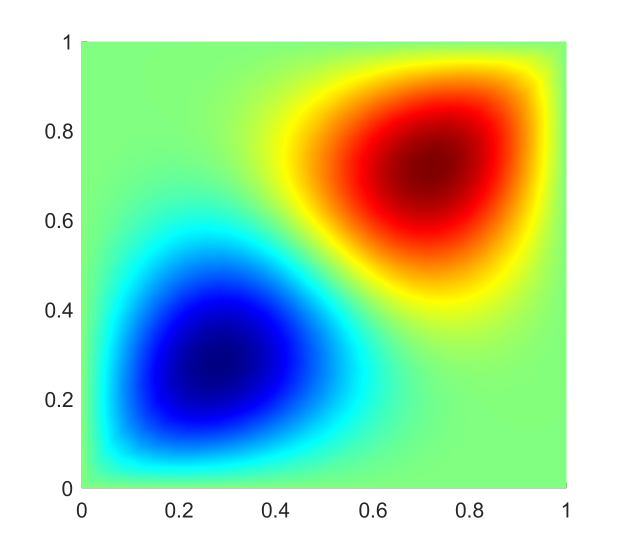}
     \end{subfigure}\\
     \begin{subfigure}{0.22\textwidth}
         \centering
         \includegraphics[height=4cm,width=4cm]{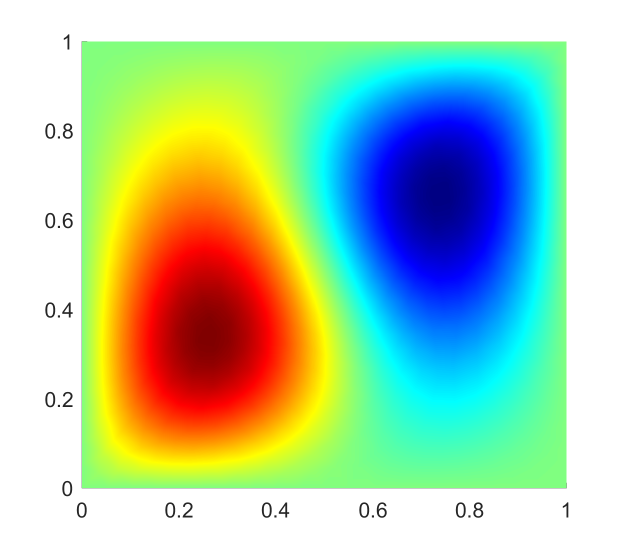}
     \end{subfigure}
     \begin{subfigure}{0.22\textwidth}
         \centering
         \includegraphics[height=4cm,width=4cm]{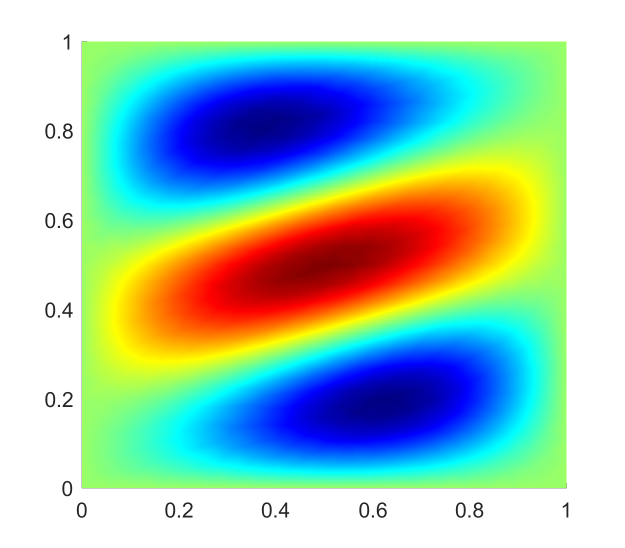}
     \end{subfigure}
     \begin{subfigure}{0.22\textwidth}
         \centering
         \includegraphics[height=4cm,width=4cm]{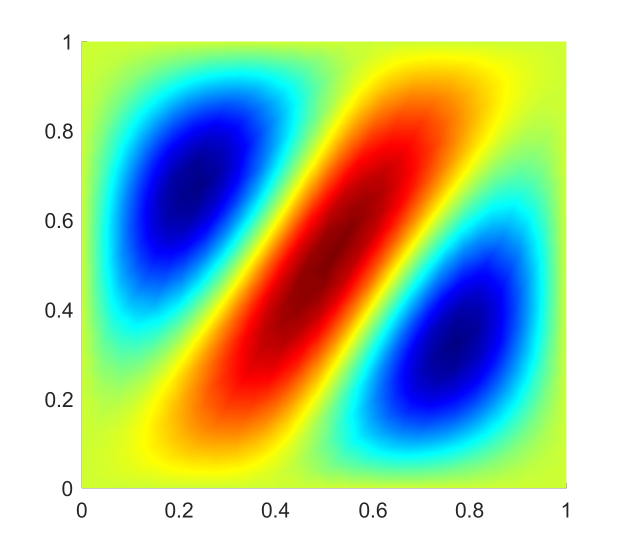}
     \end{subfigure}
     \begin{subfigure}{0.22\textwidth}
         \centering
         \includegraphics[height=4cm,width=4cm]{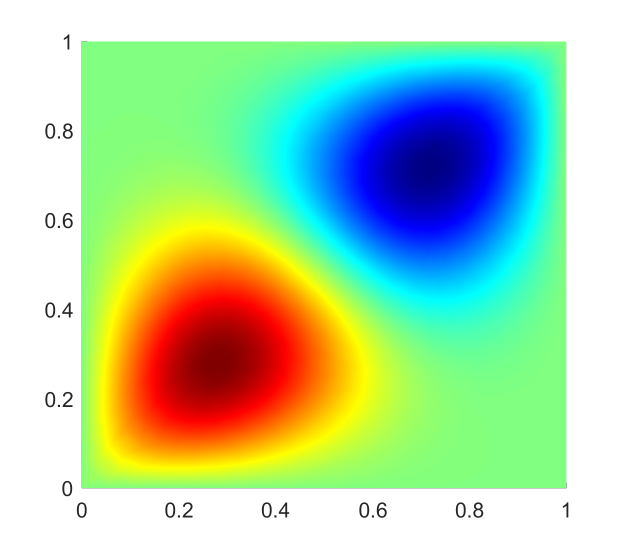}
     \end{subfigure}
        \caption{Comparison of 3rd eigenvectors using FEM (1st row) and ROM (2nd row) with snapshot based on $u_1+u_2+u_3$ at four points.}
        \label{fig:3rdev_moataj}
  \end{figure}


\subsection{Reduced order method to obtain $\lambda_4$}
In this subsection we will find the fourth eigenvalue of the problem using different choices of snapshot matrix.
\subsubsection{Results of the EVP considering only $u_4$ in the snapshot matrix}
\begin{figure}
     \centering
     \begin{subfigure}{0.4\textwidth}
         \includegraphics[height=5cm,width=5.7cm]{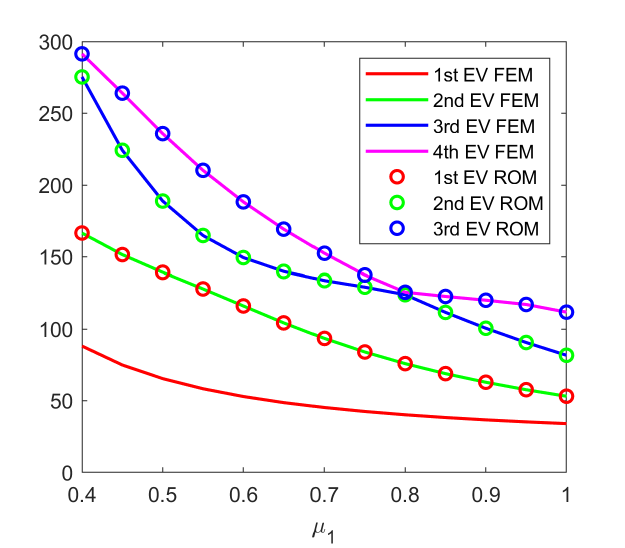}
         \caption{$\mu_2=0.6$}
     \end{subfigure}
      \begin{subfigure}{0.4\textwidth}
          \includegraphics[height=5cm,width=5.7cm]{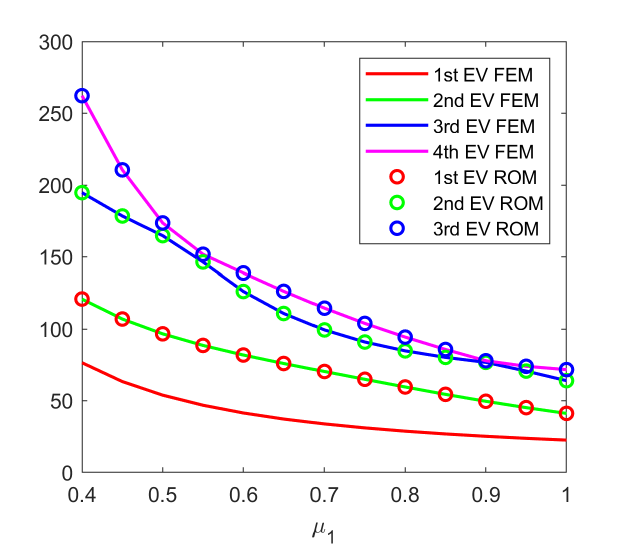}
          \caption{$\mu_2=0.8$}
     \end{subfigure}
  \caption{Approximation of $\lambda_4$ with snapshot based on $u_4$: comparison of FEM and ROM eigenvalues with varying $\mu_1$ and fixed $\mu_2$, and 49 sample points.}
        \label{evct_u4}
  \end{figure}

\begin{figure}
     \centering
     \begin{subfigure}{0.32\textwidth}
         \includegraphics[height=5cm,width=5.5cm]{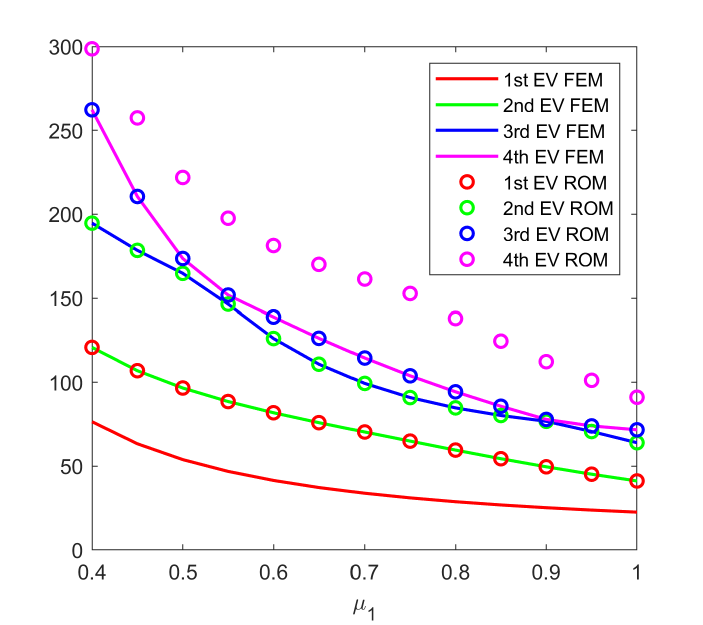}
         \caption{49 samples}
     \end{subfigure}
      \begin{subfigure}{0.32\textwidth}
          \includegraphics[height=5cm,width=5.5cm]{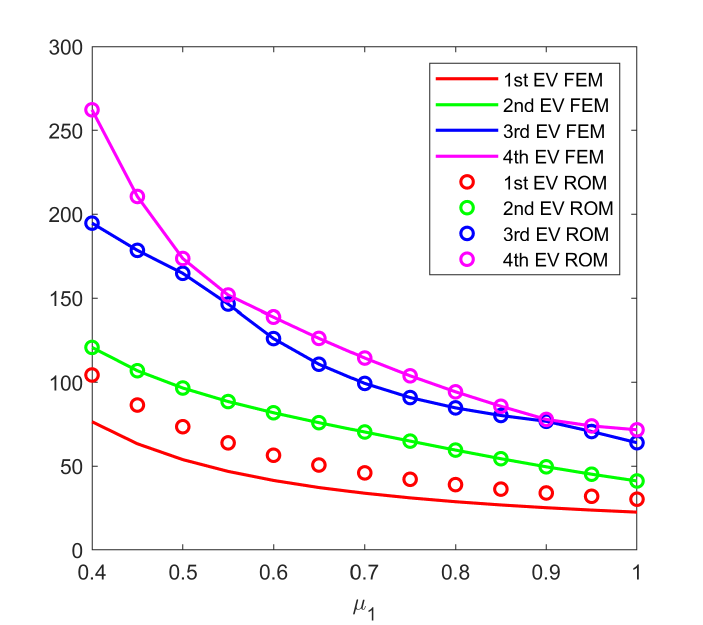}
          \caption{81 samples}
     \end{subfigure}
     \begin{subfigure}{0.32\textwidth}
          \includegraphics[height=5cm,width=5.5cm]{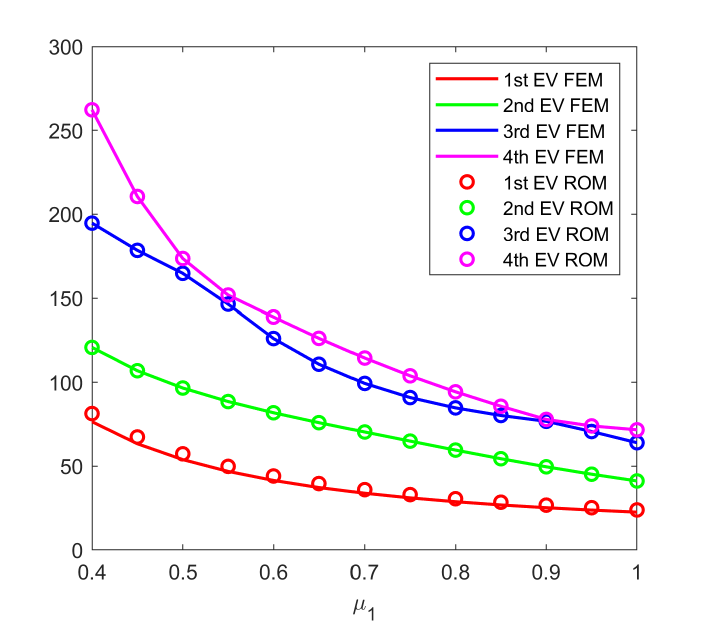}
          \caption{121 samples}
     \end{subfigure}
  \caption{Approximation of $\lambda_1,\lambda_2,\lambda_3,\lambda_4$ with snapshot based on $u_4$: comparison of FEM and ROM eigenvalues with varying $\mu_1$ and  $\mu_2=0.8$, and with different number of sample points.}
        \label{evct_u4_exp}
  \end{figure}

First, we consider the snapshot matrix containing only the fourth eigenvector at the sample points. In Figure~\ref{evct_u4} we have shown the first four eigenvalues of the FEM and the first three eigenvalues of the ROM for the parameters $\mu_1$ ranging from $0.4$ to $1$ with step $0.05$ and $\mu_2$ equal to $0.6$ and $0.8$ respectively. In this case, the fourth eigenvalue of the FEM model is not matching with the first eigenvalue of the ROM, nor with the second eigenvalue of the ROM, but with the third eigenvalue of the ROM. This is consequence of the fact that the inner products $(u_2(\mu_i),u_4(\mu_j))$ and $(u_3(\mu_i),u_4(\mu_j))$, for $i\neq j$, are not zero nor small. The first eigenvalue of the ROM coincides with the second eigenvalue of the FEM, the second eigenvalue of the ROM coincides with the third eigenvalue of the FEM, and the third eigenvalue of the ROM coincides with the fourth eigenvalue of the FEM. In Figure~\ref{evct_u4_exp} we have shown the eigenvalues of the FEM and the ROM with different number of sample points when we consider all the left singular vectors as a basis. We can see that as we increase the number of sample points the first eigenvalue of the ROM converges to the first eigenvalue of the FEM. Thus all four ROM eigenvalues follow the order of the FEM.

\begin{figure}
     \centering
     \begin{subfigure}{0.4\textwidth}
          \includegraphics[height=5cm,width=5.7cm]{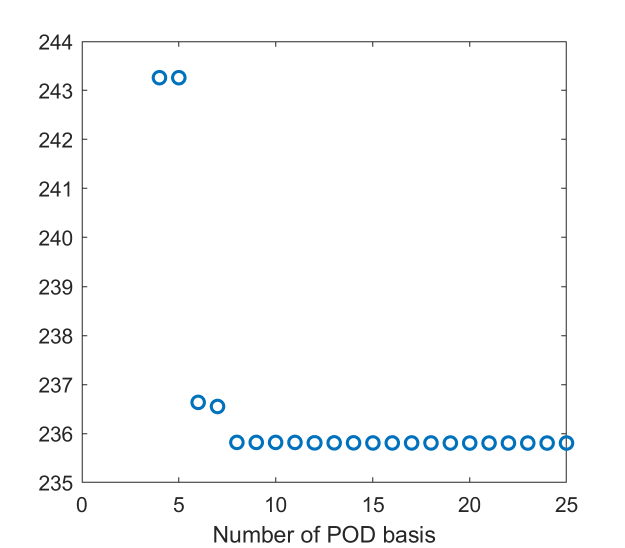}
         \caption{25 sampling points}
     \end{subfigure}
      \begin{subfigure}{0.4\textwidth}
          \includegraphics[height=5cm,width=5.7cm]{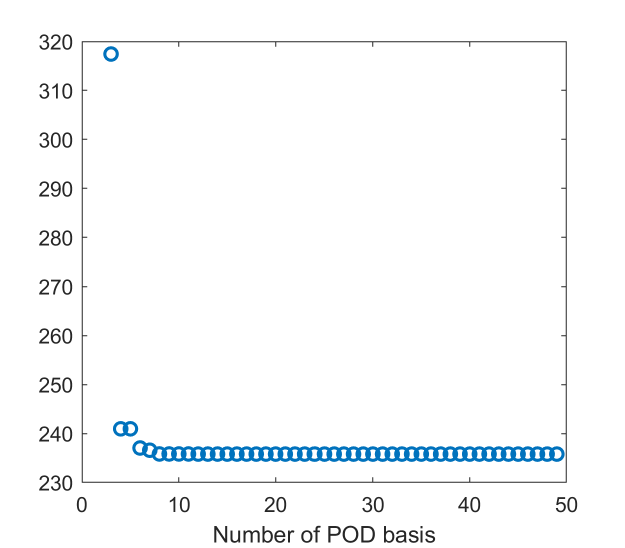}
         \caption{49 sampling points}
     \end{subfigure}
  \caption{Approximation of $\lambda_4$ at $\mu=(0.5,0.6)$ with snapshot based on $u_4$: varying the number of POD basis.}
\label{evpod_u4}
  \end{figure}

The fourth eigenvalue of the FEM and the corresponding ROM at the four test points are presented and the relative errors are presented in Table~\ref{table_u4}. The dimension of the reduced system obtained using the criterion \eqref{criterion} is also mentioned in the same table. The maximum relative error among the four test points is $10^{-6}$. Note that when we increase the number of sample points then the results are improving. We have shown the plot for the ROM at the point $(0.5,0.6)$ with the different ROM dimensions in Figure~\ref{evpod_u4} and one can see that the ROM eigenvalue is converging to the exact eigenvalue when the ROM dimension is more than $4$.

 \begin{table}
 	 	\centering
 	\begin{tabular}{|c|c|c|c|c|c|c|c|} 
 		\hline
 	Sampling&  {\begin{tabular}[c]{@{}c@{}} ROM dim \end{tabular}} &
		  {\begin{tabular}[c]{@{}c@{}} $\pmb{\mu}$ \end{tabular}} &
 		 {\begin{tabular}[c]{@{}c@{}} 4th EV(FEM) \end{tabular}} & 
 		{\begin{tabular}[c]{@{}c@{}} 3rd EV(ROM)  \end{tabular}} &
    {\begin{tabular}[c]{@{}c@{}}  Rel. Error \end{tabular}}\\
  \hline
  25&22
  &(0.5,0.6)&235.80616993&235.80648771 &$1.3 \times 10^{-6}$\\
& &(0.5,0.8)&173.56836538 &173.56839422&$1.6 \times 10^{-7}$	\\
& &(0.8,0.6)&125.56104288 &125.56186825&$6.5\times 10^{-6}$\\
&  &(0.8,0.8)&94.35030619&94.35031462	&$8.9 \times 10^{-8}$ \\ 
\hline
49&28
  &(0.5,0.6)&235.80616993 &235.80617183 &$8.0 \times 10^{-9}$\\
&  &(0.5,0.8)&173.56836538&173.56838071&$8.8 \times 10^{-8}$	\\
&  &(0.8,0.6)&125.56104288 &125.56105296 &$8.0 \times 10^{-8}$\\
&  &(0.8,0.8)&94.35030619&94.35030733 &$1.2 \times 10^{-8}$ \\ 
  \hline
 \end{tabular}
\caption{Approximation of $\lambda_4$ with snapshot based on $u_4$: comparison of FEM and ROM.}
 \label{table_u4}
 \end{table}
 
 \subsubsection{Results of the EVP considering $u_3$ and $u_4$ in the snapshot matrix}
\begin{figure}
     \centering
     \begin{subfigure}{0.4\textwidth}
         \includegraphics[height=5cm,width=5.7cm]{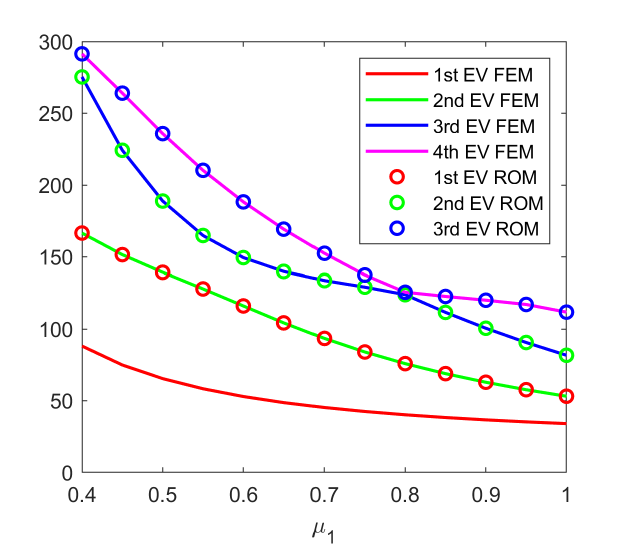}
         \caption{$\mu_2=0.6$}
     \end{subfigure}
      \begin{subfigure}{0.4\textwidth}
          \includegraphics[height=5cm,width=5.7cm]{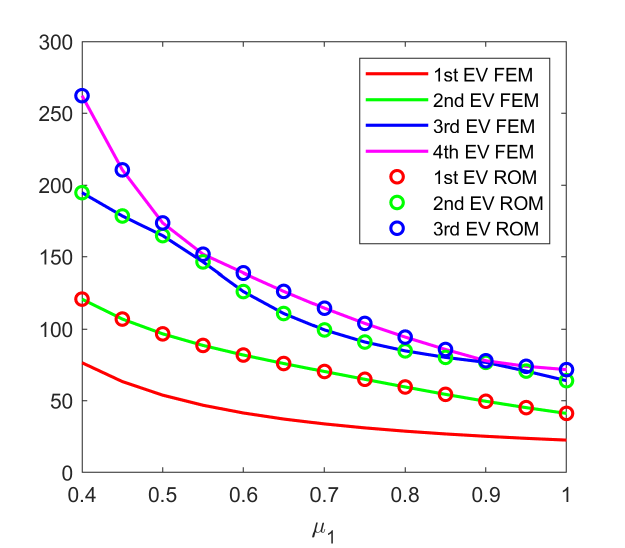}
          \caption{$\mu_2=0.8$}
     \end{subfigure}
  \caption{Approximation of eigenvalues with snapshot based on $u_3$ and $u_4$: comparison of FEM and ROM eigenvalues with varying $\mu_1$ and fixed $\mu_2$, and 49 sample points.}
        \label{evct_u34}
  \end{figure}

\begin{figure}
     \centering
     \begin{subfigure}{0.32\textwidth}
         \includegraphics[height=5cm,width=5.5cm]{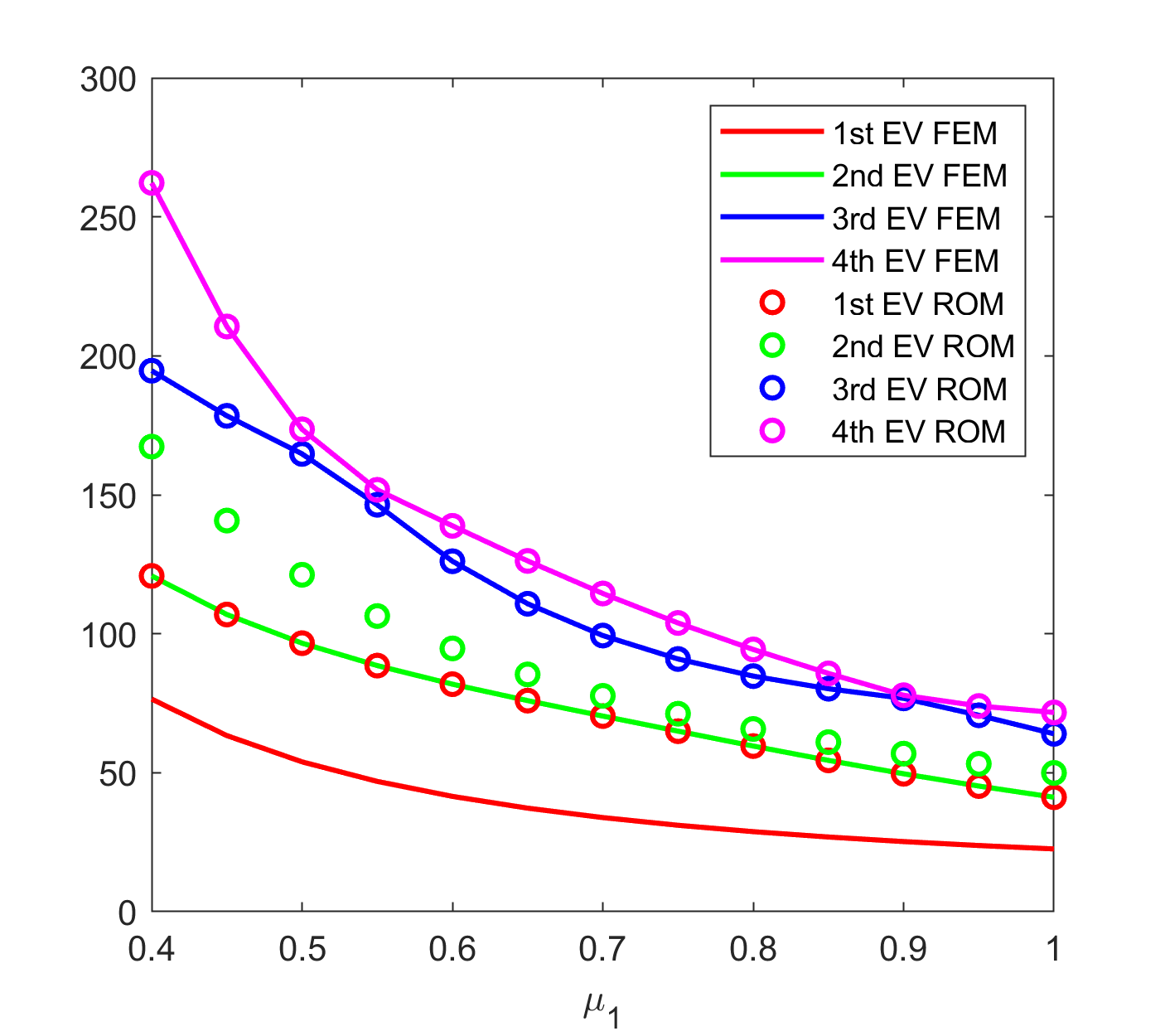}
         \caption{25 samples}
     \end{subfigure}
      \begin{subfigure}{0.32\textwidth}
          \includegraphics[height=5cm,width=5.5cm]{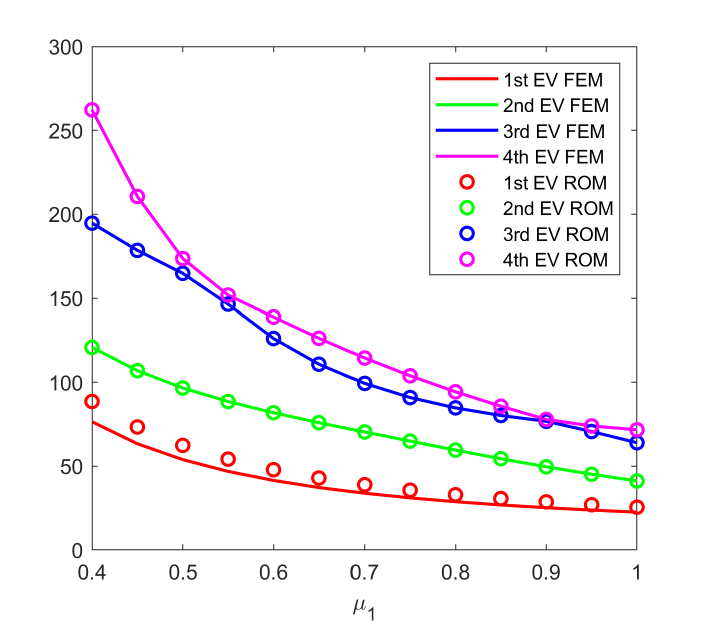}
          \caption{49 samples}
     \end{subfigure}
     \begin{subfigure}{0.32\textwidth}
          \includegraphics[height=5cm,width=5.5cm]{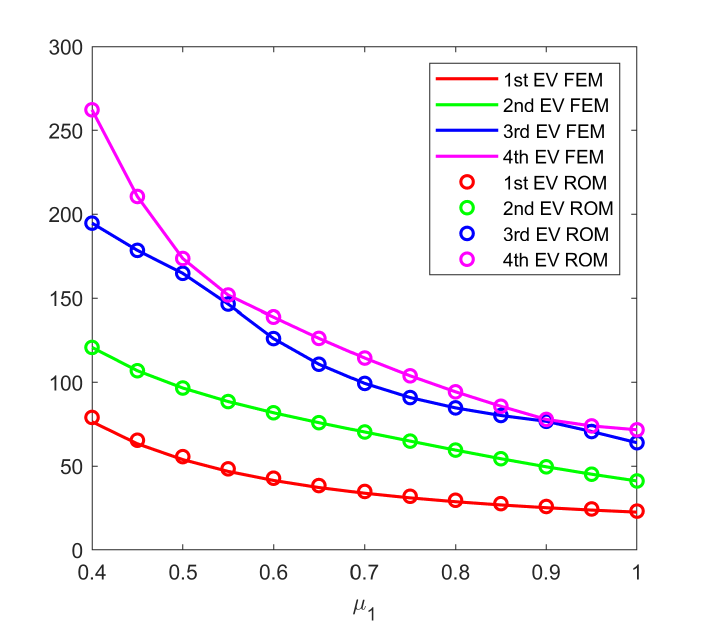}
          \caption{81 samples}
     \end{subfigure}
  \caption{Approximation of eigenvalues with snapshot based on $u_3$ and $u_4$: comparison of FEM and ROM eigenvalues with varying $\mu_1$ and $\mu_2=0.8$, and with different number of sample points.}
        \label{evct_u34_exp}
  \end{figure}

 Since the third and fourth eigenvectors are intersecting, let us put $u_3$ and $u_4$ in the snapshot matrix and observe the result. Now the results behave as in the case when we considered only $u_4$ in the snapshot matrix, that is the fourth eigenvalue of the FEM matches with the third eigenvalue of the ROM and so on, see Figure~\ref{evct_u34}). In Figure~\ref{evct_u34_exp} we have shown the eigenvalues of the FEM and the ROM with different number of sample points and consider all the left singular vectors of the snapshot matrix as a basis. We can see that as we increase the number of sample points the first eigenvalue of the ROM converges to the first eigenvalue of the FEM. Thus all four ROM eigenvalues follow the order of the FEM.
 The fourth eigenvalue of the FEM and the corresponding ROM at the four test points are presented in Table~\ref{table_u34} and the maximum relative error is $10^{-7}$.

 \begin{table}
 	 	\centering
 	\begin{tabular}{|c|c|c|c|c|c|c|c|} 
 		\hline
 	Sampling&{\begin{tabular}[c]{@{}c@{}} ROM dim \end{tabular}} &
		  {\begin{tabular}[c]{@{}c@{}} $\mu$ \end{tabular}} &
 		 {\begin{tabular}[c]{@{}c@{}} 4th EV(FEM) \end{tabular}} & 
 		{\begin{tabular}[c]{@{}c@{}} 3rd EV(ROM)  \end{tabular}}&
    {\begin{tabular}[c]{@{}c@{}}  Rel. Error \end{tabular}}\\
  \hline
  25 &30
  &(0.5,0.6)&235.80616993 &235.80618623&$6.9 \times 10^{-8}$  \\
 & &(0.5,0.8)&173.56836538  &  173.56838092 &$8.9 \times 10^{-8}$	\\
 & &(0.8,0.6)&125.56104288  & 125.56105844 &$1.2 \times 10^{-7}$	\\
 & &(0.8,0.8)&94.35030619 & 94.35030784&$1.7 \times 10^{-8}$\\ 
  \hline
 49 &31
  &(0.5,0.6)&235.80616993   &235.80617885 &$3.7 \times 10^{-8}$\\
 & &(0.5,0.8)&173.56836538   & 173.56838605 &$1.1 \times 10^{-7}$	\\
 & &(0.8,0.6)&125.56104288  &125.56105806 &$1.2 \times 10^{-7}$	\\
 & &(0.8,0.8)&94.35030619 &94.35030780&$1.7 \times 10^{-8}$\\ 
  \hline
 \end{tabular}
 \caption{Approximation of $\lambda_4$ with snapshot based on $u_3,u_4$: comparison of FEM and ROM.}
\label{table_u34}
 \end{table}

  \begin{figure}
     \centering
     \begin{subfigure}{0.4\textwidth}
          \includegraphics[height=5cm,width=5.7cm]{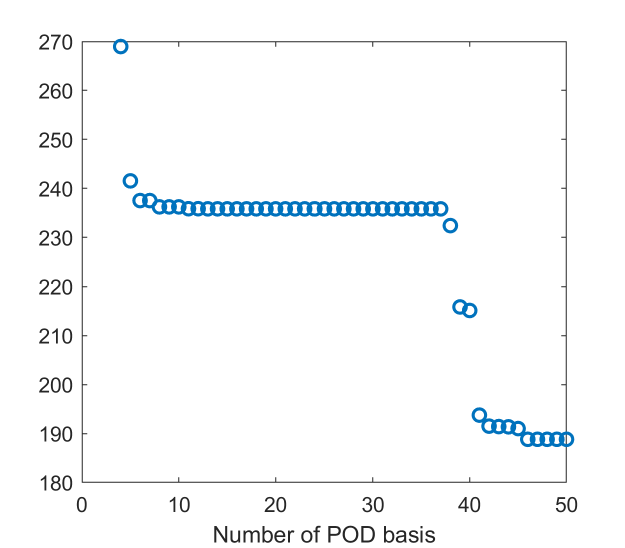}
         \caption{25 sampling points}
     \end{subfigure}
      \begin{subfigure}{0.4\textwidth}
          \includegraphics[height=5cm,width=5.7cm]{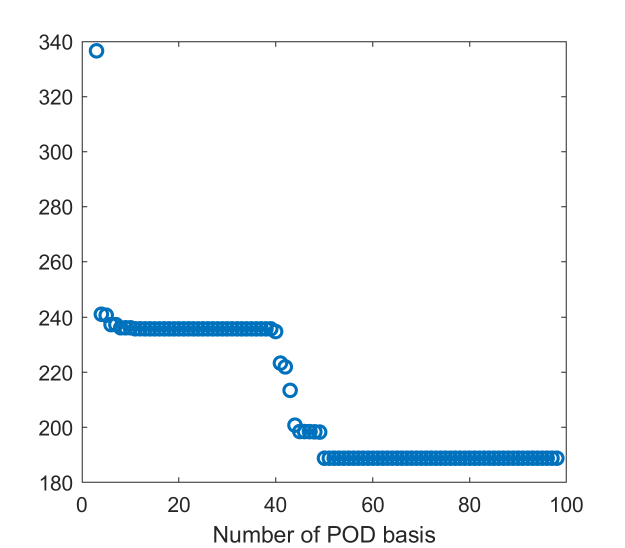}
         \caption{49 sampling points}
     \end{subfigure}
   \caption{Approximation of $\lambda_4$ at $\mu=(0.5,0.6)$ with snapshot based on $u_3,u_4$: varying the number of POD basis.}
        \label{evpod_u34}
  \end{figure}
 
\subsubsection{Results of the EVP considering $u_1$, $u_2$, $u_3$, and $u_4$ in the snapshot matrix}
\begin{figure}
     \centering
     \begin{subfigure}{0.4\textwidth}
         \includegraphics[height=5cm,width=5.7cm]{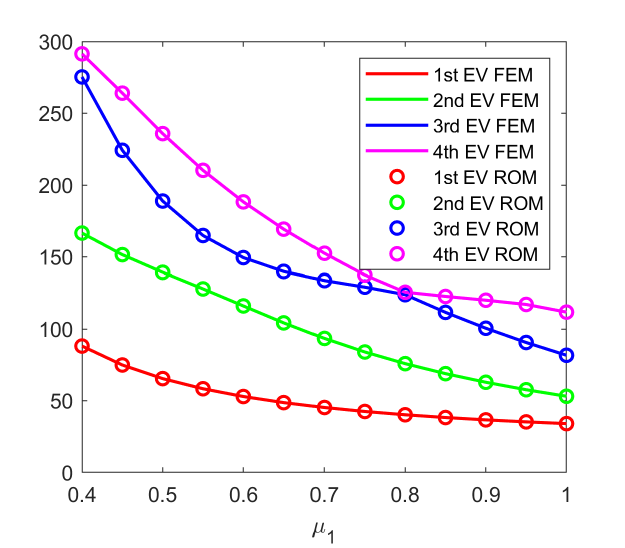}
         \caption{$\mu_2=0.6$}
     \end{subfigure}
      \begin{subfigure}{0.4\textwidth}
          \includegraphics[height=5cm,width=5.7cm]{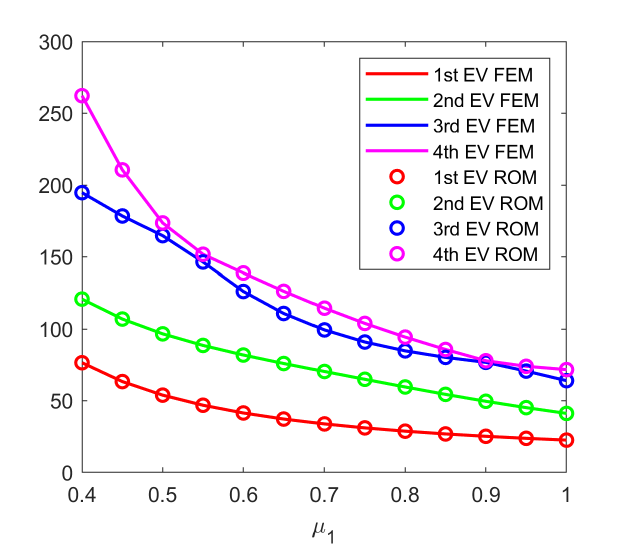}
          \caption{$\mu_2=0.8$}
     \end{subfigure}
  \caption{Approximation of eigenvalues with snapshot based on $u_1$, $u_2$, $u_3$, and $u_4$: comparison of FEM and ROM eigenvalues with varying $\mu_1$ and fixed $\mu_2$, and 49 sample points.}
        \label{evct_u14}
  \end{figure}

\begin{figure}
     \centering
     \begin{subfigure}{0.32\textwidth}
         \includegraphics[height=5cm,width=5.5cm]{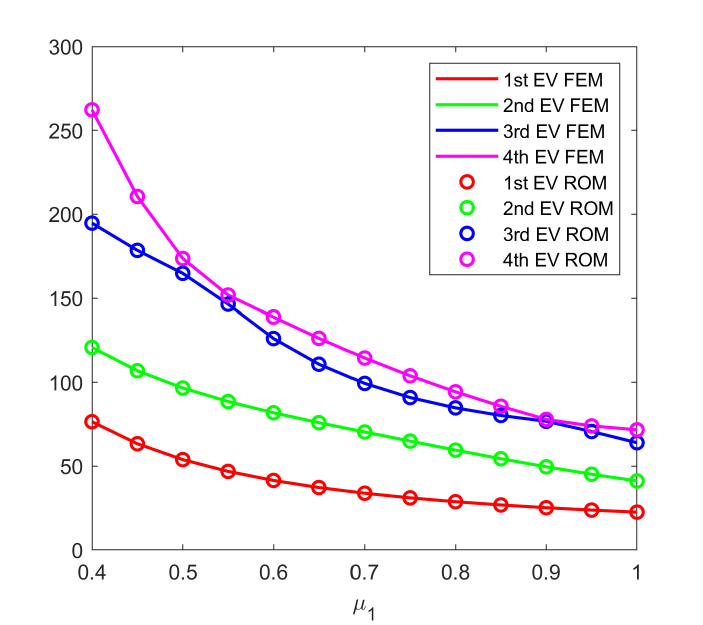}
         \caption{25 samples}
     \end{subfigure}
      \begin{subfigure}{0.32\textwidth}
          \includegraphics[height=5cm,width=5.5cm]{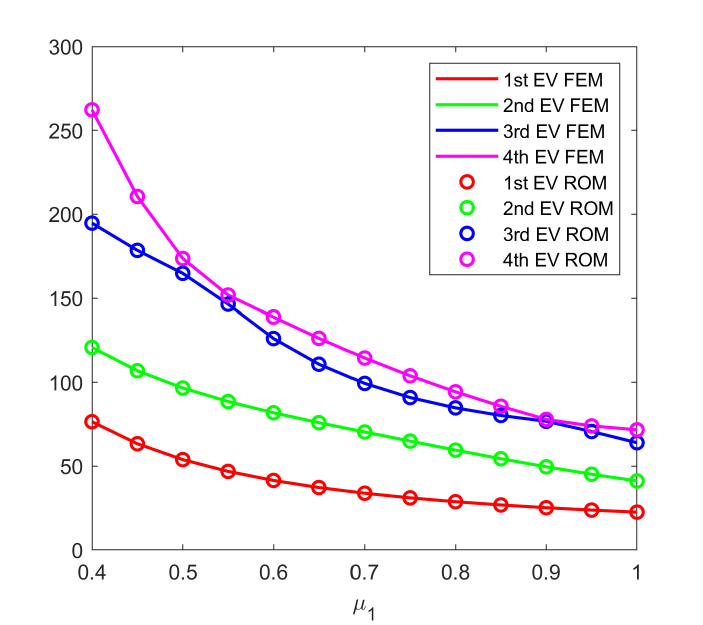}
          \caption{49 samples}
     \end{subfigure}
     \begin{subfigure}{0.32\textwidth}
          \includegraphics[height=5cm,width=5.5cm]{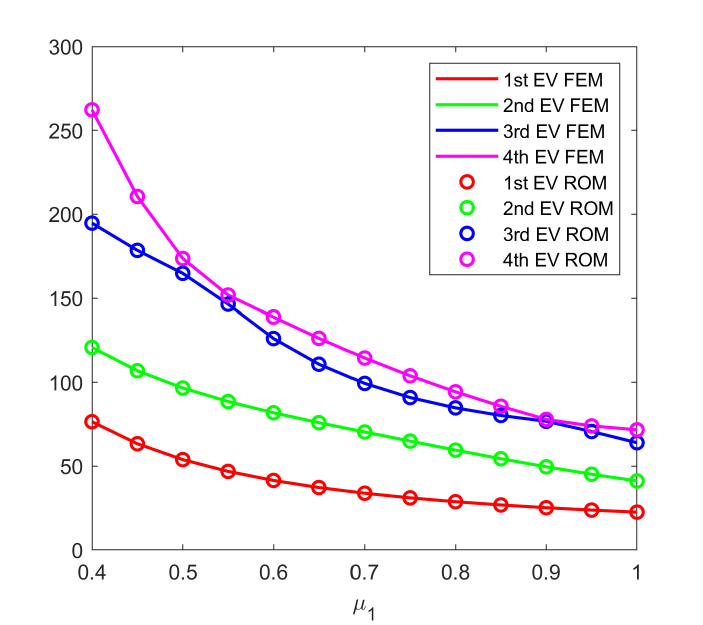}
          \caption{81 samples}
     \end{subfigure}
  \caption{Approximation of eigenvalues with snapshot based on $u_1$, $u_2$, $u_3$, and $u_4$: comparison of FEM and ROM eigenvalues with varying $\mu_1$ and $\mu_2=0.8$ for different sample points.}
        \label{evct_u14_exp}
  \end{figure}

Then we consider the snapshot matrix containing the first fourth eigenvectors at the sample points. The first four eigenvalues of the FEM and ROM are plotted in Figure~\ref{evct_u14} for the parameters $\mu_1$ ranging from $0.4$ to $1$ with step $0.05$ and $\mu_2$ equal to $0.6$ and $ 0.8$, respectively. In this case, all the first four eigenvalues of the ROM match with the first four eigenvalues of the FEM. In this case, considering all the left singular vectors of the snapshot matrix as a basis, the eigenvalues are stable (see Figure~\ref{evct_u14_exp}).

\begin{figure}
     \centering
     \begin{subfigure}{0.4\textwidth}
          \includegraphics[height=5cm,width=5.7cm]{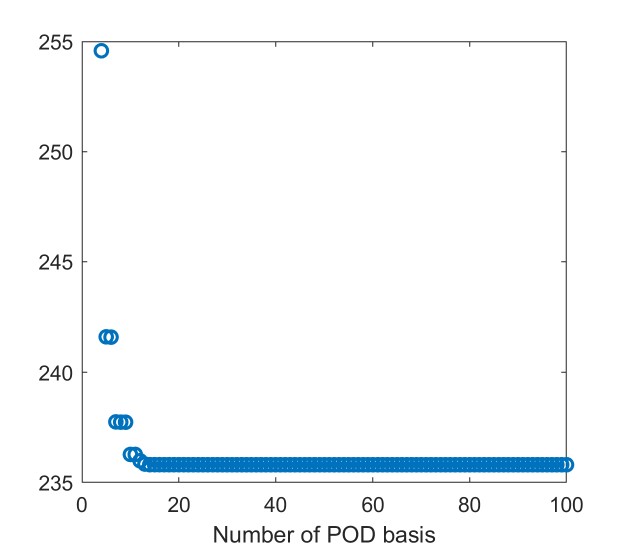}
         \caption{25 sampling points}
     \end{subfigure}
      \begin{subfigure}{0.4\textwidth}
          \includegraphics[height=5cm,width=5.7cm]{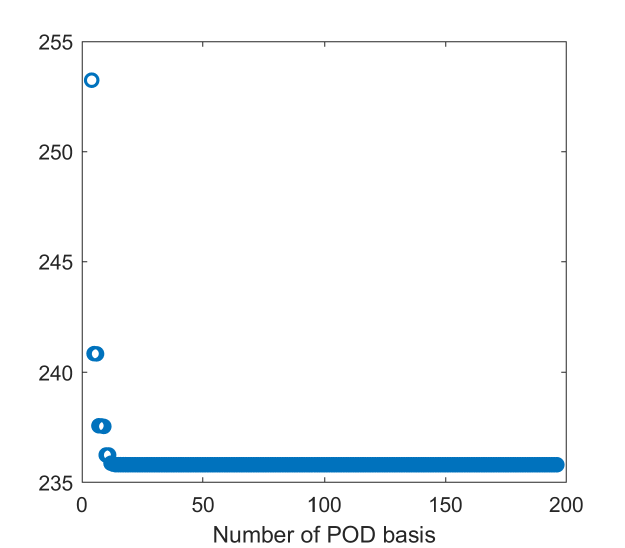}
         \caption{49 sampling points}
     \end{subfigure}
  \caption{Approximation of $\lambda_4$ at $\mu=(0.5,0.6)$ with snapshot based on $u_1,u_2,u_3,u_4$: varying the number of POD basis.}
        \label{evpod_u14}
  \end{figure}

The fourth eigenvalue of the FEM and the fourth eigenvalue of the ROM at the four test points and their relative errors are presented in Table~\ref{table_u14}. The dimension of the reduced system obtained using the criterion~\eqref{criterion} is also mentioned in the same table. The maximum relative error among the four test points is $10^{-7}$. We have shown the plot for the fourth eigenvalue of the ROM at the point $(0.5,0.6)$ with the different ROM dimensions in Figure~\ref{evpod_u14} and one can see that the ROM eigenvalue is converging to the exact eigenvalue when the ROM dimension is more than $15$.

 \begin{table}
 	 	\centering
 	\begin{tabular}{|c|c|c|c|c|c|c|c|} 
 		\hline
         Sampling& {\begin{tabular}[c]{@{}c@{}} ROM dim \end{tabular}} &
		  {\begin{tabular}[c]{@{}c@{}} $\mu$ \end{tabular}} &
 		{\begin{tabular}[c]{@{}c@{}} 4th EV(FEM) \end{tabular}} & 
 		{\begin{tabular}[c]{@{}c@{}} 4th EV(ROM)  \end{tabular}}&
    {\begin{tabular}[c]{@{}c@{}}  Rel. Error($\lambda_4$) \end{tabular}}\\
  \hline
 25 &39
  &(0.5,0.6)&235.80616993  &235.80620424&$1.4 \times 10^{-7}$ \\
&  &(0.5,0.8)&173.56836538& 173.56838242 &$9.8 \times 10^{-8}$ 	\\
&  &(0.8,0.6)&125.56104288 & 125.56109105&$3.8 \times 10^{-7}$  \\
&  &(0.8,0.8)&94.35030619&  94.35032541 &$2.0 \times 10^{-7}$ \\ 
  \hline
 49 &40
  &(0.5,0.6)&235.80616993 &235.80619246&$9.5 \times 10^{-8}$ \\
&  &(0.5,0.8)&173.56836538& 173.56837272 &$4.2 \times 10^{-8}$	\\
&  &(0.8,0.6)&125.56104288 &125.56107025&$2.1 \times 10^{-7}$ \\
&  &(0.8,0.8)&94.35030619& 94.35032558&$2.0 \times 10^{-7}$ \\ 
  \hline
 \end{tabular}
\caption{Approximation of $\lambda_4$ with snapshot based on $u_1,u_2,u_3,u_4$: comparison of FEM and ROM.}
\label{table_u14}
 \end{table}

\subsubsection{Results of the EVP considering $u_1+u_2+u_3+u_4$ in the snapshot matrix}
\begin{figure}
     \centering
     \begin{subfigure}{0.4\textwidth}
         \includegraphics[height=5cm,width=5.7cm]{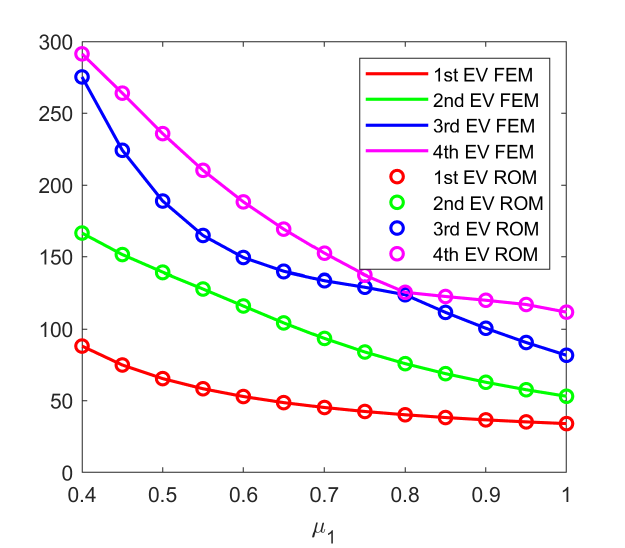}
         \caption{$\mu_2=0.6$}
     \end{subfigure}
      \begin{subfigure}{0.4\textwidth}
          \includegraphics[height=5cm,width=5.7cm]{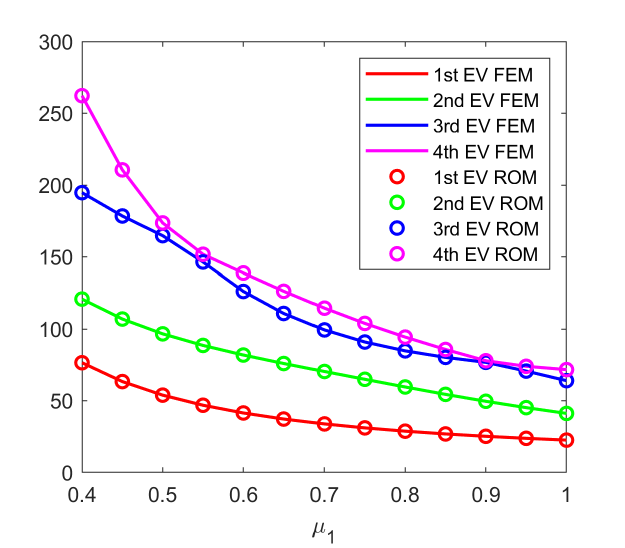}
          \caption{$\mu_2=0.8$}
     \end{subfigure}
  \caption{Approximation of $\lambda_1,\lambda_2,\lambda_3,\lambda_4$ with snapshot based on $u_1+u_2+u_3+u_4$: comparison of FEM and ROM eigenvalues with varying $\mu_1$ and fixed $\mu_2$, and 49 sample points.}
        \label{evct_au4}
  \end{figure}

The results corresponding to the case when we consider all the first four eigenvectors in the snapshot matrix are good and the order of the eigenvalues of the ROM and the FEM are the same, but the number of snapshots is four times than the sample points. In order to reduce the computational cost, we add the first four eigenvectors and choose the resulting vector as a snapshot to control the number of snapshots. Also in this case all the ROM eigenvalues match the corresponding FEM eigenvalues, as it is shown in Figure~\ref{evct_au4}.

The fourth eigenvalue of the FOM and the fourth eigenvalue of the ROM at the four test points and their relative errors are presented in Table~\ref{table_au4}. The dimension of the reduced system obtained using the criterion~\eqref{criterion} is also mentioned in the same table. The maximum relative error among the four test points is $10^{-5}$. We have shown the plot for the fourth eigenvalue of ROM at the point $(0.5,0.6)$ with the different ROM dimensions in Figure~\ref{evpod_u14} and one can see that the ROM eigenvalue is converging to the exact eigenvalue when the ROM dimension is more than $10$. We use the sum of four eigenvectors we preserve the order of the eigenvalues of ROM and FEM at a price of getting a relative error which is higher than in the case where we use all the eigenvectors.
 
   \begin{table}
 	 	\centering
 	\begin{tabular}{|c|c|c|c|c|c|c|c|} 
 		\hline
        Sampling& {\begin{tabular}[c]{@{}c@{}}ROM dim\end{tabular}} &
		  {\begin{tabular}[c]{@{}c@{}} $\mu$ \end{tabular}} &
 		{\begin{tabular}[c]{@{}c@{}} 4th EV(FEM) \end{tabular}} & 
 		{\begin{tabular}[c]{@{}c@{}} 4th EV(ROM)  \end{tabular}}&
    {\begin{tabular}[c]{@{}c@{}}  Rel. Error($\lambda_4$) \end{tabular}}\\
  \hline
 25 &24
  &(0.5,0.6)&235.80616993  &235.80707671&$3.8 \times 10^{-6}$ \\
 & &(0.5,0.8)&173.56836538&173.57125574 &$1.6 \times 10^{-5}$ 	\\
&  &(0.8,0.6)&125.56104288 &125.56717711&$4.8 \times 10^{-5}$  \\
&  &(0.8,0.8)&94.35030619&94.35107680  &$8.1 \times 10^{-6}$ \\ 
  \hline
 49 &40
  &(0.5,0.6)&235.80616993 &235.80631238& $6.0 \times 10^{-7}$ \\
&  &(0.5,0.8)&173.56836538&173.56860618 & $1.3 \times 10^{-6}$	\\
 & &(0.8,0.6)&125.56104288 &125.56118125&$1.1 \times 10^{-6}$ \\
 & &(0.8,0.8)&94.35030619& 94.35036064& $5.7 \times 10^{-7}$ \\ 
  \hline
\end{tabular}
\caption{Approximation of $\lambda_4$ with snapshot based on $u_1+u_2+u_3+u_4$: comparison of FEM and ROM.} 
\label{table_au4}
 \end{table} 
 
 \begin{figure}
 \centering
\begin{subfigure}{0.4\textwidth}
\includegraphics[height=5cm,width=5.5cm]{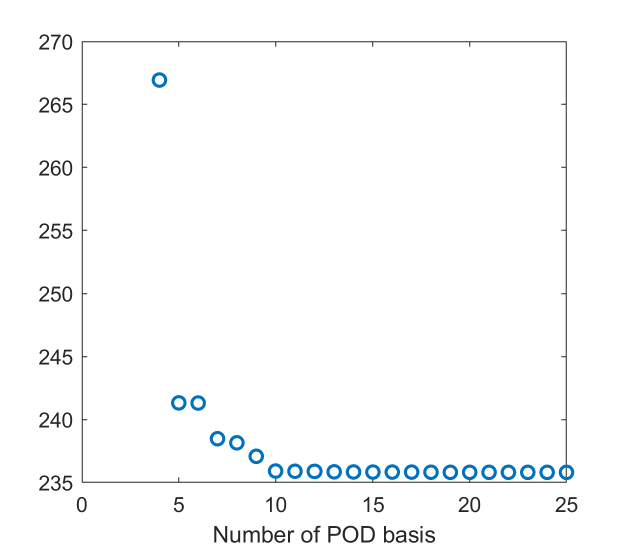}
\caption{25 sampling points}
\end{subfigure}
\begin{subfigure}{0.4\textwidth}
\includegraphics[height=5cm,width=5.5cm]{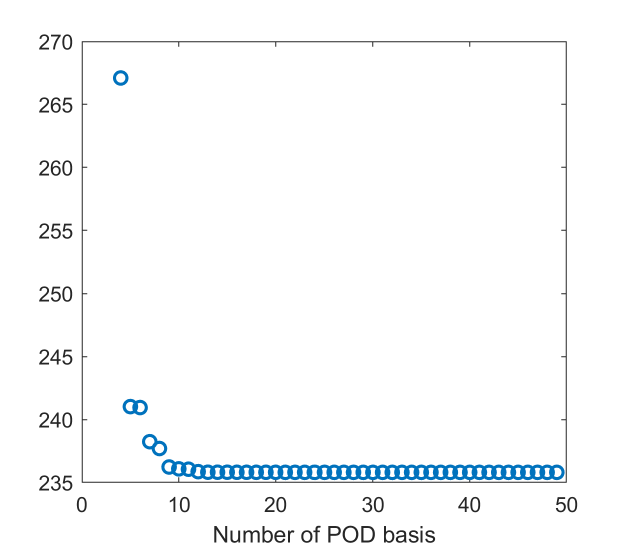}
\caption{49 sampling points}
\end{subfigure}
\caption{Approximation of $\lambda_4$ at $\mu=(0.5,0.6)$ with snapshot based on $u_1+u_2+u_3+u_4$: varying the number of POD basis.}
\label{evpod_au4}
\end{figure}

\section*{Acknowledgments}
This research was supported by the Competitive Research Grants Program CRG2020 ``Synthetic data-driven model reduction methods for modal analysis'' and by the Opportunity Fund Program OFP2023 ``Model reduction methods for modal analysis in computational mechanics'' awarded by the King Abdullah University of Science and Technology (KAUST).
Daniele Boﬃ is member of the INdAM Research group GNCS.

\section*{Conflict of interest}

Authors have no conflict of interest to declare.

\bibliographystyle{plain}
\bibliography{mybib}
\end{document}